\documentclass{acta_2004}
\usepackage{harvard_acta}
\usepackage{amssymb}
\usepackage{graphicx}
\usepackage{amsmath}
\usepackage{pdfsync}
\pagerange{\pageref{firstpage}--\pageref{lastpage}}
\usepackage{algorithmicx,algorithm,algpseudocode}

\newcounter{myenumii}

\usepackage{amsfonts}
\usepackage{amstext}
\usepackage{amsopn}
\usepackage{amsxtra}
\usepackage{mathrsfs}

\newcommand{\ignore}[1]{}

\newcommand{\rd}{\mathrm{d}}

\DeclareMathOperator{\im}{Im}

\newcommand{\ex}{\mathrm{e}}

\def\half{\frac 1 2}
\newcommand{\p}{\partial}
\newcommand{\fl}[2]{\frac{#1}{#2}}
\newcommand{\btu}{\Delta}

\newcommand{\ift}{\infty} 

\newcommand{\bR}{{\mathbb R }}

\newcommand{\bN}{{\mathbb N }}

\newcommand{\bC}{{\mathbb C}}
\newcommand{\e}{{\varepsilon }}

\newcommand{\eps}{ {\varepsilon} }

\newcommand{\E}{\mathbb{E}}
\newcommand{\bbP}{\mathbb{P}}

\newcommand{\be}{\begin{equation}}
\newcommand{\ee}{\end{equation}}
\newcommand{\ba}{\begin{array}}
\newcommand{\ea}{\end{array}}
\newcommand{\bea}{\begin{eqnarray}}
\newcommand{\eea}{\end{eqnarray}}
\newcommand{\beas}{\begin{eqnarray*}}
\newcommand{\eeas}{\end{eqnarray*}}

\newcommand{\FF}{\mathcal{F}}

\newcommand{\LL}{\mathcal{L}}
\newcommand{\EE}{\mathcal{E}}

\def\Dt{\partial_t}

\def\({\left(}
\def\){\right)}
\def\<{\left\langle}
\def\>{\right\rangle}

\newcommand{\Id}[1]{{\rm I\kern-2pt I_{#1}}}
\renewcommand{\hbar}{{\displaystyle\bar{\phantom{x}}\kern-6pt h}}

\numberwithin{equation}{section}

\newcommand{\vep}{\varepsilon}

\newtheorem{theorem}{Theorem}[section]

\newtheorem{remark}[theorem]{Remark}
\newtheorem{example}[theorem]{Example}
\newtheorem{assumption}{Assumption}

\newcommand\bx{\mathbf{x}}

\newcommand\ds{\displaystyle}
\newcommand{\RR}{\mathbb{R}}
\newcommand\dd{\textrm{d}}
\begin{document}

\label{firstpage}

\title[Asymptotic-Preserving Schemes for Multiscale  Problems]{Asymptotic-Preserving Schemes for Multiscale Physical Problems}

\author[Shi Jin]
{Shi Jin\\
School of  Mathematical Sciences, Institute of Natural Sciences, and MOE-LSC\\
 Shanghai Jiao Tong University,
Shanghai 200240, China.\\ 
Email: {shijin-m@sjtu.edu.cn}
}

\maketitle

%\vspace*{4cm}

\begin{abstract}
We present the asymptotic transitions from microscopic to macroscopic physics, their computational  challenges and the Asymptotic-Preserving (AP) strategies to efficiently compute multiscale physical problems. Specifically, we will first
study the asymptotic  transition from quantum to classical mechanics, from classical mechanics to kinetic theory, and then from kinetic theory to hydrodynamics. We then review some representative
AP schemes  that mimic, at the discrete level, these asymptotic transitions, hence can be used crossing scales and, in particular,   capture the macroscopic behavior without resolving numerically the microscopic physical scale.  
\end{abstract}

\footnote{S. Jin was partially supported by National Key R\&D Program of China Project No. 2020YFA0712000, NSFC grant No. 12031013, and Shanghai Municipal Science and Technology Major Project 2021SHZDZX0102.}

%\newpage 

\tableofcontents

\section{Introduction}\label{sec1}
\setcounter{equation}{0}
%This paper is devoted to the construction of numerical schemes making the  
%connection between kinetic theory and macroscopic fluid  dynamics. 
%For instance, 

Ignoring relativistic effects, quantum mechanics
is considered to be enough for one to understand the physical properties of matters. Since solving the Schr\"odinger equation analytically is impossible, one relies upon computer simulations to solve the equation. However, there are essential computational bottlenecks in  simulation at the quantum level. First is the curse of dimensionality.  For common molecules like carbon dioxide CO$_2$, which consists of 3 nuclei and 22 electrons, the full time-dependent
Schr\"odinger equation is defined in $75$ space dimensions! 
The benzene molecule consists of $12$ nuclei and $42$ electrons, hence one needs to solve  the Schr\"odinger equation in $162$ dimensions.
Another challenge is that quantum mechanics is valid at spatial scales of Angstrom, which is $10^{-10}$ meter, and time scale of femtoseconds, or $10^{-15}$ second. To simulate such a small
scale system to any physical scales of interest, for example,
micrometers to millimeters, or microseconds to milliseconds, is 
computationally formidable by today's computer.

\begin{figure}[ht]
\centering
%\subfigure[]{
  \includegraphics[width=0.6\textwidth]{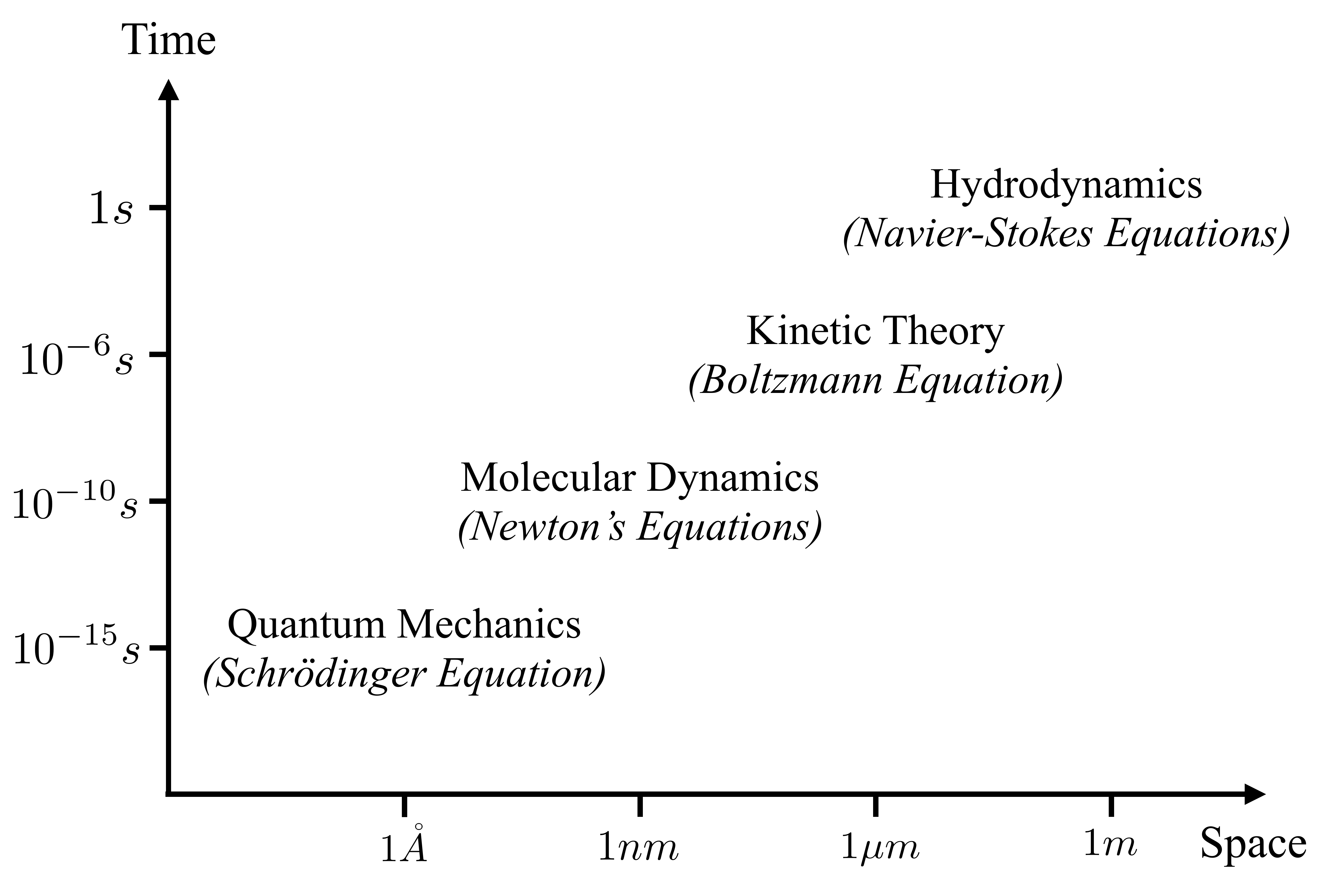}
%}
\caption{Multiscale diagram.}                                                        
%$\rho_{\rm in}(x)=|a_{\rm in}|^2(x)=e^{-(x-1/2)^2}\,$,                                             
\label{fig-2}
\end{figure}

Physical models at larger scales, such as classical mechanics, statistical mechanics, and hydrodynamics, are computationally much less expensive compared to a quantum simulation, but they are valid in certain time and spatial scales, see  Fig. \ref{fig-2}. When one deals with problems that go across different scales, either due to the nature of the problems or computational needs, multiscale computation becomes a viable tool when one cannot afford to resolve the smallest scales.

Understanding the transition from one scale to another is a central topic in mathematical physics and partial differential equations. They are  related to Hilbert's sixth problem \cite{corry}. These asymptotic transitions are  not only of great mathematical interest,
they also guide the design of multiscale computational methods, as will be reflected in this survey.

If the asymptotic or macroscopic equations are uniformly valid in the
entire domain of interest, it is much more efficient just to solve the problem  at the macroscopic level, which are computationally more economical. 
However, there are many problems where the macroscopic models break down
in part of the domain, or one lacks information or data on the macro models, thus the microscopic models
are needed, at least locally. 
Therefore a multiscale and multiphysics approach, that 
hybridizes the microscopic and the macroscopic models in a domain-decomposition or heterogeneous multiscale 
framework, becomes necessary, see for examples
\cite{BPT} \cite{KNS} \cite{DJ05} for multiscale kinetic
problems, and \cite{EE-CMS} \cite{E-Acta} \cite{Kev} for  broader areas of
multiscale modeling and simulation.

Central to the design of multiscale computational methods  is to identify the critical physical scales in the system and the connections between microscopic and macroscopic models. The Schr\"odinger equation is valid at the scale of Angstrom,
which is exceedingly small compared to the scale of interest.
The Newton equations in classical mechanics often involve the number of particles simply formidably too large. 
Kinetic equations often contain small mean free path or time, or Knudsen number, the average
distance or time between two collisions of particles.  When the characteristic scales become small,  tremendous computational challenges arise since one needs to
numerically resolve these small scales which can be prohibitively
expensive. A main difficulty in most multiscale and multi-physics 
 type methods is that one has to couple models at different
scales through an interface or buffer zone where one has to match
two different models. While it is often easy to generate macroscopic data from the micro ones through, for examples, ensemble averages or taking moments, it is difficult to convert macroscopic data to the microscopic ones, since most of the time this conversion is not unique. 
 The coupling locations may also be difficult to
determine in a dynamic problem.

On the other hand, asymptotic expansions on these 
small parameters for a microscopic model usually give rise to the   macroscopic equations. One hopes such a transition can guide the design, 
and help to analyze, effective and efficient multiscale computational methods. 
 
This paper surveys one multiscale framework--the {\it asymptotic-preserving (AP) schemes}. This approach has its origin in capturing steady-state
solution for neutron transport in the diffusive regime \cite{LMM1} \cite{LMM2}.
Since the 90s of last century, the AP schemes 
have been developed for a wide range of time-dependent kinetic and 
hyperbolic 
equations, and far beyond.  The basic idea is to develop
numerical methods that {\it preserve the asymptotic limits from the microscopic
to the macroscopic models, in the discrete setting}. 
Comparing with  multi-physics domain decomposition type methods,
the AP schemes solve one set of equations--the microscopic ones--thus avoid the
coupling of different models. An AP scheme switches from a microscopic solver to the macroscopic solver {\it automatically}.
Specifically, if one numerically resolves the small physical scales then the scheme is a micro solver. Otherwise it effectively becomes
a macro solver, when the physically small scales are {\it not numerically resolved}.

\begin{figure}[ht]
\centering
%\subfigure[]{
  \includegraphics[width=.85\textwidth]{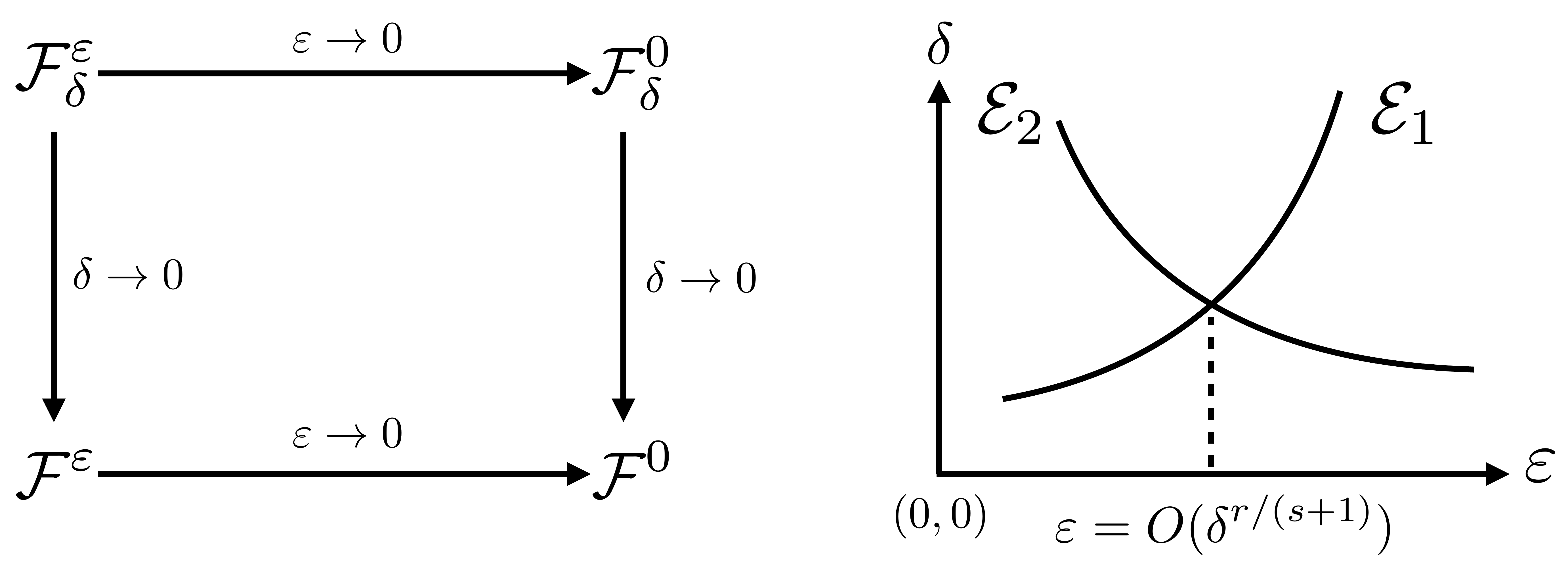}
%}
\hspace{8pt}
%\subfigure[Second order]{
 % \includegraphics[width=0.4\textwidth]{fig-2.pdf}
%}

\caption{Left: Illustration of AP schemes; right: illustration of uniform convergence of AP schemes.}                                                        
%$\rho_{\rm in}(x)=|a_{\rm in}|^2(x)=e^{-(x-1/2)^2}\,$,                                             
\label{fig-1}
\end{figure} 
The idea of AP can be illustrated by the left figure in  Fig. \ref{fig-1}. Assume one starts
with a microscopic model $\FF^\vep$ which depends on the scaling parameter $\vep$.  As $\vep \to 0$ the model asymptotically approaches the  macroscopic model $\FF^0$ which is independent of $\vep$. 
 Denote the numerical discretization of $\FF^\vep$  by $\FF^\vep_\delta$, 
where $\delta$ is the 
numerical parameter (such as mesh size and/or time step). 
The asymptotic limit of $\FF^\vep_\delta$, as $\vep \to 0$ 
(with $\delta$ fixed), if exists, is
denoted by $\FF^0_\delta$.
{\it  Scheme $\FF^\vep_\delta$ is called
AP if $\FF^0_\delta$ is a good (consistent and
stable) approximation of $\FF^0$, }

Error estimates on an AP scheme can be obtained by the following argument \cite{GJL99}.
Typically, 
\begin{equation}
  \| \FF^\vep - \FF^0 \| = O(\vep)\,,
\label{error1}
\end{equation}
under some suitable, problem dependent, norm. 
Assume $\FF^\vep_\delta$ is an $r$-th order approximation to $\FF^\vep$ for fixed $\vep$.
Due to the presence of the small parameter $\vep$, a classical numerical
analysis typically gives the following error estimate
\begin{equation}
  \EE_1 =  \| \FF^\vep_\delta - \FF^\vep\| = O( \delta^r /\vep^s)\,,
\quad 1\le s \le r\,.
\label{error2}
\end{equation}
The error is large when $\delta^r \gg \vep^s$,  namely, when the small physical scales are not numerically resolved (if one uses coarse meshes or large time steps relative to $\vep$). 
An AP scheme usually requires
\begin{equation}
  \| \FF^\vep_\delta - \FF^0_\delta\| = O(\vep)
\qquad {\text {uniformly in}} \quad \delta\,,
\label{error3}
\end{equation}
{\it and}
\begin{equation}
  \| \FF^0_\delta - \FF^0\| = O(\delta^r) \,.
\label{error4}
\end{equation}
   Clearly, if one
adds up the errors in (\ref{error1}), (\ref{error3}) and (\ref{error4}), by the triangle
inequality,  the
following error estimate can be obtained
\begin{equation}
  \EE_2 =  \| \FF^\vep_\delta - \FF^\vep\| \le \| \FF^\vep_\delta - \FF^0_\delta\|+\| \FF^0_\delta- \FF^0\|+\| \FF^0 - \FF^\vep\|= O(\vep+\delta^r)\,.
\label{error5}
\end{equation}
This error is small for $\vep\ll 1$. Clearly both estimates on $\EE_1$ (\ref{error2}) and $\EE_2$ (\ref{error5}) are mathematically valid and can hold simultaneously. By comparing the two  estimates, 
\begin{equation}
  \| \FF^\vep_\delta - \FF^\vep\| = \min(\EE_1, \EE_2)\,, 
\nonumber
\end{equation}
which has an upper bound around $\vep=O(\delta^{r/(s+1)})$, as shown
by the right figure in Fig.\ref{fig-1}. This gives
\begin{equation}
  \| \FF^\vep_\delta - \FF^\vep\| = O(\delta^{r/(s+1)})\,, \qquad {\text {uniformly in}} \quad \vep \,.
\label{error6}
\end{equation}
This argument shows that an AP scheme is {\it convergent uniformly in
$\vep$}.  Indeed, if one resolves $\vep$ by $\delta$ (with $\delta
=o(\vep^{s/r})$) one gets a good
approximation to the microscopic model $\FF^\vep$, as shown by
(\ref{error2}).  If $\vep$ is not resolved
by $\delta$ then one obtains a good approximation to the macroscopic
model $\FF^0$. This transition is done {\it automatically} by the code. 

There have been a few earlier reviews on AP schemes, for examples for multiscale kinetic equations
\cite{jin2010asymptotic}\cite{DP-Acta}\cite{degond-JCPReview}\cite{HJL-AP}, and for semiclassical computation of the Schr\"odinger equation \cite{JMS-Acta}\cite{LaLu20}. This survey, however, is unique in that it covers
in a more comprehensive way the topics in essentially {\it all} important physical regimes, from quantum to classical mechanics, from classical mechanics to kinetic theory, and then from kinetic theory to hydrodynamics. It has also included the most recent
advances, including new directions,  in this topic. 

%Since As  We refer to AP schemes for kinetic equations in the fluid dynamic or
%diffusive regimes \cite{LMM1, LMM2, CP, CJR, JPT1, JPT2, JP, K, K2, GT,
%BLM, LM}. The AP framework has also been extended in
%\cite{CDV1, CDV2,DDN} for the study of the quasi-neutral limit of Euler-Poisson 
%and Vlasov-Poisson systems, and in \cite{DJL, DT, HJL} for 
%all-speed (Mach number) fluid equations bridging the passage from 
%compressible flows to the incompressible flows. 
%One should note that under-resolved computation may not yield accurate or
%even physically correct approximations in areas with sharp transitions,
%such as shock and boundary layers. In these areas one may want to use
%resolved calculations. The AP schemes allow one to use suitable mesh size
%and time step at needed domains with one first-principle equation, thus is 
%especially
%suitable for problems with localized sharp transitions where macroscopic
%simulation is necessary.

%Figure of AP scheme here.

\begin{figure}[ht]
\centering
%\subfigure[]{
  \includegraphics[width=0.6\textwidth]{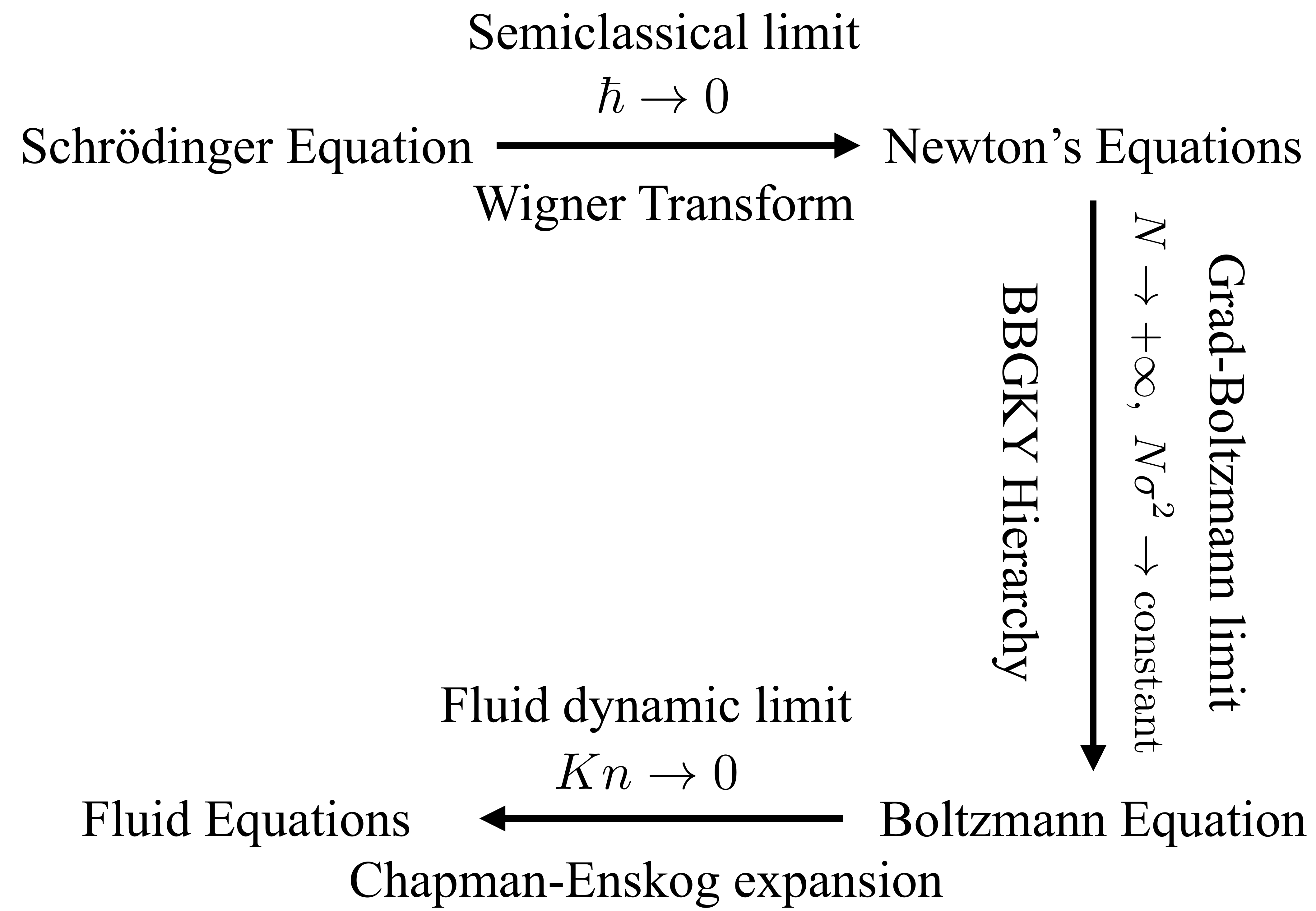}
%}
\caption{Scaling limits from microscopic to macroscopic models.}
\label{scaling}
\end{figure}

%Theorem by Golse-Jin-Levermore here
Since the design of AP schemes relies upon a good understanding of the asymptotic transitions from the microscopic to the macroscopic models, in the next section we first review such transitions for some of the most fundamental physical equations and scalings shown in Fig. \ref{fig-2}.  They are summarized in Fig. \ref{scaling}.

\section{Asymptotic transitions from microscopic to macroscopic physics}\label{sec12}
\setcounter{equation}{1}

\subsection{From quantum mechanics to classical mechanics}

Consider the dimensionless Schr\"odinger equation from quantum mechanics:
\begin{equation}
\I \e \partial _t u^\e =  -\frac{\varepsilon^2}{2} \Delta u^\e  +V(x) u^\e; 
\quad u^\e(0,x) = u^\e_{\rm in}(x)\,.
\label{Sch-eq}
\end{equation}
Here $u^\e=u^\e(t,x)\in \bC$ is complex-valued  quantum mechanical wave function, $(t,x)\in \bR \times \bR^d$, with $d\in \bN$ denoting the spatial dimension. In addition, $\e>0$ denotes the small {\it semiclassical parameter} (the scaled Planck's constant $\hbar$), describing 
the microscopic/macroscopic scale ratio. In quantum mechanics for $N$ particles,  $V(x)$ is the Coulomb potential, but here it is left as a general function of $x$.

The physical observables are real-valued quadratic quantities of $u^\e$.  They include 
  the {\it position density}
\begin{equation}\label{rho}
\rho^\e(t,x):= |u^\e(t,x)|^2, 
\end{equation}
 the {\it current density}
\begin{equation}\label{J}
\quad j^\e(t,x) := \e \im\big(\overline{u^\e}(t,x)\nabla u^\e(t,x)\big),
\end{equation}
and the {\it energy density}
\begin{equation}\label{en}
e^\e(t,x):=  \frac{1}{2} | \e \nabla u^\e (t,x) |^2 + V(x) \rho^\e(t,x).
\end{equation}
Simple analysis shows these observables are governed by the following dynamics:
\begin{align}\label{con}
%\left\{
 \begin{array}{l}
\displaystyle  \partial_t \rho^\e + \nabla \cdot j^\e = 0, \\[8pt]
\displaystyle \partial_t j^\e + \nabla \cdot \left[ \frac{j^\e \otimes j^\e}{\rho^\e} \right] + \rho\nabla V= \frac{\e^2}{2}\rho^\e \nabla 
\left( \frac{1}{\sqrt{\rho^e}} \Delta \sqrt{\rho^\e}\right),  \\[8pt]
 \displaystyle\partial_t e^\e + \nabla \cdot \left( \frac{j^\e}{\rho^\e}(e^\e+\rho^\e V-V\rho^\e)\right) = \frac{\e^2}{4}\nabla \cdot \left[ \frac{j^\e\Delta \rho^\e}{\rho^\e}-\frac{\nabla\cdot j^\e \nabla\rho^\e}{\rho^\e}\right].
\end{array}
\end{align}

From here one can easily deduce the conservation in time of total mass and energy:
\begin{equation}
\label{conservation}
\partial_t \int_{\bR} \rho^\e\, dx=0, \qquad
\partial_t \int_{\bR} e^\e\, dx=0.
\end{equation}

The two main computational challenges to the Schr\"odinger equation
are:
\begin{itemize}
\item[(i)] {\it small $\e$}.  $u^\e$ oscillates with frequency $1/ \e$ in {\it both} space {\it and} time, hence one needs to numerically resolve these oscillations, both spatially and temporally.
\item[(ii)] {\it Large $d$}. For a system consisting of $N$ particles, $d=3N$. Typically $N$ is large. For example for the carbon dioxide molecule $d=75$. For benzene molecule $N=162$. This
causes the curse-of-dimensionality. Totally different 
techniques need to be used for such high dimensional problems and we shall not elaborate on these issues in this paper. 
\end{itemize}

In this survey we focus on the first challenge, namely how one
numerically deals with the small $\e$ problem efficiently.
To this aim, we first review the so-called "{\it semi-classical }" approximation.

%%%%%%%%%%%%%%%%%%%%%%%%%%%%%%%%%%%%%%%%%%%%%%%%%%%%%%%%%%%%%%%%%%%%%%%%%%%%%%%%%

\subsubsection{The WKB analysis }\label{sec: analysis}

Consider the initial data of the following form (the so-called WKB initial data)
\begin{equation}
    u^\e(0,x)=A_0(x) e^{i S_0(x)/\e}.
\end{equation}
The WKB analysis assumes that the solution remains the same form at later time:
\begin{equation}
    u^\e(t,x)=A(t,x) e^{i S(t,x)/\e}.
\end{equation}
Here $A$ is the amplitude, and $S$ is the phase. Applying this {\it ansatz}, which is also called the Madelung transform,  to the Schr\"odinger equation (\ref{Sch-eq}), and separating the real part from the imaginary part,
one gets
\begin{eqnarray}
\label{WKB-parts}
 \begin{aligned}
 & A \partial_t S = \frac{\e^2}{2} \Delta A - \frac{1}{2}A|\nabla S|^2 -AV\,,\\
& \partial_t A = - \nabla A \cdot \nabla S - \frac{1}{2} A \Delta S\,.
 \end{aligned}
 \end{eqnarray}
 Ignoring the $O(\e^2)$ terms,  one gets
  \begin{eqnarray}
\label{Eic-Transp}
 \begin{aligned}
 & \partial_t |A|^2 + \nabla (|A|^2 \nabla S) =0\,,\\
 &\partial_t S + \frac{1}{2}|\nabla S|^2 +V=0\,.\\
 \end{aligned}
 \end{eqnarray}
 The first equation above is called the {\it transport equation}, while the second is the {\it eiconal equation}.  The eiconal equation is a Hamilton-Jacobi equation which admits solutions $S$ with discontinuous derivatives. This can be easily seen once one
 takes a gradient on the equation to get (by letting $u=
 \nabla S$)
 \begin{equation}\label{Burgers}
     \partial_t u + u \cdot \nabla u + \nabla V=0\,.
 \end{equation}
 This is the inviscid Burgers' equation (with forcing term $-\nabla V$) that admits discontinuous
 solution (shocks)  to $u$ even if the initial data of $u$ is smooth.
 Consequently the gradient of $S$ becomes discontinuous, a point
 usually refereed to as the {\it caustic}.
 
 By defining $\rho=|A|^2$, system (\ref{Eic-Transp}) can be written
 as
   \begin{eqnarray}
\label{pressureless}
 \begin{aligned}
 & \partial_t \rho + \nabla \cdot ( \rho u ) =0\\
 &\partial_t (\rho u) + \nabla \cdot (\rho u\otimes u) + \rho  V=0.\\
 \end{aligned}
 \end{eqnarray}
 This is the pressureless gas system. Clearly the system is decoupled.
 One can solve the second equation (which is actually \eqref{Burgers}) for $u$ and then obtain $\rho$
 from the first equation. When $u$ becomes discontinuous,
 $\rho$ becomes a Dirac Delta function, usually called a {\it delta shock} \cite{TanZheng}. Thus at a caustic, the amplitude $A$ blows up (becomes infinity).
 
Beyond the caustic, one notion of mathematical solution to the Hamilton-Jacobi equation is the {\it viscosity solution}, introduced by Crandall and Lions \cite{CL83}. This notion,
however, cannot be used here since system (\ref{pressureless}) is in fact the $\e \to 0$ limit of system (\ref{con}), which is a
{\it zero dispersion limit}.  Zero dispersion limit is drastically
different from the zero dissipation limit, as studied  for the Korteweg–de Vries (KdV) Equation \cite{LaxLev}. For semiclasical limit of the defocusing nonlinear Schr\"odinger equation see \cite{JLM99}).  Thus the WKB analysis is only valid up to the time when the first caustic forms.  Beyond caustics, the solution becomes
{\it multi-valued} \cite{SpMaMa}.

In contrast to that, the Wigner transform technique, which we study next, yields the Liouville equation on {\it phase space}, in the semiclassical limit $\e \to 0$, whose solution does not exhibit caustics, hence is valid globally in time.

\subsubsection{Classical limit via the Wigner transform} \label{sec: wig}

The {\it Wigner transform} of $u^\varepsilon$ is defined as \cite{Wig} 
\begin{equation} \label{Wig1}
w^\varepsilon [u^\e](x,\xi) :=
\frac{1}{(2\pi )^d}\int_{\bR^d}u^\varepsilon \left(x+\frac{\varepsilon}{2} \eta \right)\overline{{u^\varepsilon }}
\left(x-\frac{\varepsilon}{2} \eta \right)\ex^{i\xi\cdot\eta }\D \eta 
 \end{equation}
which is the Fourier transform of the density matrix.

It is easy to see that 
the  Wigner transform $w^\e $ is real-valued, but in general not necessarily positive. The moments of $w^\e$ give
 the quantum mechanical  physical
observables. For examples,  
the particle density \eqref{rho} can be computed via
\[
\rho^\e(t,x) = \int_{\bR^d}
w^\varepsilon (t, x,\xi) \D \xi ,
\]
the current density \eqref{J} can be obtained  by
\[
j^\e(t,x) = \int_{\bR^d} \xi
w^\varepsilon (t, x,\xi) \D \xi ,
\]
while the energy density \eqref{en} is just
\[
e^\e(t,x) =  \int_{\bR^d} \left(\frac{1}{2}|\xi|^2 +V(x)\right)
w^\varepsilon (t, x,\xi) \D \xi.
\]

Applying the Wigner transformation to the Schr\"odinger equation \eqref{Sch-eq},
one obtains the Wigner equation (also called the quantum Liouville equation): 
\begin{eqnarray} \label{wignereq}
\partial _tw^\varepsilon+\xi \cdot \nabla_x w^\varepsilon -\Theta ^\varepsilon
[V]w^\varepsilon  =   0, \quad
w^\varepsilon (0, x,\xi)  =   w^\e_{\rm in}(x,\xi),
\end{eqnarray}
where $\Theta^\varepsilon  [V]$ is given by
\begin{equation}
\Theta^\varepsilon  [V]f(x,\xi) := 
\frac{i}{ (2\pi )^d}\iint_{\bR^{d}\times \bR^d}
\delta V^\e(x,y) 
f(x,\xi') \ex^{i\eta (\xi -\xi')}\D \eta \D \xi ' \,,
\label{theta-def}
\end{equation}
with
$$
\delta V^\e:=\frac{1}{\e} \left(V \left(x-\frac{\varepsilon}{2}y \right)-V \left(x+\frac{\varepsilon}{2}y \right)\right).
$$
When $\e \to 0$
$$\delta V^\e \stackrel{\e\rightarrow 0
}{\longrightarrow} y \cdot \nabla_x V,$$  then \eqref{wignereq} formally becomes the classical {\it Liouville equation} on phase space:
\be\label{lio1}
\partial_t w + \xi  \cdot \nabla_\xi w - \nabla_x V(x) \cdot \nabla_\xi w = 0\,.
\ee
This is the classical limit of the Schr\"odinger equation as $\e\to 0$, valid {\it globally} in time, even beyond the caustic \cite{LiPa,GMMP}, in contract to the WKB analysis. Note that the Liouville equation 
(\ref{lio1}) is linear, which unfolds the singularity, and the linear superposition and time reversibility of the Schr\"odinger equation are also preserved.

The (bi)characteristic equations for (\ref{lio1}) is given by 
\begin{equation*}
  \dot x = \xi ,  \quad \dot \xi= - \nabla_x V(x),
\label{H-system}
\end{equation*}
which is exactly Newton's equation. 
This system can be written as a Hamiltonian system:
\begin{equation}
\label{hamode}
\left \{
\begin{aligned}
& \,  \dot x = \nabla_\xi H(x,\xi), \\
& \, \dot \xi= - \nabla_x H(x,\xi), 
\end{aligned}
\right. 
\end{equation}
with the Hamiltonian $H$ (in classical mechanics) given by 
\be\label{clham} 
H(x,\xi) = \frac{1}{2}|\xi|^2 +V(x).
\ee

For $x, \xi \in \mathbb{R}^{dN}$, and $V$ the potential for $N$-particles, (\ref{hamode}) is the particle system to be studied in the next subsection.

\subsection{From classical mechanics to kinetic equations}

\subsubsection{From hard sphere particles to the Boltzmann equation}

Consider N particles of hard sphere,
\begin{eqnarray}\nonumber
&\dot{ x}_i =  v_i,\\
\nonumber & \dot{ v}_i = 0, 
\end{eqnarray}
where $(x_i, v_i)\in {\mathbb{R}^{d}}\times {\mathbb{R}^{d}},\, (1\le i \le N)$ denote the position and velocity of particle $i$.
 Assume each particle has the same  diameter $\sigma$, then they satisfy the exclusion condition
\begin{equation}
    |x_i(t)-x_j(t)|>\sigma.
\end{equation}
 Assume particles $i$ and $j$ collide elastically when  $|x_i-x_j|=\sigma$,  then the  post-collisional velocities, denoted by 
$v_i'$ and $v_j'$ respectively, are given by
\begin{align} \label{velocities}
%\left\{
\begin{array}{l}
\displaystyle v_i'=v_i-[(v_i-v_j)\cdot\omega] \omega\,, \quad
\displaystyle v_j'=v_j+[(v_i-v_j)\cdot \omega]\omega\,,
\end{array}
\end{align}
where $\omega=(x_j-x_i)/|x_j-x_i|$. 

Define 
\begin{equation}
Z_N=(z_1, \cdots, z_N)=(x_1, v_1, \cdots, x_N,  v_N).
\end{equation}
Let $W^N(t,Z_N)$ be the probability distribution of the particle system.  Then it solves the N-body Liouville equation
\begin{equation}
\partial_t W^N(t)={\cal L}_N W^N(t),  \quad {\cal L}_N=-\sum_{i=1}^N (v_i\cdot \nabla_{x_i} ),
\end{equation}
which is defined on the domain
$$
{\mathcal {D}}_N=\left\{ Z_N \in {\mathbb R}^{2dN} \Big| \,|x_i-x_j|>\sigma,  \quad{\hbox {for}}\,  i\not= j \right\}.
$$
At the boundary, where $|x_i-x_j|=\sigma$, one has $W^N(t, Z_N')=W^N(t,Z_N)$. 

Assume all particles are {\it identical and indistinguishable}, namely
$$
W^N(z_1, \cdots, z_N)= W^N(z_{\sigma_1}, \cdots, z_{\sigma_N})
$$
for any $\{\sigma_1, \cdots, \sigma_N\}$, a random permulation of set $\{1, \cdots, N\}$.
Furthermore, assume the so-called {\it molecular chaos} condition:
$$
W^N(z_1, \cdots z_N)=W_1^N(x_1, v_1) \cdots W_N^N(x_n, v_N).
$$
Under the above assumptions, 
the Grad-Boltzmann limit of classical particles  can be derived by letting  $\sigma \to 0$, and $N \to \infty$, and under the assumption
$$
N \sigma^2 \to {\hbox {constant}}\,,
$$
then the one-particle distribution $W_1^N(x_1, v_1)$ approaches formally to the Boltzmann equation for hard spheres \cite{BGP}:
\begin{equation}
\partial_t f +v\cdot\nabla_x f 
=  \int |(v-v_*)\cdot \omega|  \{ f(v')f(v'_*)-(f(v)f(v_*)\} d\omega dv_*\,,
\label{tone-1}
\end{equation}
where $v$ and $v_*$ are pre-collisional velocites with corresponding post-collisional velocities $v'$ and $v_*'$ determined by (\ref{velocities}) (with $v_i=v, v_j=v_*$).

The proof of the Grad-Boltzmann limit is extremely challenging.  So far the only rigorous results are available for a very short duration of time--a fraction of a mean free time, see \cite{Lanford} \cite{GST}.

\subsubsection{Mean-field limit of particle systems}

The Newton type equations also arise in microscopic modeling of a vast number of important phenomena in physical, social, and biological sciences, 
\cite{vicsek1995novel} \cite{cucker2007emergent} \cite{motsch2014}\cite{Ha-Review} .
These problems can all be modelled by interacting particle systems of first order
\begin{gather}\label{eq:Nparticlesys-0}
dX^i=b(X^i)\,dt+\alpha_N\sum_{j: j\neq i} K_1(X^i-X^j)\,dt+\eta\, dW^i,~~i=1,2,\cdots, N,
\end{gather}
or second order 
\begin{eqnarray}
\label{eq:Nbody2nd-0}
\begin{aligned}
& dX^i=V^i\,dt,\\
& dV^i=\Big[ b(X^i)+\alpha_N \sum_{j:j\neq i}K(X^i-X^j)-\gamma V^i \Big]\,dt+\eta\, dW^i.
\end{aligned}
\end{eqnarray}
Here, $(X^i, V^i)\in \bR^d \times \bR^d$, loosely speaking,  represent the position and velocity of the $i$-th particle,  and $b(\cdot)$ is  the external field. The stochastic processes $\{W^i\}_{i=1}^N$
are i.i.d. Wiener processes, or the standard Brownian motions.
 The function $K(\cdot): \mathbb{R}^d\to \mathbb{R}^d$ is the interaction kernel. For the molecules in the heat bath,  $\eta$ and $\gamma$ satisfy the so-called ``fluctuation-dissipation relation"
\begin{gather}
\eta=\sqrt{2\gamma/\beta},
\end{gather}
where $\beta$ is the inverse of the temperature (we assume all the quantities are scaled and hence dimensionless so that the Boltzmann constant is absent). The first order system \eqref{eq:Nparticlesys-0} can be viewed
as  the  over-damped limit, namely rescale $t$ to $\gamma t$ and let $\gamma \to \infty$, of the second order system \eqref{eq:Nbody2nd-0} \cite{stanley1971}\cite{georges1996}\cite{lelievre2016partial}.

The mean field limit is usually taken by choosing
\begin{gather}
\alpha_N=\frac{1}{N-1}.
\end{gather}
Define the empirical distribution 
\begin{equation}
    \label{emprical}
 \mu^{(N)}:=\frac{1}{N}\sum_{i=1}^N \delta(x-X^i)\otimes\delta(v-V^i).
 \end{equation}
 For the second order system \eqref{eq:Nbody2nd-0}, as $N\to\infty$, $\mu^{(N)}$  converges almost surely in the weak topology to the solutions of the (mean-field) Fokker-Planck equation
\begin{gather}\label{eq:limitvlasov}
\partial_t f=-\nabla_x\cdot(vf)-\nabla_v\cdot\Big((b(x)+K*_xf-\gamma v)f\Big)+\frac{1}{2}
\eta^2\Delta_v f.
\end{gather}
The mean field limit corresponding to the first order system 
(\ref{eq:Nparticlesys-0}) is (\cite{mckean1967} \cite{golse2003mean} \cite{jabin2017mean})
\begin{equation}
    \label{eq:eq:limitvlasov-1st}
\partial_t f=-\nabla \cdot\Big((b(x)+K_1*f)f\Big)+\frac{1}{2}
\eta^2\Delta f.
\end{equation}

These mean field limits can also be derived from taking the limit $N\to \infty$ of the $N$-body distribution with molecular chaos assumptions, like in the derivation of the Boltzmann equation from the $N$-body Newton's equations described above. 

\subsection{From kinetic equations to hydrodynamics}

\subsubsection{Hydrodynamic limit of the Boltzmann equation}
The Boltzmann equation describes the probability density function $f(t,x,v)$ of particles  that undergo transport and binary collisions \cite{Cercig}: 
\begin{equation} \label{CBE}
\partial_t f +v\cdot \nabla_x f=\frac{1}{\varepsilon}\mathcal{Q}(f), \quad x, v\in \mathbb{R}^n,
\end{equation}
where the collision term $\mathcal{Q}(f)$ is a nonlinear integral operator:
\begin{equation} 
\mathcal{Q}(f)(v)=\frac{1}{\e}\int_{\mathbb{R}^{d}}\int_{S^{d-1}}B(v-v_*,\omega)\left[f(v')f(v_*')-f(v)f(v_*)\right]\,\rd{\omega}\,\rd{v_*}.
\end{equation}
Here $(v,v_*)$ and $(v',v_*')$ are the velocity pairs before and after an elastic collision,  which conserve the momentum and energy. They are related by 
\begin{align*}
%\left\{
\begin{array}{l}
\displaystyle v'=v-[(v-v_*) \cdot \omega] \omega\,, \quad 
\displaystyle v_*'=v_*+[(v-v_*)\cdot\omega] \omega\,,
\end{array}
\end{align*}
with the parameter $\omega \in S^{d-1}$, the unit sphere on $\mathbb{R}^d$. $B(v-v_*, \omega)$ is the (non-negative) collision kernel depending only on $|v-v_*|$ and cosine of the deviation angle $\frac{\sigma\cdot (v-v_*)}{|v-v_*|}$. $\e$ is the Knudsen number,
the dimensionless mean free path.

The hydrodynamic quantities $\rho, u$ and $T$, the density, marcoscopic velocity and temparature respectively, are defined as
 the moments of $f$:
\begin{eqnarray}
\label{moments-def}
&&\rho=\int_{\mathbb{R}^{d}} f \,\rd{v}
  = \int_{\mathbb{R}^{d}} \mathcal{M} \,\rd{v},\quad 
  u=\frac{1}{\rho}\int_{\mathbb{R}^{d}} vf \,\rd{v}
  = \int_{\mathbb{R}^{d}} v\mathcal{M} \,\rd{v},\\
&&
T=\frac{1}{d\rho}\int_{\mathbb{R}^{d}} |v-u|^2 f \,\rd{v}
  = \frac{1}{d\rho}\int_{\mathbb{R}^{d}} |v-u|^2 \mathcal{M} \,\rd{v},
\end{eqnarray}
where the local Maxwellian
\begin{equation}
\label{local-Max}
    \mathcal{M}= \frac{\rho}{(2 \pi T)^{d/2}} {\text {exp}}\left(-\frac{|u-v|^2}{2T}\right)\,.
\end{equation}
The collision operator $\mathcal{Q}(f)$ conserves mass, momentum, and energy:
\begin{equation} \label{consv1}
\int_{\mathbb{R}^{d}}\mathcal{Q}(f)\phi(v)\,\rd{v}=0, \quad \phi(v)=(1,v,|v|^2/2)^T,
\end{equation}
with the momentum $m=\rho u$, and the total energy $E=\frac{1}{2} \rho u^2 + \rho T$.

One of the most important properties of $\mathcal{Q}$ is  the celebrated Boltzmann's {\it H}-theorem:
\begin{equation*}
\partial_t \int_{\mathbb{R}^{d}} f \log f\,\rd{v}= \int_{\mathbb{R}^{d}}\mathcal{Q}(f)\ln f\,\rd{v} \le 0.
\end{equation*}
The functional $f\log f$ is the entropy of the system. Boltzmann's $H$ theorem implies that any equilibrium distribution
function, i.e., any function which is a maximum of the entropy, has the form of a local Maxwellian distribution
\begin{equation} \label{equi}
\int_{\mathbb{R}^{d}}\mathcal{Q}(f)\ln f\,\rd{v} =0 \Longleftrightarrow \mathcal{Q}(f)=0\Longleftrightarrow f=\mathcal{M}.
\end{equation}

When $\e \to 0$, $\mathcal{Q} \to 0$,  (\ref{equi}) implies that $f=\mathcal{M}$. Consequently the moments of $f$ solve the compressible Euler equations:
\begin{align} \label{Euler-limit}
\left\{
\begin{array}{l}
\displaystyle \partial_t \rho+\nabla_x\cdot (\rho u)=0,\\[8pt]
\displaystyle \partial_t (\rho u)+\nabla_x\cdot (\rho u\otimes u +pI)=0,\\[8pt]
\displaystyle \partial_tE+\nabla_x\cdot ((E+p)u)=0.
\end{array}\right.
\end{align}
 Via the Chapman-Enskog expansion, one can derive the Navier-Stokes (NS) equations by retaining $O(\varepsilon)$ terms \cite{BGP00}:
\begin{align} \label{NS}
\left\{
\begin{array}{l}
\displaystyle \partial_t \rho+\nabla_x\cdot (\rho u)=0,\\[8pt]
\displaystyle \partial_t (\rho u)+\nabla_x\cdot (\rho u\otimes u +pI)=\varepsilon \nabla_x\cdot (\mu \sigma(u)),\\[8pt]
\displaystyle \partial_tE+\nabla_x\cdot ((E+p)u)=\varepsilon \nabla_x\cdot (\mu \sigma(u)u+\kappa \nabla_xT),
\end{array}\right.
\end{align}
where $\sigma(u)=\nabla_xu+\nabla_xu^T-\frac{2}{d}\nabla_x\cdot u I$, $I$ is the identity matrix, $\mu$ and $\kappa$ are the viscosity and heat conductivity, determined through the linearized Boltzmann collision operator, and usually depend on $T$.
%\begin{remark}
%The BGK model does not yield the satisfactory result at the NS level in that it gives the wro%ng Prandtl number (the ratio of $\mu$ and $\kappa$) \cite{Cercignani}. 
%\end{remark}

\subsubsection{Diffusion limit of transport equation}

In many applications, such as neutron transport and radiative transfer,
the collision operator is linear. The interesting scaling is  often
the diffusive scaling where the scattering rate is large.  A typical such equation has the form of
\begin{equation}
\varepsilon\,\partial_t f +v\cdot\nabla_x f 
= \frac{1}{\varepsilon} \int b(v,w) \{ M(v)f(w)-M(w)f(v)\} dw\,, \quad x,v \in \mathbb{R}^d
\label{tone-2}
\end{equation}
with the normalized Maxwellian $M$ defined by
$$
  M(v) = \frac1{(2\pi)^{d/2}} \exp ( - | v|^2/2) \,.
$$
The (anisotropic) scattering kernel $b$ is  rotationally invariant,  satisfying 
$$
  b(v,w) =   b(w,v)  > 0 \,.
$$
Define the collision frequency $\lambda$ as
$$
\lambda(v) = \int b(v,w)M(w)\, dw  \,.
$$

  As $\varepsilon \to 0$, $f\to \rho (x, t) M(v)$, 
 where $\rho(t,x)=\int f(v)\, dv$
satisfies the diffusion equation \cite{BSS84} \cite{MRS-book} 
\begin{equation}
\partial_t\rho= \nabla_x \cdot (D \nabla_x \rho)
\label{televen}
\end{equation}
with the diffusion coefficient matrix
\begin{equation}\label{diff-coeff}
D=\int \frac{M(v)}{\lambda(v)}{v\otimes v}\, dv\,.
\end{equation}

%%%%%%%%%%%%%%%%%%%%%%%%%%%%%%%%%%%%%%%%%%%%%%%%%%%%%%%%%%%%%%%%%%%%%%%%%

\section{Numerical passages from quantum to classical mechanics}\label{sec2}
\setcounter{equation}{1}

The highly oscillatory nature of the solution, in both space and time, to the Schr\"odinger equation (\ref{Sch-eq}) poses a huge challenge in numerical computations,   especially in high dimensions, since one needs to resolve numerically,
in both space and time, the small wave length of $O(\e)$ which is computationally daunting.
If one does not use small enough  time step or mesh size, even if the numerical scheme is stable, one may get  completely wrong solutions \cite{MAPC}\cite{BJM}. To understand the numerical behavior in the semiclassical regime, in addition to standard consistency and stability--which implies convergence by Lax's equivalence theorem-- one needs new semiclassical analysis to understand the correct behavior of the numerical solutions.

Here one is interested in two questions:
\begin{itemize}

\item What kind of schemes  best suits the highly oscillatory problems? 

\item How to analyze the numerical performance when $\e$ is small?

\end{itemize}

For the first question,  when the solution is smooth but highly oscillatory,  the spectral or pseudo-spectral methods give the best performance in terms of numerical accuracy and resolution. It is worthy to point out that taking care of spatial discretization alone is not enough to achieve the best performance  for the Schr\"odinger equation (\ref{Sch-eq}). It takes a 
good combination of {\it both} spatial {\it and} temporal discretizations to achieve the most favorable mesh strategies (the largest possible ratio between the mesh size and time steps over $\e$). In this regard the time-splitting spectral methods, as studied  in \cite{BJM}, offers the best mesh strategy, while finite-difference type schemes require very fine numerical resolution of the oscillations \cite{MAPC}.

%%%%%%%%%%%%%%%%%%%%%%%%%%%%%%%%%%%%%%%%%%%%%%%%%%%

\subsection{Time-splitting spectral methods for the semiclassical Schr\"odinger equations} \label{sec: split}

For the sake of notation clarity, we shall discuss 
the method only in one space dimension ($d=1$). Generalizations to $d>1$ are straightforward for tensor  
product grids and the same conclusions hold. 

Consider the one-dimensional version of equation   \eqref{Sch-eq},
\be \label{1dschro}
\I \vep \partial_tu^\vep = - \fl{\vep^2}{2} \p_{xx}u^\vep+V(x)u^\vep, \quad u^\vep(0,x)=u_{\rm in}^\vep(x), 
\ee
for $x\in [0,1]$, with periodic boundary conditions
$$ 
u^\vep(t, 0)=u^\vep(t, 1), \quad  \p_x u^\vep (t, 0)=\p_x u^\vep(t, 1), \quad \forall \, t\in \bR.
$$  
We choose the spatial mesh size $\Delta x =1/M$ for some large positive integer $M$, and  time-step $\btu t >0$. The spatio-temporal grid-points 
are then given by
\[ x_j:=j \Delta x,\ j =1, \dots, M, \qquad t_n := n \Delta t, \  
n\in \bN . \] 
Let $u^{\vep,n}_j$ be the numerical approximation of $u^\vep(x_j,t_n)$, for $j =1, \dots, M$. 
  
%\subsubsection{First-order time-splitting spectral method (SP1)} 

The Schr\"{o}dinger equation \eqref{1dschro} is solved by  a time splitting method:\\

\textbf{Step 1.} From time $t=t_n$ to time $t=t_{n+1}$  first solve the free Schr\"odinger equation
\be  \label{fstep}
\I \vep \p_tu^\vep+ \fl{\vep^2}{2} \p_{xx}u^\vep=0.
\ee 

\textbf{Step 2.} Also on  $t \in [t_n, t_{n+1}]$,  solve the ordinary differential equation
\be \label{sstep}
\I \vep \p_tu^\vep- V(x)u^\vep=0, 
\ee 
with the solution $u^{\vep, *}$ obtained from Step 1  as initial data. 

Note (\ref{sstep}) can be solved {\it exactly},
$$
u(t_{n+1},x) = u(t_n,x) \ex^{-\I V(x)\Delta t/\e}.
$$
In Step 1,  equation (\ref{fstep}) will be discretized  by 
a (pseudo-)spectral method in space and consequently integrated in time {\it exactly} in the Fourier space.  More precisely, 
\begin{equation*}
u_j^{\vep, *}=\fl{1}{M}\sum_{\ell =-M/2}^{M/2-1} 
  \ex^{\I \vep \Delta t\gamma_\ell^2/2}\;\widehat{u}^{\vep,n}_\ell \;  
 \ex^{\I \gamma_\ell (x_j-a)},  
\end{equation*} 
where $ \gamma_\ell=2\pi l$ and $ \widehat{u}^{\vep,n}_\ell$ is the Fourier coefficients of $u^{\vep,n}$, i.e.   
\begin{align*}
\quad \widehat{u}^{\vep,n}_\ell= 
  \sum_{j=1}^{M} u^{\vep,n}_j\;\ex^{-i\gamma_\ell x_j},  
\quad \ell =-\fl{M}{2},\dots, \fl{M}{2}-1\,.
\end{align*} 
Note that in both steps the time integration is {\it exact}. The only time discretization error of this method is the splitting error, which is \emph{first order} in $k$, for any \emph{fixed} $\vep>0$. We will refer this method as TSSP.

%\subsection{Strang splitting (SP2)}

 The second order in time (for fixed $\e>0$) can be obtained via the Strang splitting method. 
Extensions to higher order (in time) splitting schemes can also be done, see e.g. \cite{BaSh}. See also its extension to the case of vector potential \cite{ZhouZN}.

\subsection{Numerical analysis in the semiclassical regime} 
 
Classical numerical analysis,  based on consistency and stability,  does not provide accurate assessment of the numerical performance  when $\e\ll 1$. Wigner analysis, on the other hand, gives more
insight about the behavior of numerical solutions, for {\it physical observables}, in the semiclassical regime.  

Assume that the potential $V(x)$ is periodic in domain $[0,1]$, smooth, and satisfies
\be\label{A} 
\Big \|\fl{d^{m}}{d x^{m}}V\Big \|_{L^\ift[a,b]} 
\le C_m, 
\ee
for some constant $C_m>0$, and furthermore,
\be\label{B} 
\Big \|\fl{\p^{m_1+m_2}}{\p t^{m_1}\p x^{m_2}}u^\vep 
\Big \|_{C([0,T];L^2[a,b])} 
\le \fl{C_{m_1+m_2}}{\vep^{m_1+m_2}}, 
\ee
for all $ m, m_1,\ m_2\in \mathbb{N}\cup \{0\}$, namely the solution  oscillates in space 
and time with wavelength $\vep$. 
The following estimate was given in \cite{BJM}.

\begin{theorem}\label{thm: SP1} 
Let $V(x)$ satisfy assumption \eqref{A} and $u^\vep(t,x)$ be a 
solution of \eqref{1dschro} satisfying \eqref{B}. Denote by 
$u_{\mathrm{int}}^{\vep,n}$ the interpolation of the discrete approximation obtained via TSSP.  Then, for $t_n\in [0,T]$,
\be 
\label{eeu1} 
\left\|u^\vep(t_n)-u_{\mathrm{int}}^{\vep,n}\right\|_{L^2(0,1)}\le 
G_m\fl{T}{\Delta t}\left(\fl{\Delta x}{\vep}\right)^m
+\fl{CT\Delta t}{\vep}, 
\ee 
where $C>0$ is independent of $\vep$ and $m$ and $G_m>0$ is independent of $\vep$, $\Delta x$, $\Delta t$. 
\end{theorem} 

Clearly, (\ref{eeu1}) implies that, to get an accurate $u^\e$ one needs the following mesh strategy:
$$
\Delta t=o(\e), \quad \Delta x=o(\e)\,.
$$
Hence  the oscillations need to be resolved both spatially and temporally.

\subsection{Accurate computation of  physical observables} 

 If one is just interested in obtaining accurate physical observables, it was observed in \cite{BJM} that the time step can be much relaxed. This cannot be understood from the above classical numerical analysis, rather  the Wigner picture of quantum dynamics will offer the clue. 
 
Let $u^\vep(t,x)$ be the solution of \eqref{1dschro} and $w^\varepsilon(t,x,\xi)$ the corresponding Wigner transform.
It is easy to see that the  splitting scheme (\ref{fstep})-(\ref{sstep}) corresponds to  
the following time-splitting scheme for the Wigner equation \eqref{wignereq}:\\

{\bf Step 1.} For $t \in [t_n, t_{n+1}]$, first solve the linear transport equation
\be
\label{fwige} 
\p_t w^\vep+\xi\, \partial_x w^\vep=0\,.
\ee 

{\bf Step 2.} On the same time-interval, solve the scattering term
\be 
\label{swige} 
\p_t w^\vep-\Theta^\vep[V]w^\vep=0\,,
\ee 
with initial data obtained from Step 1.\\

Since in each step of the splitting, the time integration is {\it exact}, without any discretization error, thus one can take 
 $\vep\to 0$ limit in each step,  consequently obtain the following limiting scheme:\\

{\bf Step 1.}  For $t \in [t_n, t_{n+1}]$ solve
\begin{equation}
\label{lvla-1} 
\p_t w+\xi\, \partial_xw^0=0.
\end{equation} 

{\bf Step 2.} Using the outcome of Step 1 as initial data, solve, on the same time-interval:
\begin{equation}
\label{lvla-2} 
\p_t w- \partial_xV\, \partial_\xi w^0=0.
\end{equation}

This is {\it exactly} the time-splitting scheme for the  limiting Liouville equation \eqref{lio1}!
Since in the limiting process  $\Delta t$ was held fixed, hence independent of $\e$, thus when $\Delta t=O(1)$, and as $\e \to 0$,
schemes (\ref{fwige}) and (\ref{swige}) collapse to schemes
(\ref{lvla-1}) and (\ref{lvla-2}), therefore the scheme is AP in time!
 Hence one can take $\Delta t =O(1)$,  
combined with the spectral mesh-size $\Delta x=o(\e)$ to get   accurate $w^\e$, and as a consequence, {\it all physical observables}!

\begin{remark}
  While the above Wigner analysis is formal, rigorous uniform in $\e$ error estimate was obtained recently. In \cite{GJP21}, for both first and second order splittings, a uniform in $\e$ error estimates, with explicit constants, were given (for the von Neumann equation--the density operator representation of the Schr\"odinger quation which are valid even for mixed states). The errors are measured by  a pseudo-metric introduced in \cite{GP17}, which is an analogue of the  Wasserstain distance of exponent $2$ between a quantum density operator and a classical density in phase space. The regularity requirement for $V$ is $V\in C^{1,1}$.
  Sharper uniform error estimates for physical observables were also obtained for the Strang splitting \cite{LaLu20},  based on Egorov's theorem, with additional regularity requirement on $V$.
\end{remark}

\begin{example}\label{ex-1}
We take one example from \cite{BJM}. The Schr\"odinger equation \eqref{1dschro} is solved with initial condition $u_0(x)=n_0(x)\exp(-\I S_0(x)/\e)$, where  
\be 
\label{inite1} 
n_0(x)=\left(e^{-25(x-0.5)^2}\right)^2, \quad 
S_0(x)=-\fl{1}{5}\ln\left(e^{5(x-0.5)}+e^{-5(x-0.5)}\right),
\ee 
and $V(x)=10$. 
Due to the compressive initial velocity $\frac{d}{dx} 
S_0(x)$, caustics will form.  
The weak limits $n^0(x,t)$, $J^0(x,t)$  
of $n^\vep(x,t)$, $J^\vep(x,t)$  
respectively as $\vep\to0$ 
 can be computed  by evaluating 
the zeroth and first order velocity 
  moments of the solution to the Liouville  equation (\ref{lio1}).
As a reference 
we plot them  
at $t=0.54$ (after the caustics formed ) in Fig. \ref{fig:BJM-2}. 
We compare the solutions between CNSP (Crank-Nicolson in time and pseudo-spectral method in space) and TSSP2 (Strang's splitting in time and pseudo-spectral method in space).
  The  mesh size $\Delta x$ is taken in the same order as $\vep=10^{-3}$. 
 One can see that  for CNSP, 
even for $\Delta t=0.0001$, the numerical solution cannot capture the correct 
weak limit. TSSP2 can capture the physical observations correctly
with $\Delta t$ much larger than $\e$.

\begin{figure} [t]
\centering

\includegraphics[width=6cm]{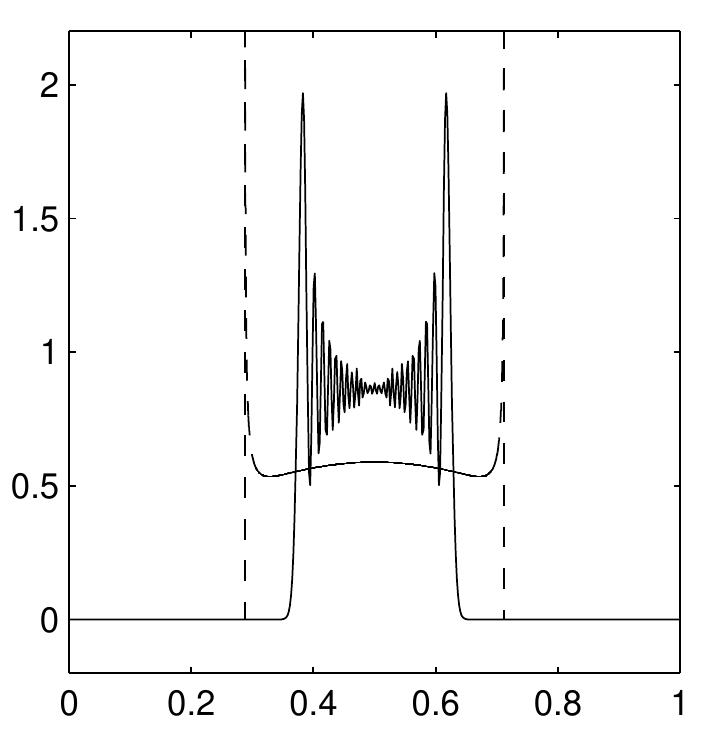}
\includegraphics[width=6cm]{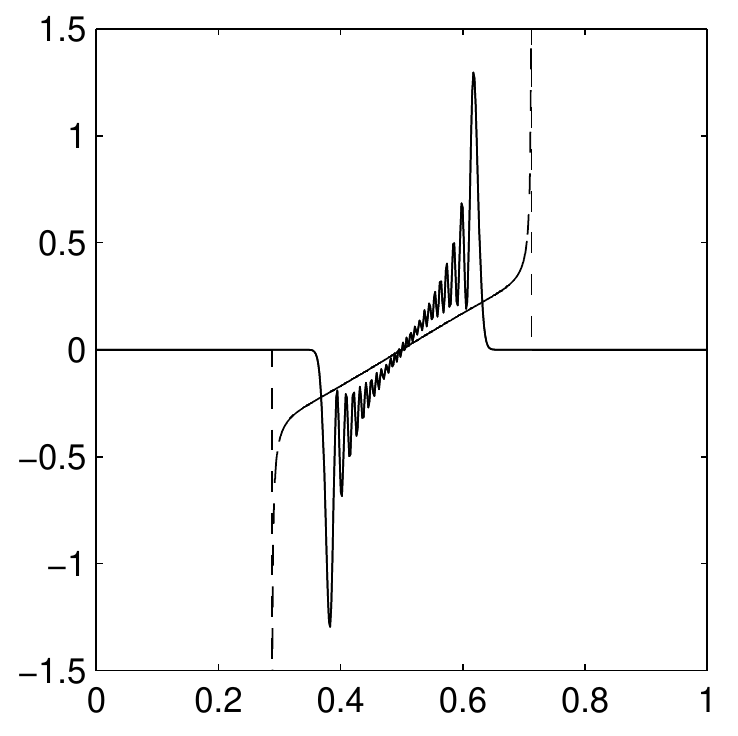}

\includegraphics[width=6cm]{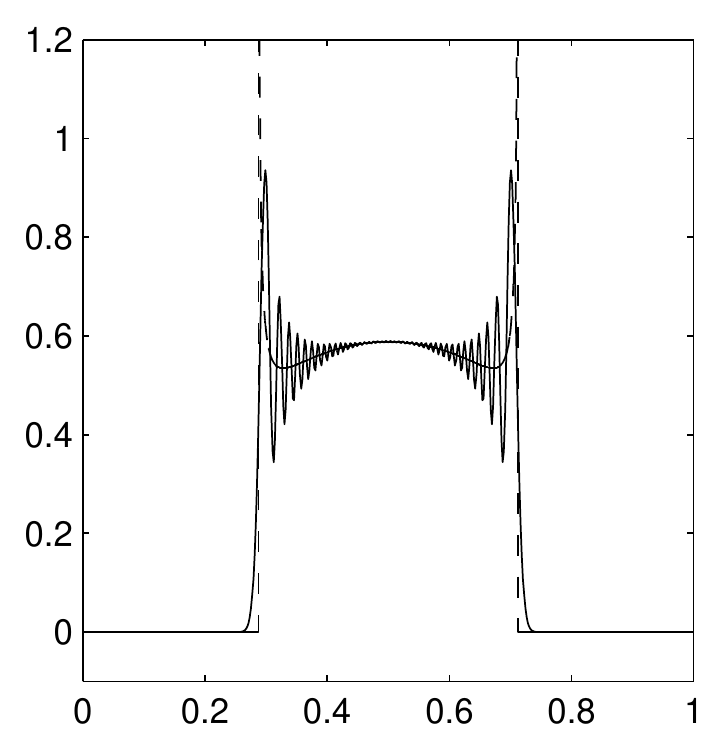}
\includegraphics[width=6cm]{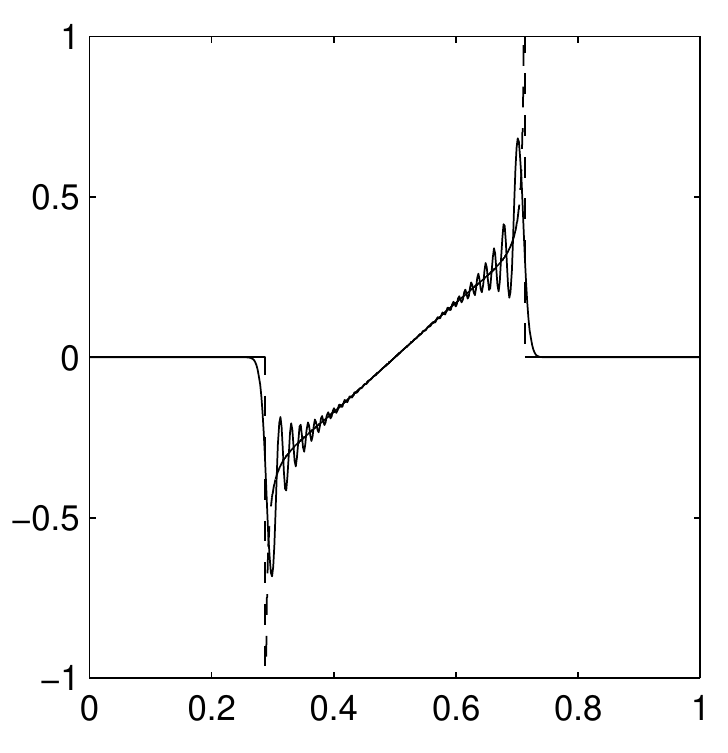}
\caption{Example \ref{ex-1}. Numerical solutions at $t=0.54$.
Top: CNSP with $k=0.0001$; Bottom: TSSP2. Left: position density; Right: current density.
$\vep=10^{-3}$, $V(x)=10$, $h=\fl{1}{512}$. 
  }
\label{fig:BJM-2}
\end{figure}
\end{example}

\bigskip
So far no numerical schemes are known to allow $\Delta x=o(1)$ for the Schr\"odinger equation (\ref{Sch-eq}). The best one can do is
to allow $\Delta x=o(\sqrt{\e})$, by using the Guassian beam or Gaussian wave packet methods, see 
\cite{Hel06}\cite{Hil90}\cite{JWY1}\cite{LQ09}\cite{RussSme}. For more recent results about Gaussian type approximations  see \cite{JMS-Acta}\cite{LaLu20}. 

\subsection{Ehrenfest dynamics}

The {\it ab initio} methods have played indispensable roles in simulating  large systems of quantum molecular dynamics. There the forces acted on the nuclei are computed   from electronic structures, a procedure  known as the
\textquotedblleft on-the-fly" calculation in  chemistry literature (for
detailed reviews, see, e.g., \cite
{tully1998mixed}\cite{marx2009ab}. The \textit{Ehrenfest dynamics} is one of  popularly used such  methods.  There one separates the quantum system into two sub-systems:
 a fast varying, quantum mechanical part for electrons and a slowly varying part for the nuclei.
Due to the large mass difference between electrons and nuclei,  the nucleonic system can be passed to the (semi-)classical limit, hence the computational cost is significantly reduced. 

Take $x\in \bR^d$ as the electronic coordinate, $y\in \bR^n$ the nucleonic coordinate, with $d,n\in \bN$, and 
denote by $\langle \cdot, \cdot\rangle_{L_x^{2}}$ and $\langle \cdot, \cdot\rangle_{L_y^{2}}$ the usual inner product in 
$L^2(\bR^d_x)$ and $L^2(\bR^n_y)$, respectively, i.e. 
\[
\left\langle f, g \right\rangle_{L_z^{2}} \equiv \int_{\mathbb{R}^m} \overline{f}(z)g(z)dz.
\]
The total Hamiltonian of the system acting on $L^2(\bR^{d+n})$ is assumed to be of the form
\begin{equation}\label{eq:ham}
H= - \frac{\vep^2}{2} \Delta_x - \frac{\delta^2}{2} \Delta_y + V(x,y) ,  
\end{equation}
where $V(x,y)\in \bR$ is some real potential.

Consider the following mixed quantum-classical system (corresponding to the limit $\delta \to 0$)
 \cite{tdscf1}\cite{tdscf4} \cite{JSZ17}:
\begin{equation}
\label{SL-1}
\left\{
\begin{aligned}
& i \e \partial _{t}\psi ^{\vep }=-\frac{\vep ^{2}}{2}\Delta _{x}\psi
^{\e }+\Upsilon^{\e }( x,t) \psi ^{\e }, \quad \psi^\e (0,x) = \psi^\e_{\rm in}(x) \\
%\label{SL-2}
& \partial _{t}\mu ^{\e }+\eta \cdot \nabla _{y}\mu ^{\e }+F^{\e
}( y,t) \cdot \nabla _{\eta }\mu ^{\e }=0,  \quad \mu^\e(0,x,\eta) = \mu_{\rm in}(y, \eta).
\end{aligned}
\right.
\end{equation}%
Here, $\mu ^{\e } (\cdot, \cdot, t) $
denotes the phase-space probability density for the slowly varying nuclei at time $t$, $F^\e = - \nabla_y V^\e_{\rm E}$ is the force obtained from the Ehrenfest potential
\begin{equation*}
V^\e_{\rm E}(y,t)=\int_{\mathbb{R}^{d}}V\left( x,y\right) \left\vert \psi^\e (x,t)\right\vert ^{2}\, dx,
\end{equation*}
and
\begin{eqnarray}
\label{potentials}
\Upsilon ^{\e }\left( x,t\right) &=&\iint_{\mathbb{R}^{2n}}V\left(
x,y\right) \mu ^{\e }\left( y,\eta ,t\right) \, dy\, d\eta .
%F^{h }\left( y,t\right) &=&-\int_{\mathbb{R}^{d}}\nabla _{y}V\left(
%x,y\right) \left\vert \psi ^{h }\right\vert ^{2}dx.  \notag
\end{eqnarray}%
This system will be called  the \textit{%
Schr\"{o}dinger-Liouville-Ehrenfest} (SLE) system.
Note that the dependence of $\mu^\e$ on $\e$  is purely from the forcing through the Ehrenfest potential $V_E^\e$ appearing in the Liouville equation.
In the case of a single particle distribution concentrated on the classical trajectories $(y(t), \eta(t))$, i.e.,
\[
\mu(t,y,\eta) = \delta(y-y(t), \eta-\eta(t)),
\]
(\ref{SL-1}) gives \cite{tully1998mixed}
 \cite{Drukk} \cite{schutte1999singular}
\cite{szepessy2011} 
\begin{equation}
\label{SN-1}
\left \{
\begin{aligned}
i\e \partial _{t}\psi^\e &=-\frac{\e ^{2}}{2}\Delta_{x}\psi^\e +V\left( x,y( t) \right) \psi^\e , \quad \psi^\e(0,x) = \psi_{\rm in}(x),\\
%\label{SN-2}
\dot{y} (t)  &= \eta(t), \quad y(0) = y_0,  \\
\dot \eta(t)& = -\nabla _{y}V^\e_{\rm E}( y (t)), \quad \eta(0)=\eta_0 .
\end{aligned}
\right.
\end{equation}
The iterated semiclassical limit ($\delta \to 0$, then $\e\to 0$) and the full classical limit ($\delta=\e\to 0$) were rigorously justified in \cite{JSZ17}.
 
Again, the main numerical difficulty for $\e \ll 1$ here
is that one needs to resolve oscillations of frequency of order $\mathcal O(1/\e)$ in both time and space, as they are
present in the solution $\psi^\e$. This requires one to use
time-steps of order $\Delta t=o(\e)$ as well as a spatial grid with $\Delta x=o(\e)$ to resolve the wave functions. As analyzed in the proceeding subsection, one may ask whether one can design a numerical method which allows the capturing of  {\it physical observables} even for time-steps
much larger than $O(\e)$. 
For {\it nonlinear} Schr\"odinger equations, in general, this is no longer true, as was numerically
demonstrated in \cite{BJM3}.
The SLE system (\ref{SL-1}) is a nonlinearly coupled system, and one therefore expects the same type of problem at first glance. Nevertheless,  an efficient numerical method for
the SLE system was introduced in \cite{FJS18} which allows large
(compared with $\e$) computational mesh-sizes in both $y$ and $\eta$ and a large time step for both the Schr\"odinger and the Liouville
equations, while still {\it correctly capture the physical observables}. While large
meshes in $y$ and $\eta$ do not seem so surprising, since they are coordinates of the nuclei,  the possibility of large time steps for solving
the Schr\"odinger equation for electrons is far from obvious, due to the nonlinear nature
of the SLE system.

\subsubsection{A  time-splitting scheme for the SLE system}\label{sec:tsscheme}

Consider the case $d=n=1$, and the domain $(x,y)\in [0,1]^2$, with
uniform mesh sizes $\Delta y,\Delta \eta $  applied to the classical
part of the SLE \eqref{SL-1}.  Set
\begin{equation*}
J=\frac{1}{\Delta y},\ K=%
\frac{1 }{\Delta \eta }, \ M=\frac{1}{\Delta x}, \ y_{j}=j\Delta y,\ \eta _{k}=
k\Delta \eta, \ x_{j}=j\Delta x.
\end{equation*}

The time-splitting  scheme, introduced in \cite{FJS18}, can then be described as
follows: From time $t=t_{n}=n\Delta t$ to $t=t_{n+1}=\left( n+1\right)
\Delta t$, the SLE system is solved in two steps. First, solve
\begin{equation}\label{step1}
\left\{
\begin{aligned}
& i\e \partial _{t}\psi ^{ \e}  =-\frac{\e ^{2}}{2}\Delta _{x}\psi
^{\e },   \\
& \partial _{t}\mu ^{\e }= -\eta \cdot \nabla _{y}\mu ^{\e }-F^{\e
}\left( y,t\right) \cdot \nabla _{\eta }\mu ^{\e },
\end{aligned}
\right.
\end{equation}%
from $t=t_{n}$ to an intermediate time $t_{\ast }=t_n +\Delta t$. Then, solve
\begin{equation}\label{step2}
\left\{
\begin{aligned}
& ih \partial _{t}\psi ^{h }=\Upsilon^{h} \left( x,t\right) \psi
^{h },  \\
& \partial _{t}\mu ^{h }=0,
\end{aligned}
\right.
\end{equation}
with initial data computed from Step 1, to obtain the solution at time $t=t_{n+1}$.

In (\ref{step1}), the Schr\"odinger equation will be discretized
in space by a spectral method using the Fast Fourier Transform, and integrated in the Fourier space in time exactly.
The Liouville equation can be solved either by a spectral method, or by a finite difference
(e.g., upwind) scheme in space, and then marching the corresponding ODE system forward in time.
An advantage of this splitting method is that in the second step, $\Upsilon ^{h }\left(
x,t\right) $ defined in (\ref{potentials}) is  {\it independent} of time, since obviously $\mu ^{h }$ is.
Hence, the time integration in \eqref{step2} can also be solved {\it exactly} as
\begin{equation*}
\psi _{j}^{\e ,n+1}=\exp \left( -\frac{i}{\e }\Upsilon ^{\e
}\left( x_{j},t_{\ast }\right) \Delta t\right) \psi _{j}^{\e ,\ast }.
\end{equation*}

As an example, consider an upwind
spatial discretization of $\mu $. 
In the first step, solve
\begin{equation}\label{ts1}
\left\{
\begin{aligned}
& i\e \partial _{t}\psi ^{\e }=-\frac{\e ^{2}}{2}\partial
_{xx}\psi ^{\e }, \\
&\frac{d}{dt}\mu _{jk}^{\e } =-\eta _{k}\left( D_{y}\mu
^{\e }\right) _{jk}-F_{j}^{\e }\left( D_{\eta }\mu ^{\e }\right)
_{jk},%
\end{aligned}
\right.
\end{equation}%
where both $D_{y}\mu ^{ \e}$ and $D_{\eta }\mu ^{\e }$ represent the
upwind discretization of the spatial derivatives.
To solve the Liouville equation, apply the forward Euler scheme for the time discretization.
Specifically,
\begin{equation} \label{upwind1}
\left\{
\begin{aligned}
& \psi _{j}^{\e ,\ast }=\frac{1}{{M}}\sum\limits_{\ell=-{{M}}/2}^{{M}/2-1}e^{-i\e
\omega _{\ell}^{2}/2}\hat{\psi}_{\ell}^{\e ,n}e^{i\omega _{\ell}
x_{j} },\quad  j=0,\dots ,{M}-1 , \\
& \frac{\mu _{jk}^{\e ,\ast }-\mu _{jk}^{\e ,n}}{\Delta t}=-\eta
_{k}\left( D_{y}\mu ^{ \e,n}\right) _{jk}-F_{j}^{\e ,n}\left(
D_{\eta }\mu ^{\e,n}\right) _{jk},
\end{aligned}%
\right.
\end{equation}%
where $w_\ell = 2 \pi \ell$.

The second step is then given by
\begin{equation}\label{ts2}
\left\{
\begin{aligned}
& i\e \partial _{t}\psi ^{\e }=\Upsilon _{d}^{\e }\left(
x,t\right) \psi ^{ \e}, \\
& \frac{d}{dt}\mu _{jk}^{\e }=0,
\end{aligned}%
\right.
\end{equation}%
where $\Upsilon_{d}^{\e }\left( x,t\right) $ is the quadrature approximation of $%
\Upsilon ^{\e }\left( x,t\right) $.
Thus, one explicitly gets
\begin{equation} \label{upwind2}
\begin{split}
\psi _{j}^{\e ,n+1}=\exp \left(-i\Upsilon _{d}^{\e ,\ast }\left(
x_{j}\right) {{ \Delta t}} /\e \right) \psi _{j}^{\e,\ast } , \quad  \mu _{jk}^{\e ,n+1}=\mu _{jk}^{\e,\ast }%
\end{split}%
\end{equation}%
where
\begin{align*}
\Upsilon _{d}^{\e ,\ast }\left( x\right)
 =\sum\limits_{j=0}^{J-1}\sum\limits_{k=0}^{K-1}V\left( x,y_{j}\right) \mu
_{jk}^{\e,\ast }\Delta y\Delta \eta
 =\sum\limits_{j=0}^{J-1}\sum\limits_{k=0}^{K-1}V\left( x,y_{j}\right) \mu
_{jk}^{\e,n+1}\Delta y\Delta \eta ,
\end{align*}
which is the trapezoidal rule for $\mu V$ with compact support.

\subsubsection{The spatial meshing strategy}\label{sec:mesh}

 We first  show that
one can take the limit $\e \to 0$, for fixed $\Delta y$ and $\Delta \eta$. 
Consider a semi-discretized version of the SLE system \eqref{SL-1} in one spatial dimension $d=n=1$ where
the Liouville equation is discretized by the  upwind scheme:
\begin{equation}\label{quad-sle}
\left\{
\begin{aligned}
& i\e \partial _{t}\psi ^{\e }=-\frac{\e ^{2}}{2}\partial
_{xx}\psi ^{\e }+\Upsilon _{d}^{\e}\left( x,t\right) \psi ^{\e
},\quad \psi ^{\e }(0,x)=\psi _{\rm in}^{\e}(x) ,   \\
& \partial _{t}\mu ^{\e }+\eta D_{y}\mu ^{\e}+F^{\e}\left(
y,t\right) D_{\eta }\mu ^{\e}=0,\quad  \mu ^{\e}(0,y, \eta)=\mu
_{\rm in}^{\e} (y,\eta).
\end{aligned}
\right.
\end{equation}
The following theorem is given in \cite{FJS18}:

\begin{theorem}
Under some suitable conditions for $V$ and initial data, for any $T>0$, the solution of semi-discretized SLE system \eqref{quad-sle} satisfies, up to extraction of sub-sequences,
\begin{equation*}
w^{\e}[ \psi ^{\e }]  \stackrel{\e\rightarrow 0_+
}{\longrightarrow}
\nu ,\quad \mu _{jk}^{\e } \stackrel{\e\rightarrow 0_+
}{\longrightarrow} \mu
_{jk}^{0},
\end{equation*}
in ${\rm w}\text{--}\ast$ topology,
where $j=0,\cdots ,J-1$ and $k=0,\cdots ,K-1$. In addition, $\nu $ and $\mu _{jk}$
solve the semi-discretized Liouville-system
\begin{equation*}
\left\{
\begin{aligned}
&\, \partial _{t}\nu +\xi \partial _{x}\nu -\partial _{x}\Upsilon _{d}^{0}\left(
x,t\right) \partial _{\xi }\nu =0, \\
&\, \frac{d}{dt}\mu _{jk}^{0}+ \eta _{k}D_{y}\mu _{jk}^{0} + F_{j}^{0}D_{\eta }\mu _{jk}^{0}=0.
\end{aligned}
\right.
\end{equation*}
\end{theorem}

The above result shows that the scheme is AP in $y, \eta$ with respect to $\e$, namely one can use $\Delta y, \Delta \eta \sim O(1)$. This is the {\it first} such result for highly oscillatory problem in spatial variables, and  more interestingly, the problem under study is {\it nonlinear}!

\begin{remark}
Numerical experiments show that the same type of behavior is true not only for mixed spectral-finite difference schemes, but also for
purely spectral schemes, see \cite{FJS18}. The proof, however, only works for the former case since it requires positivity of the energy. For spectral method the theory is still lacking.
\end{remark}

%%%%%%%%%%%%%%%%%%%%

\subsubsection{Time-discretization}\label{sec:time}

The time-discretization of the splitting scheme can also be shown to be AP.  
Notice the semiclassical limit of SLE (\ref{SL-1}), as $\e\to 0$, is \cite{JSZ17}:
\begin{equation}\label{liou liou 1}
\partial _{t}\mu+\eta \cdot \nabla _{y}\mu+F^{0}\left( y,t\right)
\cdot \nabla _{\eta }\mu =0,
\end{equation}
\begin{equation}\label{liou liou 2}
\partial _{t}\nu+\xi \cdot \nabla _{x}\nu-\nabla _{x}\Upsilon ^{0}\left(
x,t\right) \cdot \nabla _{\xi }\nu=0.
\end{equation}
As $\e\to 0$, the splitting schemes (\ref{ts1}) and (\ref{ts2}) approach respectively 
\begin{equation} \label{macro1}
\left\{
\begin{aligned}
&\partial _{t}\nu +\xi \partial _{x}\nu =0, \\
&\frac{d}{dt}\mu_{jk}+\eta _{k}\left( D_{y}\mu \right) _{jk}+F_{j}^{0}\left( D_{\eta }\mu\right) _{jk}=0,
\end{aligned}
\right.
\end{equation}
and
\begin{equation} \label{macro2}
\left\{
\begin{aligned}
&\partial _{t}\nu -\partial _{x}\Upsilon _{d}^{0}\left( x,t\right) \partial
_{\xi }\nu =0, \\
& \frac{d}{dt}\mu _{jk}=0.
\end{aligned}
\right.
\end{equation}
This is the time splitting scheme for (\ref{liou liou 1})-(\ref{liou liou 2}), which $\nu$ is the limit of the Wigner transform of $\psi^\e$ on $x$ variable. 
This shows that $\Delta t\sim O(1)$ can be chosen independent of the small parameter $\e$. In turn, 
this yields the convergence of the  scheme towards the corresponding scheme of the limiting equation, as stated in \eqref{macro1} and \eqref{macro2}, uniformly in $\Delta t$. Hence it is AP in $t$.

In summary, the scheme (\ref{ts1})-(\ref{ts2}) is AP in $t, y, \eta$ with respect to $\e$. One only needs $\Delta x=O(\e)$.

\subsubsection{Numerical experiments}

We now present some numerical experiments from \cite{FJS18}. 
The interaction potential is given by 
\[
V\left( x,y\right) =\frac{%
\left( x+y\right) ^{2}}{2}.
\] 
The one-dimensional SLE system is solved  on the interval $x\in %
\left[ -\pi ,\pi \right] $ and $y,\eta \in \left[ -2\pi ,2\pi \right] $ with
periodic boundary conditions.

\begin{example}  \label{ex1n}
The initial conditions for the SLE system \eqref{SL-1} is:
\begin{equation*}
\psi _{\rm in}\left( x\right) =\exp\left({-25\left( x+0.2\right) ^{2}}\right) \exp\left({\frac{-i\ln
\left( 2\cosh \left( 5\left( x+0.2\right) \right) \right) }{5\e }}\right),
\end{equation*}
and
\[
\ \mu
_{\rm in}\left( y,\eta \right) =\left\{
\begin{array}{cc}
C_{\rm N}\, \exp\left({-\frac{1}{1-y^{2}}}\right) \exp\left({-\frac{1}{1-\eta ^{2}}}\right), & \  \text{for $\left\vert y\right\vert <1$, $\left\vert \eta \right\vert <1$} \\
0, & \text{otherwise}.%
\end{array}
\right.
\]
Here, $C_{\rm N}>0$ is the normalization factor such that $\iint_{%
\mathbb{R}^{2}}\mu _{in} dy\, d\eta $=1.
The time-splitting method with spectral-upwind scheme (i.e., with an upwind scheme for the Liouville's equation) is used. For $\e =$ $\frac{1}{256},\frac{1}{1024},%
\frac{1}{4096}$,   $T=0.5$, choose $\Delta x=\frac{%
2\pi \e }{16}$, $\Delta y=\Delta \eta =\frac{4\pi }{128}$.
For each choice of $%
\e$,  the SLE system is solved first with
$\Delta t\text{ independent of }h$ and, second, with $
\Delta t=o\left( \e \right)$,   specifically, we compare the two cases where $\Delta t=0.01$ and $\Delta t=%
\frac{\e }{10}$, for 
 numerical values of $\mu$  (denoted as $\mu_1$ and $\mu_2$, respectively). As shown in Table \ref{table1}, the error is insensitive in $\e$, showing a uniform in $\e$ convergence in $\Delta t, \Delta \eta$ and $\Delta y$. 

\begin{table}[ht]
\centering
%\setcaptionwidth{4.6in}

\begin{tabular}{|c|c|c|c|}
\hline
$\e $ & $1/256$ & $1/1024$ & $1/4096$ \\ \hline
$\frac{\left\Vert \mu_1-\mu_2 \right\Vert_{\ell^2}}{\left\Vert \mu_2 \right\Vert_{\ell^2}} $ & 1.65e-03 & 1.69e-03 & 1.70e-03 \\ \hline
\end{tabular}%
\caption{Example \ref{ex1n}. The relative $\ell^2-$difference (defined as $\frac{\left\Vert \mu_1-\mu_2 \right\Vert_{\ell^2}}{\left\Vert \mu_2 \right\Vert_{\ell^2}} $) for various $\e$.}
\label{table1}
\end{table}
\end{example}

\begin{figure}[tbph]
%\setcaptionwidth{4.6in}
\begin{center}
\includegraphics[scale=0.6]
{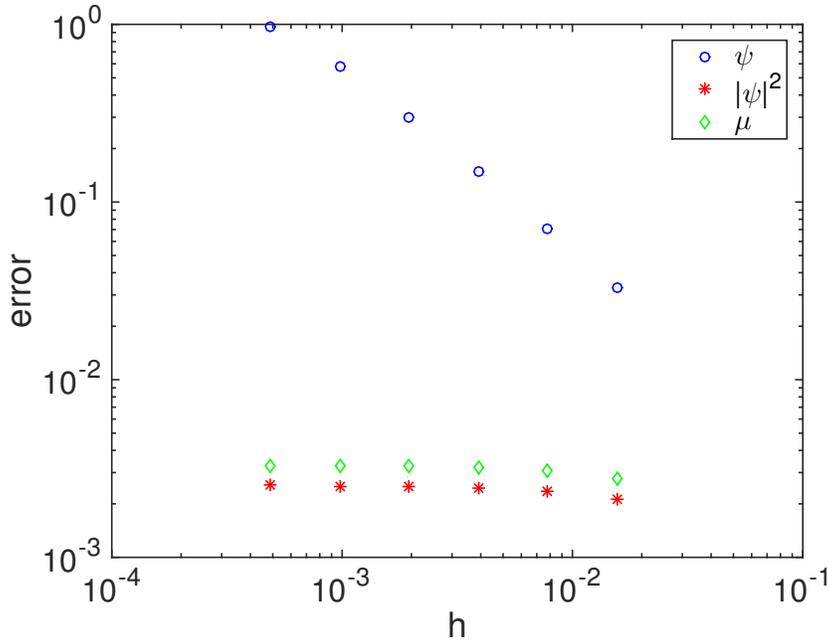}
\end{center}
\caption{Example \ref{ex2}: $\ell^2-$errors of the wave function $\psi^\e$, position density $|\psi^\e|^2$ and $\mu$ for various $\e$. Fix $\Delta t=$ $0.01$. For $h =\frac{1}{64},\frac{1}{%
128},\frac{1}{256},\frac{1}{512},\frac{1}{1024},\frac{1}{2048}$, choose $%
\Delta x=\frac{2\pi \e }{16}$ respectively. The
reference solution is computed with $\Delta t=\frac{ \e}{10}$. }
\label{delta conv}
\end{figure}

\begin{example} \label{ex2}
In this example, we choose the same initial data for $\mu_{\rm in}$ as in Example \ref{ex1n} but
\begin{equation*}
\psi _{\rm in}\left( x\right) =\exp\left({-5\left( x+0.1\right) ^{2}}\right)\exp \left({\frac{i\sin x}{%
h }}\right).
\end{equation*}%
Now, fix $\Delta t=$ $0.01$, a stopping time $T=0.4$, and $\Delta y=\Delta \eta =%
\frac{4\pi }{128}$, while $\Delta x=\frac{2\pi
\e }{16}$, for $\e =\frac{1}{64},\frac{1}{128},\frac{1}{256},%
\frac{1}{512},\frac{1}{1024},\frac{1}{2048}$, respectively. The reference solution is computed with $\Delta
t=\frac{h }{10}$. From the $\ell^2$-error plotted in Figure \ref{delta conv},
one can see that although the error in the wave function increases as $%
\e $ decreases, the error for the position density $| \psi^\e |
^{2}$ as well as for the macroscopic quantity $\mu $ does not change noticeably. This
shows that $\e-$independent $\Delta t, \Delta y$ and $\Delta \eta$ can be taken to accurately obtain
physical observables, but not the wave function  $\psi^\e$ itself.
\end{example}

\section{Numerical passages from classical mechanics to kinetic equations}\label{sec3}
\setcounter{equation}{0}

\subsection{The Random Batch Methods}\label{sec:alg}

Consider the second order interacting particle systems described by
\begin{gather}\label{eq:Nbody2nd}
\begin{split}
& d{r}_i={v}_i\,dt,\\
& d{v}_i=\Big[ b({r}_i)+\alpha_N \sum_{j:j\neq i}K({r}_i-{r}_j)-\gamma {v}_i \Big]\,dt+\sigma\,
d{W}_i,
\end{split}
\end{gather}
and the first order system, 
\begin{gather}\label{eq:Nparticlesys-2}
d{r}_i=b({r}_i)\,dt+\alpha_N\sum_{j: j\neq i} K_1({r}_i-{r}_j)\,dt+\sigma\, d{W}_i,~~i=1,2,\cdots, N.
\end{gather}
The main difficulty for the numerical simulations of particle
system  \eqref{eq:Nbody2nd} or \eqref{eq:Nparticlesys-2}  is that for lager $N$,
 the computational cost per time step is $\mathcal{O}(N^2)$. The Fast Multipole Method (FMM) \cite{rokhlin1985rapid}  reduces the complexity to $\mathcal{O}(N)$ if the interaction decays sufficiently fast. However, the implementation of FMM is quite delicate. A simple random algorithm, called the Random Batch Method (RBM), has been proposed in \cite{jin2020random}  to reduce the computation cost per time step from $\mathcal{O}(N^2)$ to $\mathcal{O}(N)$. The key idea of RBM is to use randomly chosen ``mini-batch'' in the summation term in \eqref{eq:Nparticlesys-2}  and \eqref{eq:Nbody2nd}.  Such an idea has its origin in the  stochastic gradient descent (SGD) method.
 The idea was also used for the computation of the mean-field flocking model \cite{albi2013}\cite{carrillo2017particle}.

Let $T>0$ be the simulation time, and one chooses a time step $\Delta t>0$. Pick a batch size $2\le p\ll N$ that divides $N$. Consider the discrete time grids $t_k:=k\Delta t$, $k\in \mathbb{N}$. For each sub-interval $[t_{k-1}, t_k)$, the method has two sub-steps: (1) at $t_{k-1}$, randomly group  the $N$ particles into $n:=N/p$ sub-groups (batches); (2) particles only interact with those in the same batch. This is given in 
Algorithm  \ref{alg:rbm2nd1}.
\begin{algorithm}[H]
\caption{(RBM for \eqref{eq:Nbody2nd})}
\label{alg:rbm2nd1}
\begin{algorithmic}[1]
\For {$m \text{ in } 1: [T/\Delta t]$}   
\State Divide $\{1, 2, \ldots, N=pn\}$ into $n$ batches $\mathcal{C}_q, 1\le q\le n$ randomly.
     \For {each batch  $\mathcal{C}_q$} 
     \State Update ${r}_i, {v}_i$ ($i\in \mathcal{C}_q$) by solving  for $t\in [t_{m-1}, t_m)$ the following
     \begin{gather}\label{eq:RBM2nd}
     \begin{split}
            & d{r}_i={v}_i\,dt,\\
            & d{v}_i=\Big[b({r}_i)+\frac{\alpha_N(N-1)}{p-1}\sum_{j\in\mathcal{C}_q,j\neq i}K({r}_i-{r}_j) -\gamma {v}_i\Big]\,dt+\sigma  d{W}_i.
      \end{split}
      \end{gather}
      \EndFor
 \EndFor
\end{algorithmic}
\end{algorithm}

RBM uses the random permutation, and each particle belongs to one and only one batch. 
An alternative approach, which allows replacement, is the following algorithm: 
\begin{algorithm}[H]
\caption{(RBM-r)}\label{randomreplacement}
\begin{algorithmic}[1]
\For {$m \text{ in } 1: [T/\Delta t]$}   
     \For {$k$ from $1$ to $N/p$}
     \State Pick a set $\mathcal{C}_k$ of size $p$ randomly with replacement. 
     \State Update ${r}_i$'s ($i\in \mathcal{C}_k$) by solving the following SDE for time $\Delta t$.
     \begin{gather}\label{eq:algorithmreplacement}
     \left\{
           \begin{split}
           & d{x}_i={u}_i\,dt,\\
            & d{u}_i=\Big[b({x}_i)+\frac{\alpha_N(N-1)}{p-1}\sum_{j\in\mathcal{C}_k,j\neq i}K({x}_i-{x}_j) -\gamma {u}_i\Big]\,dt+\sigma\, d{W}_i. \\
            & {x}_i(0) ={r}_i, \quad {u}_i(0)={v}_i,
            \end{split}
            \right.
      \end{gather}
 \quad\quad\quad i.e., solve \eqref{eq:algorithmreplacement} with initial values ${x}_i(0) ={r}_i, {u}_i(0)={v}_i$, and set ${r}_i\leftarrow {x}_i(\Delta t)$, ${v}_i\leftarrow {u}_i(\Delta t)$.
      \EndFor
\EndFor
\end{algorithmic}
\end{algorithm}

Different from Algorithm 1, in Algorithm 2,  for one iteration of $k$, some particles may not be updated while some may be drawn more than once. 

The random division into $n$ batches of equal size can be implemented using random permutation, which can be realized in $O(N)$ operations by Durstenfeld's modern revision of Fisher-Yates shuffle algorithm \cite{durstenfeld1964} (in MATLAB, one can use ``randperm(N)''). 
The ODE solver per particle per time step in \eqref{eq:RBM2nd} or (\ref{eq:algorithmreplacement}) requires merely $O(p)$ operations, thus for all particles, each time step costs only $O(pN)$. Since $p\ll N$ the overall cost per time step is significantly reduced from  $O(N^2)$ to basically $O(N)$. 

For RBM to really gain significant efficiency, one needs $\Delta t$ to be {\it independent} of $N$. We state an error estimate on RBM
 for the second order systems \eqref{eq:Nbody2nd} in the mean field regime (i.e., $\alpha_N=1/(N-1)$) from \cite{jin2020rbm2nd}, which was built upon the argument for the first order system in \cite{jin2020random}. 
 
 Denote $(\tilde{{r}}_i, \tilde{{v}}_i)$ the solutions to the random batch process \eqref{eq:RBM2nd} with the Brownian motion used being $\tilde{{W}}_i$. Consider the synchronization coupling:
\begin{gather}\label{eq:coupling}
{r}_i(0)=\tilde{{r}}_i(0) \sim \mu_0,~~{W}_i=\tilde{{W}}_i.
\end{gather}
 Let $\E$ denote the expected value, namely integration on $\Omega$ with respect to the probability measure $\bbP$, and consider the $L^2(\cdot)$ norm of a random variable
\begin{gather}
\|\zeta\|=\sqrt{\mathbb{E}|\zeta|^2}.
\end{gather}

For finite time interval, the error  of RBM is given by the following theorem.
\begin{theorem}\label{thm:convfinitetime}
Let $b(\cdot)$ be Lipschitz continuous, and  assume that $|\nabla^2 b|$ has polynomial growth,  and the interaction kernel $K$ is Lipschitz continuous. Then,
\begin{gather}
\sup_{t\in [0, T]}\sqrt{\E|\tilde{{r}}_i(t)-{r}_i(t)|^2+\E|\tilde{{v}}_i(t)-{v}_i(t)|^2}
\le C(T)\sqrt{\frac{\Delta t}{p-1}+(\Delta t)^2},
\end{gather}
where $C(T)$ is independent of $N$.
\end{theorem}

RBM has also been proposed for interacting particle systems used as a sampling method for the invariant measure of \eqref{eq:Nbody2nd} \cite{li2020stochastic}\cite{li2020random}\cite{J-L-QMC}. In these applications, the long-time behavior, and in particular
the convergence to the invariant measure, is of interest. 
 For such analysis   some additional contraction  assumptions are needed:
\begin{assumption}\label{ass:convexity}
$b=-\nabla V$ for some $V\in C^2(\mathbb{R}^d)$ that is bounded from below (i.e., $\inf_x V(x)>-\infty$), and there exist $\lambda_M\ge \lambda_m>0$ such that the eigenvalues of $H:=\nabla^2V$ satisfy
\[
\lambda_m \le \lambda_i(x)\le \lambda_M,~\forall~1\le i\le d, x\in \mathbb{R}^d.
\]
The interaction kernel $K$ is bounded and Lipschitz continuous. Moreover, the friction $\gamma$ and the Lipschitz constant $L$ of $K(\cdot)$ satisfy
\begin{gather}
 \gamma>\sqrt{\lambda_M+2L},~~\lambda_m>2L.
\end{gather}
\end{assumption}
Then the following uniform strong convergence estimate holds \cite{jin2020rbm2nd}:
\begin{theorem}\label{thm:longtimeconv}
Under Assumption \ref{ass:convexity} and the coupling \eqref{eq:coupling}, 
the solutions to \eqref{eq:Nbody2nd} and \eqref{eq:RBM2nd} satisfy
\begin{gather}\label{RBM-error}
\sup_{t\ge 0}\sqrt{\E|\tilde{{r}}_i(t)-{r}_i(t)|^2+\E|\tilde{{v}}_i(t)-{v}_i(t)|^2}
\le C\sqrt{\frac{\Delta t}{p-1}+(\Delta t)^2},
\end{gather}
where the constant $C$ does not depend on $ p$ and $N$.
\end{theorem}

Since the errors in {RBM-error} and (\ref{RBM-error}) are independent of $N$, in the mean-field regime,  thus RBM is AP in  particle number $N$ in this regime.

\subsubsection{An illustrating example: Dyson Browinan motion}

The following example is from \cite{jin2020random}. Consider a typical example in random matrix theory,  
where one is interested in solving the following  system of SDEs ($1\le j\le N$), called the Dyson Brownian motion:
\begin{gather}
d\lambda_j(t) =-\beta \lambda_j(t)\,dt+\frac{1}{N}\sum_{k: k\neq j}\frac{1}{\lambda_j-\lambda_k}dt
+\frac{1}{\sqrt{N}} dW_j,
\end{gather}
where $\{W_j\}$'s are independent standard Brownian motions. The system can be used to find the eigenvalues of a Hermitian-valued Ornstein-Uhlenbeck process.  The Brownian motion effect is small when $N$ is large.  The limiting equation for $N\to\infty$ is given by 
\begin{gather}\label{eq:dysonlimiteq}
\partial_t\rho(x,t)+\partial_x(\rho(u-\beta x))=0, ~~u(x, t)=\pi(H\rho)(x, t),
\end{gather}
where $\rho$ is the density for $\lambda$ as $N\to\infty$, $H(\cdot)$ is the Hilbert transform on $\mathbb{R}$, and $\pi=3.14\ldots$ is the circumference ratio.

For
$\beta=1$,
it can be shown that the corresponding limiting equation \eqref{eq:dysonlimiteq} has an invariant measure, given by the semicircle law:
\begin{gather}\label{eq:semicircle}
\rho(x)=\frac{1}{\pi}\sqrt{2-x^2}
\end{gather}

\begin{figure}
\begin{center}
	\includegraphics[width=1\textwidth]{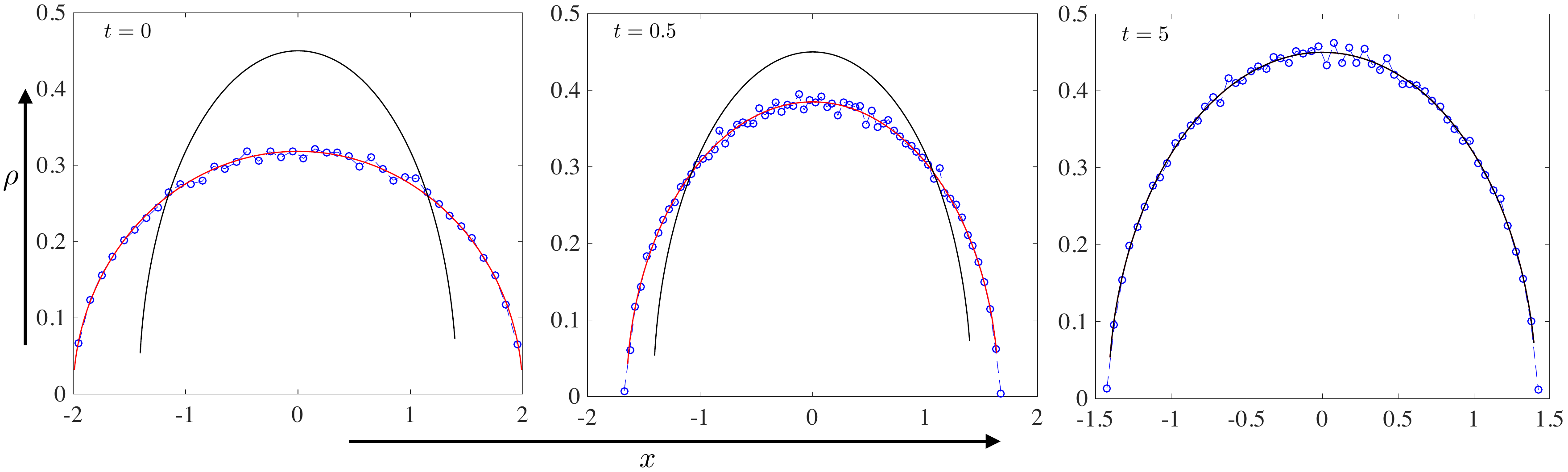}
\end{center}
\caption{The RBM solution (circles) of the Dyson Brownian motion. The empirical densities at various times are plotted. The red curve is the density distribution predicted by the analytic solution \eqref{eq:analyticaldis}. The black curve is the equilibrium semicircle law \eqref{eq:semicircle}.}
\label{fig:dysonbmDivision}
\end{figure}

To numerically test the behavior of RBM,  note an analytic solution to the limiting equation \eqref{eq:dysonlimiteq}
\begin{gather}\label{eq:analyticaldis}
\rho(x, t)=\frac{\sqrt{2\sigma(t)-x^2}}{\sigma(t)\pi},~~\sigma(t)=1+e^{-2t}.
\end{gather}
 For each iteration,  the force is singular, a  splitting strategy is adopted. Specifically, 
 define
\begin{gather}
X^{ij}:=X^i-X^j.
\end{gather}
The  RBM is implemented as follows:
\begin{itemize}
\item  \[
 \begin{split}
        Y_m^i=\frac{1}{2}(X_{m-1}^i+X_{m-1}^j)
                      +\mathrm{sgn}(X_{m-1}^{ij})\sqrt{|X_{m-1}^{ij}|^2+4\Delta t},\\
        Y_m^j=\frac{1}{2}(X_{m-1}^i+X_{m-1}^j)
                      -\mathrm{sgn}(X_{m-1}^{ij})\sqrt{|X_{m-1}^{ij}|^2+4\Delta t}.
\end{split}
\]
\item
\[
     X_m^i=Y_{m}^i-\Delta t  Y_{m}^i+\sqrt{\frac{\Delta t}{N}}z^i,~~
       X_m^j=Y^j(t_m)-\Delta t  Y_m^j+\sqrt{\frac{\Delta t}{N}}z^j.
\]
\end{itemize}
Here, $z^i,z^j\sim \mathcal{N}(0, 1)$.

Fig. \ref{fig:dysonbmDivision} shows that   RBM captures the evolution of distribution and the equilibrium semicircle law \eqref{eq:semicircle}, as desired. RBM-r also has similar behavior.

\subsubsection{The mean-field limit of RBM}

To further understanding the behavior of RBM, when $N $ is large, it will be interesting to investigate
its mean field limit. To this aim, consider RBM for the first order system \eqref{eq:Nparticlesys-2}
 with $\alpha_N=1/(N-1)$.  

Intuitively, when $N\gg 1$, the probability that two chosen particles are correlated is very small. Hence, in the $N\to\infty$ limit, two chosen particles will be independent with probability $1$. Due to the exchangeability, the marginal distributions of the particles will be identical.  Based on this observation, the following mean field limit   was derived and proved in  \cite{jin2021mean}:

\begin{algorithm}[H]
\caption{(Mean Field Dynamics of RBM  for first order system \eqref{eq:Nparticlesys-2})}\label{meanfield}
\begin{algorithmic}[1]
\State $\tilde{\mu}(\cdot, t_0)=\mu_0$.
\For {$k \ge 0$}  

\State Let $\rho^{(p)}(\cdots, 0)=\tilde{\mu}(\cdot, t_{k})^{\otimes p}$ be a probability measure on $(\mathbb{R}^{d})^{ p}\cong \mathbb{R}^{pd}$.

\State Evolve the measure $\rho^{(p)}$ to find $\rho^{(p)}(\cdots, \Delta t)$ by the following Fokker-Planck equation:
\begin{gather}\label{eq:firstalgorithm}
            \partial_t\rho^{(p)}=-\sum_{i=1}^p 
            \nabla_{x_i}\cdot\left(\Big[b(x_i)+\frac{1}{p-1}\sum_{j=1,j\neq i}^p K_1(x_i-x_j)\Big]\rho^{(p)}\right)+\frac{1}{2}\sigma^2\sum_{i=1}^p \Delta_{x_i}\rho^{(p)}.
\end{gather}

\State Set
\begin{gather}
\tilde{\mu}(\cdot, t_{k+1}):=\int_{(\mathbb{R}^{d})^{\otimes(p-1)}}
\rho^{(p)}(\cdot,dy_2,\cdots,dy_p, \Delta t).
\end{gather}

\EndFor
\end{algorithmic}
\end{algorithm}

The dynamics in Algorithm \ref{meanfield} naturally gives a nonlinear operator $\mathcal{G}_{\infty}: \mathbf{P}(\mathbb{R}^d)\to \mathbf{P}(\mathbb{R}^d)$ as
\begin{gather}\label{eq:Ginfty}
\tilde{\mu}(\cdot, t_{k+1})=: \mathcal{G}_{\infty}(\tilde{\mu}(\cdot, t_k)).
\end{gather}
Corresponding to this is the following SDE system for $t\in [t_k, t_{k+1})$
\begin{gather}\label{eq:RBMmeanfieldSDE}
dm{x}_i=b(m{x}_i)\,dt+\frac{1}{p-1}\sum_{j=1,j\neq i}^{p}K_1(m{x}_i-m{x}_j)\,dt
+\sigma\,dm{W}_i,~~i=1,\cdots, p,
\end{gather}
with $\{m{x}_i(t_k)\}$ drawn i.i.d from $\tilde{\mu}(\cdot, t_k)$.  

Hence, in the mean field limit of RBM, one starts with a  configuration in molecular chaos, then the $p$ particles evolve by interacting with  each other. One takes the first marginal of this new $p$-particle distribution, and  at the starting point of the next time interval, one imposes the  molecular chaos condition so that the particles are independent again. 

Furthermore, in \cite{jin2021mean} it was proven that this mean-field limit is $O(\Delta t)$ distance (in Wasserstein-1 sense) to   the mean-field limit of the original particle
system \eqref{eq:eq:limitvlasov-1st}, thus completing the AP diagram  in Fig. \ref{fig-1}.

\subsection{Molecular dynamics}\label{sec:md}

One of the most important interacting particle systems is molecular dynamics (MD), which simulates  dynamics or equilibrium properties  of large system of atoms and molecules using Newton's second law.  It has wide range of applications,  such as chemical physics, soft materials and biophysics \cite{ciccotti1987simulation,frenkel2001understanding}.  Here we review an interesting application of RBM  to MD simulation, called Random Batch Ewald \cite{jinlixuzhao2020rbe}, which achieves an $O(N)$ complexity  with a high parallel efficiency \cite{liangetalRBE}.

The equations of motion governing  $N$ ``molecules'' with  masses $m_i$'s are given by
\begin{gather}\label{eq:md1}
\begin{split}
& d{r}_i={v}_i\,dt,\\
& m_i d{v}_i=\Big[-\sum_{j:j\neq i}\nabla \phi({r}_i-{r}_j)\Big]\,dt
+\sigma_i d{W}_i.
\end{split}
\end{gather}
Here,  ${W}_i$ are noise or other external forcing terms, $\phi(\cdot)$ is the Coulomb potential
\[
\phi(x) = \frac{q_i q_j}{r},
\]
where $q_i$ is the charge for the $i$th particle and $r=|x|$. Another popular potential often used is the Lennard-Jones potential \cite{frenkel2001understanding}:
\begin{equation}\label{LJ-p}
%\label {L-J}
\phi(x)=4\left(\frac{1}{r^{12}}-\frac{1}{r^6}\right).
\end{equation}
Between ions, both types of potential exist and between charge-neutral molecules, the Lennard-Jones potential might be the main force.

\subsubsection{RBM with kernel splitting}\label{subsec:splitting}

Due to the singularity at  $x=0$ of the 
 Lenard-Jones potential (\ref{LJ-p}),  a direct application of RBM could give poor results. One effective strategy is to  decompose the $K$ into two parts \cite{martin1998novel,hetenyi2002multiple}, 
\begin{gather}
K(x)=K_1(x)+K_2(x).
\end{gather}
Here, $K_1$ has short range that decays quickly hence can be ignored  for $|x|\ge r_0$, for some $r_0$  chosen to be comparable to the mean distance of the particles. $K_2(x)$ is a bounded smooth function.
One then applies RBM to the $K_2$ part only \cite{li2020random}.

\subsubsection{Random Batch Ewald: importance sampling}

The Coulomb interaction is a long range interaction, which decays slowly as $1/r$, and in the mean time contains a  singularity at $r=0$. The bottleneck in MD simulation lies in the expensive simulation of the Coulomb interaction,
which has the computational complexity of $O(N^2)$.  Some popular methods include  the particle-particle particle mesh Ewald (PPPM) \cite{LDT+:MS:1994,Deserno98JCP}, and multipole type methods such as treecode \cite{BH:N:1986,DK:JCP:2000} and fast multipole methods (FMM) \cite{GR:JCP:1987,YBZ:JCP:2004}. 
These methods can reduce the complexity per time step from $O(N^2)$ to $O(N\log N)$ or even $O(N)$, and have gained big success in practice. However, some issues still remain to be resolved, e.g., the prefactor in the linear scaling can be large, or the implementation can be nontrivial, or the scalability for parallel computing is not high.

 The RBE method is based on the Ewald splitting for the Coulomb kernel with a random ``mini-batch'' type technique applied in the Fourier series for the long-range part. 

The solids or fluids with large volume are usually modeled  in a box with length $L$,  with  periodic conditions. 
Consider $N$  particles
 with net charge $q_i$ ($1\le i\le N$) under the electroneutrality condition
\begin{gather}
\sum_{i=1}^N q_i=0.
\end{gather}
The forces are computed using ${F}_i=-\nabla_{{r}_i}U$, where $U$ is the Coulomb potential energy, with  periodic boundary condition, given by
\begin{gather}\label{eq:energy}
U=\frac{1}{2}\sum_{{n}}{}'\sum_{i,j=1}^N q_iq_j \frac{1}{|{r}_{ij}+{n}L|},
\end{gather}
where ${n}\in \mathbb{Z}^3$. $\sum'$ is defined such that ${n}=0$ is not included when $i=j$.

The classical Ewald summation decomposes $1/r$ into the  long-range smooth parts and short-range singular parts: 
\begin{gather}
\frac{1}{r}=\frac{{\text{erf}}(\sqrt{\alpha}r)}{r}+\frac{{\text{erfc}}(\sqrt{\alpha r})}{r},
\end{gather}
where ${\text{erf}}(x):=\frac{2}{\sqrt{\pi}}\int_0^x \exp(-u^2)du$ is the error function and $\text{erfc}=1-\text{erf}$. 
Correspondingly, $U=U_1+U_2$ with
\begin{gather}
U_1=\frac{1}{2}\sum_{{n}}{}'\sum_{i,j}q_iq_j\frac{\text{erf}(\sqrt{\alpha}|{r}_{ij}+{n}L|)}{|{r}_{ij}+{n}L|}, \\
~~U_2=\frac{1}{2}\sum_{{n}}{}'\sum_{i,j}q_iq_j\frac{\text{erfc}(\sqrt{\alpha}|{r}_{ij}+{n}L|)}{|{r}_{ij}+{n}L|},
\end{gather}
where  $U_{2}$ corresponds to the short-range forces which is inexpensive, while $U_1$ is the long-range part that will be put into  the Fourier space
\begin{gather}
U_1=\frac{2\pi}{V}\sum_{{k}\neq 0}\frac{1}{k^2}|\rho({k})|^2
e^{-k^2/4\alpha}-\sqrt{\frac{\alpha}{\pi}}\sum_{i=1}^N q_i^2,
\end{gather}
where  
$\rho({k}):=\sum_{i=1}^N q_i e^{i{k}\cdot{r}_i}$.
Then
\begin{gather}\label{eq:force}
{F}_{i,1}=-\nabla_{r_i} U_1=-\sum_{{k}\neq 0}\frac{4\pi q_i {k}}{V k^2}
e^{-k^2/(4\alpha)}\mathrm{Im}(e^{-i{k}\cdot{r}_i}\rho({k})),
\end{gather}
where ${r}_{ij}={r}_j-{r}_i$, 
 is bounded for small ${k}$. The key idea of RBE is to do {\it importance sampling} according to the discrete Gaussian distribution 
 $e^{-k^2/(4\alpha)}$.  Denote 
\begin{gather}\label{eq:S}
S:=\sum_{{k}\neq 0}e^{-k^2/(4\alpha)}=H^3-1,
\end{gather}
where
\begin{gather}
H:=\sum_{m\in mathbb{Z}}e^{-\pi^2 m^2/(\alpha L^2)}
=\sqrt{\dfrac{\alpha L^2}{\pi}}\sum\limits_{m=-\infty}^{\infty}e^{-\alpha m^2L^2}
\approx\sqrt{\frac{\alpha L^2}{\pi}}(1+2e^{-\alpha L^2}), \label{psf}
\end{gather}
since often $\alpha L^2 \gg 1$. Hence, $S$ is the sum for all three-dimensional vectors ${k}$ except $0$. 
Then, one can regard the sum as an expectation over the probability distribution
\begin{gather}\label{eq:probexpression}
\mathscr{P}_{{k}}:=S^{-1}e^{-k^2/(4\alpha)},
\end{gather}
which, with ${k}\neq 0$, is a discrete Gaussian distribution that can be sampled efficiently {\it offline}. Once the time evolution starts one just needs to randomly
draw a few ($O(p)$) samples for each time step from this pre-sampled Guassian sequence.

Ultimately, the force ${F}_{i,1}$ in \eqref{eq:force} will be calculated  by the following mini-batch random variable:
\begin{gather}\label{eq:rbmapprox}
{F}_{i,1}\approx {F}_{i,1}^*:=-\sum\limits_{\ell=1}^p \dfrac{S}{p}\dfrac{4\pi {k}_\ell q_i}{V k_\ell^2}\mathrm{Im}(e^{-i{k}_\ell\cdot{r}_i}\rho({k}_\ell)).
\end{gather}

 The PPPM uses FFT, while RBE uses random mini-batch to speed up the computation in the Fourier space. 
The complexity of RBE for the real space part is $\mathcal{O}(N)$. 
 By choosing {\it the same batch} of frequencies for all forces \eqref{eq:rbmapprox} (i.e., using the same
${k}_{\ell}$, $1\le \ell \le p$ for all ${F}^*_{i,1}, 1\le i\le N$) in the same time step, the complexity per iteration for the frequency part is reduced to $\mathcal{O}(pN)$. Therefore the RBE method has linear complexity per time step if one chooses $p=\mathcal{O}(1)$.

Another advantage of RBE is that there are few particle interactions at each iteration. This significantly reduces the amount of message passing when many CPUs are used for parallel computing,
hence one achieves remarkable scalability \cite{liangetalRBE}.

\begin{figure}[ht]	
	\centering
	\includegraphics[width=0.48\textwidth]{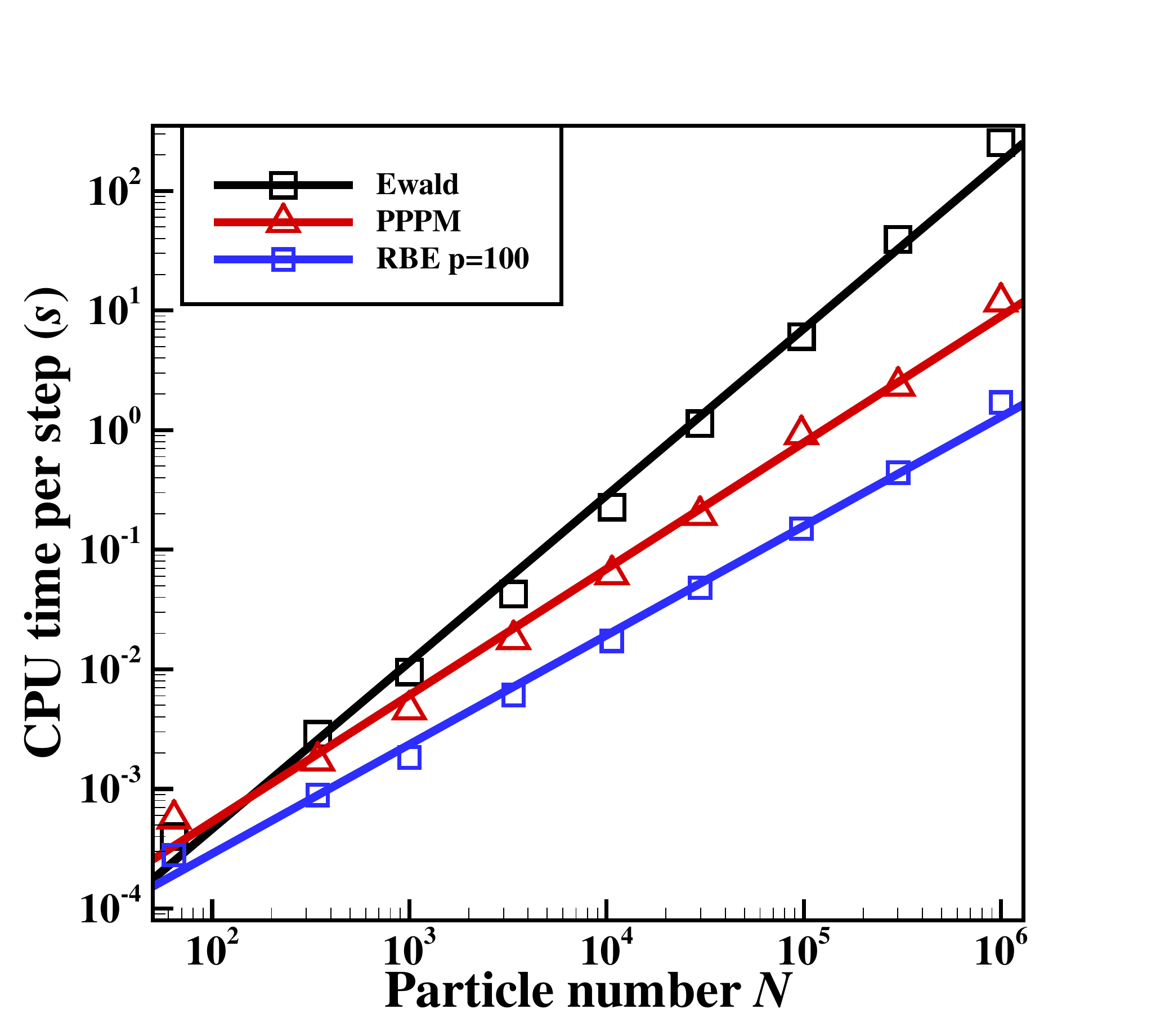}
	\caption{Comparision of the Ewald sum, the PPPM and the RBE methods}
	\label{fig:rbecompare}
\end{figure}

To illustrate the performance of the RBE method, consider an electrolyte with monovalent binary ions (first example in \cite{jinlixuzhao2020rbe}). In the reduced units (\cite[section 3.2]{frenkel2001understanding}), the dielectric constant is taken as $\varepsilon=1/4\pi$ so that the potential of a charge is $\phi(r)=q/r$ and the temperature is $T=\beta^{-1}=1$. Under the Debye-H\"uckel (DH) theory (linearized Poisson-Boltzmann equation), the charge potential outside one ion is given by
\[
-\varepsilon \Delta \phi= 
\begin{cases}
0  & r<a\\
q\rho_{\infty,+}e^{-\beta q\phi}-q\rho_{\infty,-}e^{\beta q\phi}\approx \beta q^2\rho_r \phi, & r>a
\end{cases}
\]
where $\rho_{\infty,+}=\rho_{\infty,-}=N/(2V)$ are the densities of the positive and negative ions at infinity, both being $\rho_r/2$. The parameter $a$ is the effective diameter of the ions, which is related to the setting of the Lennard-Jones potential. In the simulations, $a=0.2$ and the setting of Lennard-Jones potential can be found in \cite{jinlixuzhao2020rbe}.  This approximation gives the net charge density $\rho=-\varepsilon\Delta\phi$
for $r\gg a$,
\[
\ln(r\rho(r))\approx -1.941 r-1.144.
\]
 Fig. \ref{fig:rbecompare} shows the CPU time consumed for different particle numbers inside the box with the same side length $L=10$. Both the PPPM and RBE methods scale linearly with the particle numbers. However, even for batch size $p=100$, the RBE method consumes much less time. 
 Clearly, the RBE method has the same level of accuracy compared with the PPPM method for the densities considered.

\begin{figure}[ht]	
	\centering
	\includegraphics[width=0.48\textwidth]{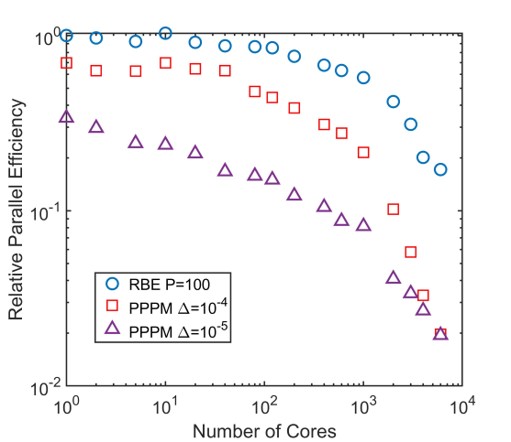}
	\includegraphics[width=0.48\textwidth]{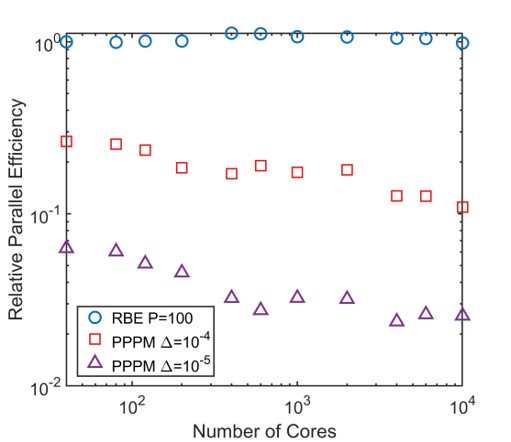}
	\caption{The parallel efficiency of the PPPM and the RBE methods for all-atom simulation of pure water system (Left)  $3\times 10^5$ atoms; (Right) $3\times 10^7$ atoms}
	\label{fig:paraefficiency}
\end{figure}

Next, in Fig. \ref{fig:paraefficiency}, the parallel efficiency of the PPPM and RBE methods from  \cite{liangetalRBE} for the all-atom simulation of pure water systems is shown.  As can be seen, due to the reduction of communications for the particles, the RBE method gains better parallel efficiency. This parallel efficiency is more obvious when the number of particles is larger.  In \cite{liangetalRBE}, the simulation results of pure water system also indicate that the RBE type methods can not only sample from the equilibrium distribution, but also compute accurately the dynamical properties of the pure water systems.

For a comprehensive review of RBM and its extensions and applications, see a recent review \cite{jin2021random}.

\section{Numerical passages from kinetic equations to hydrodynamic equations}\label{sec4}
\setcounter{equation}{1}

In nuclear reactor, neutrons may conduct significant amount of scattering, in the diffusive regime. 
In the space shuttle reentry problem, the mean free path  could vary from $O(1)$ meter in the space to $O(10^{-8})$ meter, when the vehicle
 passes from free streaming, rarefied gas (described by the Boltzmann
equation), transition to the 
hydrodynamic (described by the Euler or Navier-Stokes equations)
 regimes \cite{riv}.  It is also known that in hypersonic flows (Mach number larger than $1.4$), 
the shock profile of the Navier-Stokes equations  do not give accurate shock width, hence one needs to use  the Boltzmann equation in the shock region \cite{Foch}\cite{Aga}. In plasma physics, 
the Debye length could be small, and one needs to deal with quasi-neutral regime \cite{FO}\cite{degond-JCPReview}. In all these kinetic problems one needs to deal with multiple time and spatial scales.

Kinetic theory is the area in which the concept of AP was first introduced, and also most successfully and extensively used. Earlier effort in this direction concentrates on time-independent transport equations that have diffusive behavior \cite{LMM1}\cite{LMM2}. However for multiscale kinetic equations  the main challenges lie in time discretizations, due to the stiffness, non-locality and nonlinearity of the collision operators.  

The term "Asymptotic-Preserving" was first coined in \cite{Jin-AP}. An AP scheme typically possesses the following  key features  for multiscale kinetic equations:

\begin{itemize}
    \item Implicit time discretization that can be either explicitly or easily implemented: for example at least avoids complicated nonlinear algebraic system solvers such as  Newton's iteration;
    \item when the Knudsen number $\varepsilon \to 0$ the scheme for the kinetic equations automatically become a good scheme for the limiting hydrodynamic equations
\end{itemize}

This program is also related to the development of kinetic schemes for compressible Euler equations, which was based on discretizing a linear kinetic equation thanks to its linear convection, followed by a projection to the local Maxwellian \cite{Deshpande}\cite{perthame90}\cite{PreXu}. It is also relevant to the lattice Boltzmann approximation to incompressible Navier-Stokes equations \cite{ChenLS}\cite{ChenDoolen} \cite{QianLG}\cite{Luo1997}.  Relaxation 
schemes for nonlinear hyperbolic systems also share similar spirit
\cite{Jin-Xin}. Below we review a few representative AP schemes.

\subsection{The BGK-penalization method}

We are mainly interested in  dealing with the numerical difficulties when the Knudsen number $\varepsilon \ll 1$. The first challenge is numerical stiffness, which puts severe constraint on $\Delta t$.  In order to allow $\Delta t \gg\varepsilon$, one needs some implicit treatment for the nonlocal, nonlinear collision operator, which is numerically nontrivial. 

The {\it penalization method}, introduced by Filbet and Jin \cite{Filbet-Jin}, 
was the first AP scheme for the nonlinear Boltzmann equation that overcomes the stiffness issue of the collision operator.  
The idea is to penalize $\mathcal{Q}(f)$ by the BGK operator $\beta(\mathcal{M}-f)$:
\begin{equation} \label{penalization}
\partial_tf+v\cdot \nabla_x f=\underbrace{\frac{\mathcal{Q}(f)-\beta(\mathcal{M}-f)}{\varepsilon}}_{\text{less stiff, explicit}}+\underbrace{\frac{\beta(\mathcal{M}-f)}{\varepsilon}}_{\text{stiff, explicit}},
\end{equation}
where $\beta$ is some constant chosen properly to approximate the Frechet derivative of $\mathcal{Q}(f)$ around $\mathcal{M}$, such that terms in the first brace become less stiff or non-stiff and can be treated explicitly. The other part  is a BGK operator, which can be inverted {\it explicitly} \cite{Coron-Perthame}, thanks to the conservation properties of the collision terms on the right hand side of  \eqref{penalization}.

A first-order IMEX (Implicit-Explicit) discretization of (\ref{penalization}) can  be written as:
\begin{equation} \label{IMEXBoltz}
\frac{f^{n+1}-f^n}{\Delta t}+v\cdot \nabla_x f^n=\frac{\mathcal{Q}(f^n)-\beta(\mathcal{M}^n-f^n)}{\varepsilon}+\frac{\beta(\mathcal{M}^{n+1}-f^{n+1})}{\varepsilon}.
\end{equation}
 Taking the moments $\int_{\mathbb{R}^d} \cdot \,\phi(v)\,\rd{v}$, with $\phi(v)$ defined in \eqref{consv1} on both sides of (\ref{IMEXBoltz}),  and using the properties  (\ref{consv1}), one gets
\begin{equation}\label{moment-eqn}
  \frac{\langle f \rangle ^{n+1}-\langle f \rangle ^n}{\Delta t}+ \nabla_x \cdot \langle v f \rangle ^n=0\,,
\end{equation}
where $\langle \cdot \rangle = \int \phi(v) \cdot dv$ means the moments. 
 From \eqref{moment-eqn}  one can solve for the moments $\rho, u$ and $T$ at $t=t^{n+1}$, hence $\mathcal{M}^{n+1}$
 is obtained. Then $f^{n+1}$ can be obtained from (\ref{IMEXBoltz}) explicitly. Notice the entire process is {\it explicit}!
 
 In practice, $\beta$ can be roughly estimated as
\begin{equation*}
\beta=\sup_{v} |\mathcal{Q}^-(f)|,
\end{equation*}
where $\mathcal{Q}^-$ is the loss part of the collision operator defined such that $\mathcal{Q}(f)=\mathcal{Q}^+(f)-f\mathcal{Q}^-(f)$. $\beta$ can also be made time and spatially dependent for better numerical accuracy \cite{YanJin13}.

To capture the compressible Euler limit, a necessary condition is that, as $\e \to 0$, 
\[
f^n=\mathcal{M}^n, \quad {\hbox {for any}}\, n, \quad {\hbox {with}} \quad \Delta t, \Delta x \quad {\hbox {fixed}}.  
\]
It was formally shown in  \cite{Filbet-Jin} that
\[
 {\hbox {for}}\, \, \varepsilon \ll 1, \quad  {\hbox {if}}\,\,f^n=\mathcal{M}^n+O(\varepsilon),\,\, {\hbox {then}} \,\, f^{n+1}=\mathcal{M}^{n+1}+O(\varepsilon).
 \]
Numerical experiments in \cite{Filbet-Jin} shows that regardless of the initial condition $f^0$, there exists an integer $N>0$ such that 
\begin{equation} \label{weakAP}
f^n=\mathcal{M}^n+O(\varepsilon),  \quad \text{for any } n \geq N.
\end{equation}
Substituting (\ref{weakAP}) into (\ref{IMEXBoltz}) and taking the moments, one has
\begin{equation*} 
\frac{\langle f \rangle ^{n+1}-\langle f \rangle ^n}{\Delta t}+\nabla_x \cdot \int_{\mathbb{R}^d}v \phi(v) \,\mathcal{M}^n\,\rd{v}=O(\varepsilon), \quad \text{for any } n \geq N,
\end{equation*}
which is a consistent discretization to the limiting Euler system (\ref{Euler-limit}). This means the scheme is AP after an initial transient time. 
\begin{remark}
A possible way to remove the initial layer problem and hence achieve AP in one time step was suggested in \cite{YanJin13}, where the idea is to perform the penalization in two successive steps:
\begin{align*} 
\left\{
\begin{array}{l}
\displaystyle \frac{f^*-f^n}{\Delta t}+v\cdot \nabla_x f^n=\frac{\mathcal{Q}(f^n)-\beta(\mathcal{M}^n-f^n)}{\varepsilon}+\frac{\beta(\mathcal{M}^*-f^*)}{2\varepsilon},\\[8pt]
\displaystyle  \frac{f^{n+1}-f^*}{\Delta t}=\frac{\beta(\mathcal{M}^{n+1}-f^{n+1})}{2\varepsilon}.
\end{array}\right.
\end{align*}
\end{remark}

The idea of using a linear or simpler operator to penalize the nonlinear or complicated operator turns out to be a generic approach. For specific problems, one needs to seek appropriate penalization operator. For example,
consider the nonlinear Fokker-Planck-Landau equation,  whose collision operator is given by
\begin{equation} 
\mathcal{Q}(f)(v)=\nabla_v\cdot\int_{\mathbb{R}^{d}}A(v-v_*)\left[f(v_*)\nabla_vf(v)-f(v)\nabla_{v_*}f(v_*)\right]\,\rd{v_*},
\end{equation}
where $A$ is a semi-positive definite matrix. This equation is relevant in the study of Coulomb interactions. The diffusive nature of the collision operator introduces more stiffness. An explicit scheme would require $\Delta t=O(\varepsilon (\Delta v)^2)$, where $\Delta v$ is the mesh size in $v$,  which is even more restrictive than the Boltzmann collision operator. In \cite{JinYan11}  the following Fokker-Planck operator was proposed as a penalization:
\begin{equation*}
\mathcal{P}_{FP}(f)=\nabla_v\cdot \left(\mathcal{M}\nabla_v\left(\frac{f}{\mathcal{M}}\right)\right).
\end{equation*}

Similar approaches, with variant penalties, have been proposed for the quantum Boltzmann equation \cite{FHJ12}, the quantum Fokker-Planck-Landau equation \cite{HJY12}, and the multi-species Boltzmann equation \cite{JinLi13}.

Another AP scheme for the Boltzmann equation,  developed later 
in \cite{liu2016unified}, relies on the 
integral representation of the BGK model. The final form of the
scheme also ends up with a linear 
combination of the Boltzmann collision operator and the BGK operator, with a slightly different combination coefficients.

%\textcolor{red}{{\bf Jingwei:} I feel we may not have enough space to include these equations%... I will leave this part as it is for the moment.}
%%AP Monte-Carlo \cite{RLJ14}

\subsection{Exponential integration}

Another class of asymptotic preserving method is the {\it exponential integration method}. This method is based on a  reformulation of  the equation into an exponential form, with the Maxwellian embedded. This makes it  easier to capture the asymptotic limit and other   physical properties such as  positivity.

  For the space homogeneous  Boltzmann equation:
\begin{equation}\label{eqn:homo_Boltzmann}
\partial_t f = \frac{1}{\varepsilon}\mathcal{Q}(f)\,,
\end{equation}
Dimarco and Pareschi in~\cite{DP11} introduced the following reformulation:
\begin{equation}\label{eqn:reformulation}
\partial_t\left[(f-\mathcal{M})e^{\beta t/\varepsilon }\right] = \partial_t f e^{\beta t/\varepsilon} + \frac{\beta(f-\mathcal{M})}{\varepsilon}e^{\beta t/\varepsilon} = \frac{\mathcal{Q} - \beta(\mathcal{M}-f)}{\varepsilon} e^{\beta t/\varepsilon}\,.
\end{equation}
Here $\beta$ is an auxiliary parameter and as in the penalization method, $\beta(\mathcal{M}-f)$ is used to approximate the Frechet derivative of $\mathcal{Q}$.  $\beta$ is chosen to be the smallest value that preserves the positivity of $f$. 

Equation~\eqref{eqn:reformulation} is fully equivalent to the original problem~\eqref{eqn:homo_Boltzmann}. However, it updates the difference between $f$ and $\mathcal{M}$, and the exponential term $\exp(-\beta t/\e)$ removes the stiffness and forces numerically the convergence between $f$ and $\mathcal{M}$, an essential mechanism for the AP property.  It can be easily extended to  all explicit  Runge-Kutta methods which are not only  of high order  but also holds the  AP property automatically. 
%The scheme writes as:
%\begin{equation}
%\begin{cases}
%\displaystyle f^{(i)} =\left(1-e^{-c_i\lambda} -\sum_{j=1}^{i-1}a_{ij}\lambda %e^{\lambda(-c_i+c_j)}\right)\mathcal{M}+e^{-c_i\lambda}f^n + \sum_{j=1}^{i-1}a%_{ij}\lambda e^{\lambda(c_j-c_i)}\frac{\mathcal{P}^{(j)}}{\beta}\,,\\
%\displaystyle f^{n+1} =\left(1-e^{-\lambda} -\sum_{i}^\nu b_{i}\lambda e^{\lam%bda(-1+c_i)}\right)\mathcal{M}+e^{-\lambda}f^n + \sum^\nu_{i}b_{i}\lambda e^{\%lambda(c_i-1)}\frac{\mathcal{P}^{(i)}}{\beta}\,.
%\end{cases}
%\end{equation}
%Here $\lambda = \beta\Delta t/\varepsilon$. $f^{(i)}$ denotes the solution at %the $i$-th substage in the $n$-th time step of the Runge-Kutta method (we negl%ect $n$ in the notation $f^{(i)}$ here for conciseness). We also used the nota%tion $\mathcal{P}=\mathcal{Q}+\beta f$. 

The need to convect $\mathcal{M}$  makes it difficult to extend the scheme to the 
nonhomogenous case. 
%Generally speaking, for the nonhomogeneous Boltzmann equation with a flux term time accuracy %has to be sacrificed. A straightforward method, as suggested in~\cite{DP}, is to perform oper%ator splitting and compute:
%\begin{equation}
%\begin{cases}
%\displaystyle \partial_t f + v\cdot\nabla_x f = 0\,,\\
%\displaystyle \partial_t f = \frac{1}{\varepsilon}\mathcal{Q}(f)
%\end{cases}
%\end{equation}
%iteratively in each time step. The drawback is immediate. As commented in Remark~\ref{remark:%high_order}, the splitting automatically deteriorates the high order of accuracy, making the %Runge-Kutta in the inner loop a waste.
In ~\cite{li2014exponential}, Li and Pareschi use an evolving  Maxwellian function 
\emph{within} each time step. 
%This provides the scheme the most accurate equilibrium function, and allows on%e to combine the treatment for the collision and for the transport terms, and %therefore get rid of the splitting that restricts the time accuracy order. The% method is strongly AP, namely achieve $f^n=\mathcal{M}^n+O(\varepsilon)$ in o%ne time step (compare with (\ref{weakAP})), and can reach arbitrary high order% of accuracy in time. 
%%The high order spatial accuracy is provided by the standard WENO.
They reformulate the Boltzmann equation as 
%ng Maxwellian idea. Basically one needs to rewrite the equation entirely keepi%ng in mind that the Maxwellian evolves together with the distribution function%. The reformulation reads:
\begin{equation}\label{eqn:exp_reform}
\partial_t\left[ (f-\mathcal{M})\exp{(\beta t/\varepsilon)}\right] = \left(\frac{\mathcal{P}-\beta\mathcal{M}}{\varepsilon} - v\cdot\nabla_xf -\partial_t \mathcal{M}\right)\exp(\beta t/\varepsilon)\,,
\end{equation}
while the moment equations are obtained after taking the moments of the original Boltzmann equation (\ref{CBE}):
\begin{equation}\label{eqn:exp_moments}
\partial_t  \langle f \rangle  +\nabla_x \cdot \langle  \phi v f\rangle = 0\,.
\end{equation}
%Having a $\partial_t\mathcal{M}$ term on the right hand side, the equation~\eq%ref{eqn:exp_reform} does not look easy to compute at the first sight. In fact,% the two equations~\eqref{eqn:exp_reform} and~\eqref{eqn:exp_moments} are coup%led and numerically need to be updated simultaneously. This means at each time% step, one applies the Runge-Kutta method on both, and in each substage, equat%ion~\eqref{eqn:exp_moments} is computed first for the updated moments, and the%n~\eqref{eqn:exp_reform} is computed for $f$. 
To compute $\partial_t\mathcal{M}$, 
note that
\begin{equation}\label{eqn:Maxwellian_time}
\partial_t \mathcal{M} = \partial_\rho \mathcal{M}\partial_t\rho + \nabla_u\mathcal{M}\cdot \partial_tu +\partial_T \mathcal{M}\partial_t T\,,
\end{equation}
where $\partial_\rho \mathcal{M}$, $\nabla_u\mathcal{M}$ can be expressed 
analytically and $\partial_T\mathcal{M}$ are all explicit. The time 
derivatives of the  other three 
macroscopic quantities $\rho, u, T$ can be obtained from (\ref{eqn:exp_moments}). 

With this formulation, one can just use the Runge-Kutta time discretization.

This method preserves positivity, high order accuracy, and strong AP properties.

\subsection{Micro-macro decomposition}

The ``micro-macro" decomposition decomposes the density distribution function into the local
Maxwellian, plus the deviation
\begin{equation}\label{eqn:mm_decomposition}
f =  \mathcal{M} + \varepsilon g\,,\quad\text{with}\quad\int\phi[f- \mathcal{M}]\rd{v} = 0\,.
\end{equation}
One early approach of using such a decomposition to design an AP scheme for the radiative heat transfer equations 
was used by Klar and Schmeiser in \cite{KlarSch}, and it was also used by Liu and Yu
in~\cite{LiuYou} for analyzing the shock propagation of the Euler equations in passing the fluid limit of the Boltzmann equation.
%The idea here is to reformulate the equation such that it is easier to design schemes that ac%hieve AP property. As usual, the Maxwellian function is embedded in the formulation, to enfor%ce the convergence, and the same as in the exponential method, here the Maxwellian function a%lso evolves together with the distribution function.
Its application to the nonlinear Boltzmann equation started with
the work of 
Bennoune-Lemou-Mieussens in~\cite{bennoune2008uniformly}. 
%For application to  
% the linear transport equation see Lemou-Mieussens in~\cite{LM_MM}.

Define the linearized collision operator around $\mathcal{M}$, as:
\begin{equation*}
\mathcal{L}_{\mathcal{M}} g = \mathcal{Q}[ \mathcal{M},g]+\mathcal{Q}[ g,\mathcal{M}]\,.
\end{equation*}
With some calculation, one gets
\begin{equation}\label{eqn:micro_macro}
\begin{cases}
\displaystyle \partial_t g + (\mathrm{I} - \mathrm{\Pi}_ \mathcal{M})(v\cdot\nabla_xg) -  \mathcal{Q}[g,g] = \frac{1}{\varepsilon}\left[\mathcal{L}_ \mathcal{M}g-(\mathrm{I}-\mathrm{\Pi}_ \mathcal{M})(v\cdot\nabla_x \mathcal{M})\right]\,,\\
\displaystyle  \partial_t\int \phi  \mathcal{M}\rd{v} + \int \phi v\cdot\nabla_x  \mathcal{M}\rd{v} + \varepsilon\nabla_x\cdot\langle v\phi g\rangle = 0\,.
\end{cases}
\end{equation}
Here $\mathrm{\Pi}_ \mathcal{M}$ is the  projection operator that maps arbitrary $\mathcal{M}$-weighted $L_2$ function into the null space of $\mathcal{L}_ \mathcal{M}$, namely, for any $\psi\in L_2(\mathcal{M}\rd{v})$:
\begin{equation}
\mathrm{\Pi}_\mathcal{M}(\psi)\in\text{Null}\mathcal{L}_\mathcal{M}=\text{Span}\{ \mathcal{M},v \mathcal{M},|v|^2 \mathcal{M}\}\,.
\end{equation}
For the Boltzmann equation,  the projection operator can be written explicitly as:
\begin{eqnarray}\nonumber
&& \Pi_\mathcal{M}(\psi) = \frac{1}{\rho}\left[\langle\psi\rangle 
 + \frac{(v-u)\cdot\langle(v-u)\psi\rangle}{T}\right.\\ 
&& \qquad \qquad \left.+\left(\frac{|v-u|^2}{2T}-\frac{d}{2}\right)\frac{2}{d}\langle\left(\frac{|v-u|^2}{2T}-\frac{d}{2}\right)\psi\rangle\right]\mathcal{M}\,,
\end{eqnarray}
where $\langle \cdot \rangle$ is the integration over $v$.

In the original Boltzmann equation, the stiff term $\mathcal{Q}[f,f]$ is quadratic in $f$, hence difficult to invert. The two stiff terms here are both linear thus their implicit
discretization can be inverted more easily. 
In \cite{bennoune2008uniformly}, the following discretization is taken:
\begin{equation}\label{eqn:MM_scheme}
\begin{cases}
\displaystyle \frac{g^{n+1} - g^n}{\Delta t}+\left(\mathrm{I} - \mathrm{\Pi}_{ \mathcal{M}^n}\right)(v\cdot\nabla_xg^n) -  \mathcal{Q}[g^n,g^n] = \frac{1}{\varepsilon}\left[\mathcal{L}_{ \mathcal{M}^n}g^{n+1} \right.\\
\qquad \left. - \left(\mathrm{I} - \mathrm{\Pi}_{ \mathcal{M}^n}\right)(v\cdot\nabla_x \mathcal{M}^n)\right],\\
\displaystyle  \int\phi  \mathcal{M}^{n+1}\rd{v} +\Delta t\varepsilon\int \phi v\cdot\nabla_xg^{n+1}\rd{v} = \int\phi  \mathcal{M}^n\rd{v}- \Delta t\int \phi v\cdot\nabla_x \mathcal{M}^{n}\rd{v}
\end{cases}\,.
\end{equation}
The only term that needs to be inverted is $\mathrm{I} - \frac{\Delta t}{\varepsilon}\mathcal{L}_\mathcal{M}$ in the first equation. It is a linear operator, and the negative spectrum of $\mathcal{L}$ guarantees the invertibility.
 The quadratic operator $\mathcal{Q}[f,f]$ is no longer stiff thus is treated explicitly.

The following AP property was proved in~\cite{bennoune2008uniformly}:
\begin{theorem}
The scheme is AP, more specifically:
\begin{itemize}
\item[(i)] The time discretization~\eqref{eqn:MM_scheme} of the Boltzmann equation~\eqref{eqn:micro_macro} gives in the limit $\varepsilon\to 0$ a scheme consistent to the compressible Euler system (\ref{Euler-limit}).
\item[(ii)] For small $\varepsilon$,  scheme~\eqref{eqn:MM_scheme} is asymptotically equivalent, up to $\mathcal{O}(\varepsilon^2)$, to an explicit time discretization of the Navier-Stokes equations (\ref{NS}).
\end{itemize}
\end{theorem}

 In \cite{gamba2019micro} the BGK-penalization method was used in the micro-macro decomposition framework to further avoid the inversion of the linearized collision operator $\mathcal{L}$.

\subsection{Linear transport equations}

\subsubsection{Parity equations-based AP schemes}

We now consider the linear transport equation in diffusive regime \eqref{tone-2}. 
Let
\[
\mathcal{L}(f)=  \int b(v,w) \{ M(v)f(w)-M(w)f(v)\} dw\,.
\]
Split (\ref{tone-2}) into two equations, one for $v$ and one for $-v$:
\begin{eqnarray}\label{thtwo}
\begin{aligned}
\vep\,\partial_t f(v) +v\cdot\nabla_x f(v)
&= \frac{1}{\vep}\LL(f)(v) ,
\\
\vep\,\partial_t f(-v) -v\cdot\nabla_x f(-v) 
&= \frac{1}{\vep}\LL(f)(-v) .
\end{aligned}
\end{eqnarray}
Define the even- and odd-parities  as
\begin{eqnarray}
\begin{aligned}
  r(t,x, v) &= \frac12 \, [ f(t,x, v)+f(t,x,-v)],
\\
  j(t,x, v) &= \frac1{2\vep} \, [ f(t,x, v)-f(t, x,-v)].
  \end{aligned}
\end{eqnarray}
Adding and subtracting the two equations in (\ref{thtwo}) lead to
\begin{eqnarray}
\label{4.7}
\partial_t r+ v\cdot \nabla_x j 
&=&\frac{1}{\vep^2} \LL(r)\,, 
\\
\partial_t j+\frac{1}{\vep^2}  v \cdot \nabla_x r 
&=&-\frac{1}{\vep^2}\,\lambda j \,,
\label{thsix}
\end{eqnarray}
where we used the property that
$$
  \int b(v,w)j(w)\, dw = 0
$$
since $j(w)$ is an odd-function in $w$.

\begin{remark}
 If $b(v, w)= b(|v|, |w|)$, then it is possible
to use the even and odd parities only for the positive components of
$v$ and $w$, hence reduces the computational domain, as is the case for neutron-transport equation \cite{LM84}.
\end{remark}

Since now the convection term is also stiff, the idea of \cite{JPT2} was to rewrite (\ref{4.7}) and (\ref{thsix}) into the following
form
\begin{eqnarray}
\label{r-eq}
&\partial_t r+ v\cdot \nabla_x j 
=\frac{1}{\vep^2} \LL(r) ,
\\
&\partial_t j+ v \cdot \nabla_x r
= -\frac{1}{\vep^2}\,
\left[\lambda j + (1-\vep^2\psi) v \cdot \nabla_x r 
\right]\,,
\label{thseven}
\end{eqnarray}
where $\psi=\psi(\vep)$ is a free parameter satisfying $0 \leq \psi \leq 1/\vep^2$. Hence the characteristic speeds on the right hand side are now independent of $\e$. The simplest
choice of $\psi$ is 
\[
\psi(\vep) = \min\left\{\,1,\,\frac1{\vep^2}\,\right\}.
\label{thphi}
\]
(A related approach in \cite{Klar-AP} moves all the stiff terms in (\ref{thsix}) to the right hand side).

  One can easily derive the diffusion equation from (\ref{r-eq}) and
(\ref{thseven}).
As $\vep\to 0$, they   give 
\begin{eqnarray}
  \LL(r) &=& 0 \,,
\label{Q-eq}
\\
  \lambda j & =& -v \cdot \nabla_x r  \,.
\label{thq}
\end{eqnarray}
Solving  (\ref{Q-eq}) gives
\begin{equation}
  r = \rho(x, t) M(v) \,,
\label{thr}
\end{equation}
where 
\[
  \rho(x, t) =\langle f(x,\cdot,t)\rangle= \langle r(x, \cdot, t) \rangle \,.
\label{thrho}
\]
With \eqref{thr}, equation  (\ref{thq}) gives
\begin{equation}
  j = \frac{M(v)}{\lambda(v)} [ -v \cdot \nabla_x \rho] \,.
\label{thj}
\end{equation}
Applying (\ref{thr}) and (\ref{thj}) in
(\ref{r-eq}), and integrating
over $v$, one gets the diffusion equation (\ref{televen}) with \eqref{diff-coeff}.
Thus (\ref{r-eq}) and  (\ref{thseven}) set the foundation for  AP
schemes.
  One can  split
 the stiff relaxation step
\begin{eqnarray}
\partial_t r
&=&\frac{1}{\vep^2} \LL(r),
\label{4.11}
\\
\partial_t j
&=& \frac{1}{\vep^2}\,
\left[-\lambda j - (1-\vep^2\phi)( v \cdot \nabla_x r )
\right]\,,
\label{theight}
\end{eqnarray}
from the non-stiff  transport step
\begin{eqnarray}\label{4.13}
\begin{aligned}
\partial_t r+ v\cdot \nabla_x j 
&=0,
\\
\partial_t j+ v \cdot \nabla_x r
&=0.
%\label{thnine}
\end{aligned}
\end{eqnarray}
Equations (\ref{4.13})  can be solved using an explicit scheme, whereas for step 
(\ref{4.11})-(\ref{theight}) one uses  an implicit scheme. 

  The key is how to solve the collision step (\ref{4.11}) implicitly in an efficient
way.
In the case of neutron transport, where $L(r)=\rho-r$, the implicit
term can be integrated explicitly \cite{JPT2}. Otherwise, one can use the
penalty method of Filbet-Jin \cite{Filbet-Jin}, see \cite{Deng}.

As far as spatial discretization is concerned, one can use any upwind type
 scheme for convection terms in (\ref{4.13}), while on the right hand side of (\ref{theight}), it was
suggested in \cite{JPT2} to use center difference for the gradient of
$r$. When $\vep \to 0$, these spatial discretizations become consistent
and stable discretization of (\ref{televen}), thus is AP spatially.
However, the limiting discrete diffusion equation is not compact.
In 1d it is a five-point rather than a three-point discretization
of the diffusion equation. This problem can be fixed  by using staggered
grid for $r$ and $j$, as pointed out in \cite{JP-Proc} and then extended to two space dimensions in
 \cite{FrankJin}.

One AP scheme developed in \cite{JiangXu} allows one to get a compact three point scheme in the limit. 

\subsubsection{Micro-macro decomposition based AP schemes}\label{sec:MM}

The micro-macro decomposition approach, 
 proposed by Lemou and Mieussens \cite{lemou2008new}, begins with the decomposition
\begin{equation}
f=\rho M + \vep g\,.
\label{f-decomp}
\end{equation}
 Clearly  $\langle g \rangle=0$.
Applying (\ref{f-decomp}) in (\ref{tone-2}) gives
\begin{equation}
\label{rho-eq}
\vep M \partial_t \rho + \vep^2 \partial_t g + v \cdot M\nabla_x \rho
+\vep v \cdot \nabla_x g = \mathcal{L}g \,.
\end{equation}
Integrating this equation  with respect to $v$ gives the
following continuity equation:
\begin{equation}
 \partial_t \rho +  \nabla_x \cdot \langle vg \rangle\, = 0\,.
\label{mM-1}
\end{equation}
Define operator $\Pi: \Pi(\cdot)(v):=M \langle\cdot\rangle$, and $I$  the identity operator. 
Applying the orthogonal
projection $I-\Pi$ to (\ref{rho-eq})  gives the equation for  $g$: 
\begin{equation}
\vep^2 \partial_t g+ \vep(I-\Pi)(v \cdot \nabla_x g) + v \cdot M\nabla_x
\rho=\LL g \,.
\label{mM-2}
\end{equation}
(\ref{mM-1}) and (\ref{mM-2}) constitute the micro-macro formulation of
(\ref{tone-2}).

We first consider the time discretization.
The following  was used in \cite{lemou2008new}:
\begin{equation}
\frac{g^{n+1}-g^n}{\Delta t}+\frac{1}{\vep}(I-\Pi)(v \cdot \nabla_x g^n)
=\frac{1}{\vep^2} \LL g^{n+1}-\frac{1}{\vep^2} v \cdot M\nabla_x \rho^n\,.
\label{g-dis}
\end{equation}
In the continuity equation (\ref{mM-1}) there is no stiff term, but to
recover the correct diffusion limit, the flux of $g$ is taken at
$t_{n+1}$, which gives
\begin{equation}
\frac{\rho^{n+1}-\rho^n}{\Delta t} + \nabla_x \cdot \langle vg^{n+1} \rangle\,
=0\,.
\label{rho-dis}
\end{equation}
As $\vep \to 0$, (\ref{g-dis}) gives
\[
\LL g^{n+1}=v \cdot M \nabla_x \rho^n\,,
\]
which implies
\[
g^{n+1}(v)=\LL^{-1}(vM) \cdot \nabla_x \rho^n
=\frac{M(v)}{\lambda(v)} \left[ \int b(v,w)g^{n+1}(w) \, dw - v \cdot \nabla_x \rho^n \right]\,.
\]
Applying this to (\ref{rho-dis}), and using the rotational invariance of $\sigma$, yield
 the following time explicit discretization
of the diffusion equation \eqref{televen}
\[
\frac{\rho^{n+1}-\rho^n}{\Delta t} + \nabla_x \cdot \langle D \nabla \rho^{n}\rangle\,
=0\,.
\label{discrete-rho1}
\]
Thus this time discretization is  AP.

Now consider the case of one space dimension.
A staggered grid can be used.  Define $x_{i+1/2}=(i+1/2)\Delta x$.
Now the macroscopic density $\rho$ will be defined at grid point $x_i$,
while $g$ is defined at $x_{i+1/2}$. Using upwind discretization for
the space derivative, one arrives at
\begin{eqnarray}
\label{LM-1}
&&
  \frac{\rho_i^{n+1}-\rho_i^n}{\Delta t} + \left< v \frac{g^{n+1}_{i+1/2}
-g^{n+1}_{i-1/2}}{\Delta x}\right>=0\,,
\\
&&
\nonumber
   \frac{g^{n+1}_{i+1/2}
-g^{n}_{i+1/2}}{\Delta x} +\frac{1}{\vep \Delta x}(I-\Pi)
\left(  v^+(g^n_{i+1/2}-g^n_{i-1/2} ) \right. \\
\nonumber
&& \qquad  \left. +
 v^-(g^n_{i+3/2}-g^n_{i+1/2} ) \right) 
\\
\label{LM-2}
&&=\frac{1}{\vep^2} \LL g^{n+1}_{i+1/2} 
 -\frac{1}{\vep^2}v M \frac{\rho^n_{i+1}-\rho_{i}^n}{\Delta x} \,,
\end{eqnarray}
where $v^\pm = (v\pm|v|)/2$.

As $\vep \to 0$, (\ref{LM-2}) gives
\[
g^{n+1}_{i+1/2} = \LL^{-1} (vM) \frac{\rho^n_{i+1}-\rho^n_{i}}{\Delta x}
\]
which when applied to (\ref{LM-1}) gives the following scheme
\[
  \frac{\rho^{n+1}-\rho^n}{\Delta t} + D 
\frac{\rho^{n}_{i+1}-2\rho^{n}_{i} +\rho^{n}_{i-1}}
{(\Delta x)^2} = 0\,.
\]
This is the classical three point explicit discretization of the diffusion
equation (\ref{televen}) and \eqref{diff-coeff}.

The uniform stability condition ($\Delta t \le C (\Delta x)^2$, uniformly in $\vep$)
of this method was proved in \cite{LM10}.

Among all the above approaches, in the limit $\vep \to 0$, the 
discrete time discretization  is explicit for the limiting diffusion equation. This imposes the numerical stability condition like $\Delta t=O((\Delta x)^2)$. Consider the case of $\mathcal{L}=r-\rho$, hence $\lambda=1$, in the parity formulation \eqref{4.7} and \eqref{thsix}. In \cite{Russo-CFL}, the authors proposed to reformulate the system into
\begin{eqnarray}
 &&   \partial_t r = \underbrace{-v \partial_x \left( j+ \frac{\nu(\varepsilon) v \partial_x r}{\sigma} \right)}_{\text{explicit}} - \underbrace{\frac{\sigma}{\varepsilon^2}(r-\rho)  + \nu(\varepsilon)v^2 \frac{\partial_{xx} r}{\sigma}}_{\text{implicit}},\\
 && \partial_t j = -\underbrace{\frac{1}{\varepsilon^2} \left(j+ \frac{v\partial_x r}{\sigma}\right)}_{\text{implicit}},
\end{eqnarray}
where $\mu(\e)\in [0, 1]$  is a free parameter such that $\mu(0)=1$.
$\mu=1$ guarantees the largest stability region.
When $\vep \to 0$, one gets an {\it implicit} discretization of the diffusion equation, enabling a stability condition like
$\Delta t=O(\Delta x)$.

\subsection{Stochastic AP schemes for linear transport equation with uncertainties}

Kinetic models usually have  {\it uncertainties} that  can arise in  collision kernels, scattering coefficients, initial
or boundary data, geometry, source or forcing terms \cite{Bird}\cite{BHW86}\cite{KM91}.
  Understanding the impact of, quantify and even control  these uncertainties, in the sense of uncertainty quantification (UQ),  is crucial to the simulations of the complex kinetic systems in order to verify,  validate and improve these models, and to conduct risk management.  

The uncertainty is usually modelled by  a random vector $z\in \mathbb{R}^n$ in a properly defined probability space $(\Sigma, \mathcal{A}, \mathbb{P})$, whose event space is $\Sigma$ and  equipped with $\sigma$-algebra $\mathcal{A}$ and probability measure $\mathbb{P}$. We also assume the components of $z$ are mutually independent random variables with known probability $\omega(z): I_{z} \longrightarrow \mathbb{R}^+$, obtained already through some dimension reduction technique, e.g., Karhunen-Lo\`{e}ve (KL) expansion \cite{Loeve}.

\subsubsection{The linear transport equation with isotropic scattering}

Consider the linear transport equation in one dimensional slab geometry with random input:
\begin{align}\label{eq-onegroup}
  & \eps \delta_t f + v \delta_x  f = \frac{\sigma}{\eps}{\cal L} f,  \quad  t>0, \  x\in [0,1], \  v\in[-1,1], \ z\in I_{z},
  \\
  \label{L}
  & {\cal L} f(t,x,v,z) = \frac{1}{2} \int_{-1}^1 f(t, x,v',z) d v' - f(t,x,v,z)\,,
\end{align}
with the  initial condition
\begin{equation}\label{IC}
  f(0,x,v,z) = f^0(x,v,z).
\end{equation}
This equation arises in neutron transport, radiative transfer, etc. and describes particles (for example neutrons) transport in a background media (for example
nuclei). $v=\Omega\cdot e_x = \cos \theta$ where
$\theta$ is the angle between the moving direction and $x$-axis.
%$$
%\eps \dt f + \Omega \cdot \nabla f = \frac{\sigma}{\eps}{\cal L} f - \eps \sig%ma_a f + \eps S, \quad
%\sigma \ge \sigma_{min} > 0
%$$
Assume
\begin{equation}
\sigma(x, z) \ge \sigma_{\text{min}} > 0.
\end{equation}
 
%The Dirichlet boundary conditions are given  in the incoming direction by
%\begin{equation}\label{BC}
%  \begin{split}
%  & f(t,0,v,z) = f_{\mathrm{L}}(t,v,z),  \qquad \mbox{for }v\ge 0\,,
%  \\
%  & f(t,1,v,z) = f_{\mathrm{R}}(t,v,z),  \qquad\mbox{for } v\le 0\,,
%  \end{split}
%\end{equation}

Denote
\begin{equation}\label{ave}
  \langle \phi \rangle =\half\int_{-1}^1\phi(v) dv
\end{equation}
as the average of a velocity dependent function $\phi$.

Let $\rho = \langle f \rangle$. For each  fixed $z$,  as $\varepsilon \to 0$, $\rho$ solves the following diffusion equation:
\begin{equation} \label{eq-diff1}
  \partial_t \rho =  \partial_{x} \left( \frac{1}{3}\sigma(x,z)^{-1} \partial_{x}\rho\right).
\end{equation}

In order to understand the property of numerical methods for uncertain kinetic equations, it is importnt to study the regularity and long-time behavior in the random space of the linear transport equation (\ref{eq-onegroup})-(\ref {IC}).
Consider the Hilbert space of the random variable
\begin{equation}
  H(I_z;\;\omega dz) = \Big\{\,f\mid I_z \rightarrow \mathbb{R}^+, \;\int_{I_z} f^2(z)\omega(z) dz < +\infty\,\Big\},
\end{equation}
equipped with the inner product and norm defined as
\begin{equation}\label{inner-norm}
  \langle f, g \rangle_\omega= \int_{I_z} fg\,\omega(z) dz,\quad \| f \|_{\omega}^2 = \langle f, f \rangle_{\omega}\,.
\end{equation}
 Define the $k$th order differential operator with respect to $z$ as
\begin{equation}
  D^k f(t,x,v,z) := \partial^k_z f(t,x,v,z),
\end{equation}
and the Sobolev norm in $z$ as
\begin{equation}
  \|f(t,x,v,\cdot)\|_{H^k}^2 := \sum_{\alpha\le k} \|D^\alpha f(t,x,v,\cdot)\|_\omega^2.
\end{equation}
Finally, introduce norms in space and velocity as follows,
\begin{align}
  & \|f(t,\cdot,\cdot,\cdot)\|_{\Gamma}^2:=\int_Q \|f(t,x,v,\cdot)\|_{\omega}^2\, dx\,dv,\qquad t\geq 0,
  \\
% & \|f(t,\cdot,\cdot,\cdot)\|_{\Gamma^k}^2:=\int_Q %\|f(t,x,v,\cdot)\|_{H^k}^2\,\diff x\diff v,\qquad t\geq 0,
\end{align}
where $Q=[0,1]\times [-1,1]$ denotes the domain in the phase space. The following results
 were established in \cite{JLM}.

\begin{theorem}[{\bf Uniform regularity}]\label{rgl1}

  If for some integer $m\ge 0$,
  \begin{equation}
    \|D^k\sigma(z)\|_{L^\infty}\le C_\sigma, \qquad \|D^k f_{0}\|_{\Gamma}\le C_0,\qquad k=0,\dots,m,
  \end{equation}
  then the solution $f$ to the linear transport equation (\ref{eq-onegroup})--(\ref{IC}), with  periodic boundary condition in $x$, satisfies,
  \begin{equation}
    \|D^k f(t, \cdot, \cdot, \cdot)\|_{\Gamma} \le C, \qquad k=0,\cdots,m, \qquad \forall t>0,
  \end{equation}
  where $C_\sigma$, $C_0$ and $C$ are constants independent of $\eps$.
\end{theorem}

The above theorem shows that, under some smoothness assumption on $\sigma$, the regularity of the initial data is preserved  in time and the Sobolev norm
of the solution  is
bounded uniformly in $\eps$.

\subsubsection{Stochastic Galerkin approximation}

An interesting and important scenario  is when the uncertainty  and small scaling are both present in the equation. Among various UQ methods \cite{Xiu}\cite{GWZ-Acta}, we consider the stochastic Galerkin (SG) method, which is suitable for our AP analysis thanks to its Galerkin formulation.

Assume the  complete orthogonal polynomial basis in the
Hilbert space $H(I_z; \omega(z) dz)$ corresponding to
the weight $\omega(z)$ is $\{\phi_i(z), i=0,1, \cdots,\}$, where
$\phi_i(z)$ is a polynomial of degree $i$ and satisfies the orthonormal condition:
$$
   \langle \phi_i, \phi_j \rangle_\omega=\int \phi_i(z)\phi_j(z)\omega(z) dz=\delta_{ij}.
$$
Here $\phi_0(z)=1$, and $\delta_{i j}$ is the Kronecker delta function. Since the solution $f(t,\cdot,\cdot,\cdot)$ is defined in
 $L^2\big( [0,1]\times [-1,1]\times \mathbb I_z; d\mu)$, one has the generalized polynomial chaos expansion \cite{xiu2002wiener}
$$
   f(t,x,v,z) = \sum_{i=0}^{\infty} f_i(t,x,v) \, \phi_i(z),  \quad
   \hat f = \big(\, f_i\, \big)_{i=0}^\infty:=\big(\bar f, \hat f_1\big).
$$
The mean and variance of $f$ can be obtained from the expansion coefficients
as
$$
  \bar f= E (f) = \int_{I_z}  f   \omega(z)\, dz = f_0, \quad \mbox{ var }(f) = |\hat f_1|^2\,.
$$

Denote the SG solution by
\begin{equation}
\label{gPC}
   f^K  = \sum_{i=0}^{K} f_i \, \phi_i,  \quad  \hat f^K = \big(\, f_i \, \big)_{i=0}^K
   := \big(\bar f, \hat f_1^K\big),
\end{equation}
from which one can extract the  mean and  variance of $f^K $ from
the expansion coefficients as
$$
   E (f^K ) = \bar f, \quad \mbox{ var }(f^K ) = |\hat f_1^K|^2\,.
$$
Furthermore,  define
\begin{eqnarray}
\nonumber
&& \sigma_{ij}= \big\langle \phi_i,\, \sigma \phi_j \big\rangle_\omega, \quad
\Sigma = \big( \,\sigma_{ij} \, \big)_{M+1, M+1}; \\
&& \sigma^a_{ij}= \big\langle \phi_i, \, \sigma^a \phi_j \big\rangle_\omega, \quad
\Sigma^a = \big(\, \sigma^a_{ij} \,\big)_{M+1, M+1},
\end{eqnarray}
for $0\le i,j \le M$. Let $\mbox{ I }$ be the $(K+1)\times (K+1)$ identity
matrix. $\Sigma, \Sigma^a$ are symmetric positive-definite matrices satisfying
$$
\Sigma \ge \sigma_{\text{min}} \mbox{ I }.
$$

If one applies the polynomial chaos ansatz~(\ref{gPC}) into the transport equation~(\ref{eq-onegroup}), and conduct the Galerkin projection, one obtains
%linear transport in PC form
%\begin{equation}  \label{eq-onegroup}
%\eps \dt \hat f + v \dx \hat f = - \frac{1}{\eps}(I-[ \cdot ]) \Sigma  \hat f
%\end{equation}

%linear transport in $m$ term truncated PC (spectral Galerkin approximation)
\begin{equation}  \label{eq-gPC-1}
\eps \partial_t \hat f + v \partial_x \hat f = - \frac{1}{\eps}(I-[ \cdot ]) \Sigma  \hat f .
\end{equation}

%Denote $A=(I-[ \cdot ]) \Sigma_m$, then the solution to
%(\ref{eq-onegroup}) can be expressed
%as
%$$
%   \hat f(x,v,t) = e^{-A t/ \eps^2} \hat f_{in}(x-\frac{vt}{\eps}, v)
%$$
%As a direct consequence, if each component of $\hat f_{in}(x-vt, v)$ is nonnega%tive,
%then each component of $\hat f(x,v,t)$ is also nonnegative. ({\bf {Question:
%%o we need this? the coefficients are not positive anyway}})

Note the SG method
makes the random transport equations into  deterministic systems \eqref{eq-gPC-1} which are
vector analogue of the original scalar deterministic  transport equations. Therefore
one can naturally utilize the deterministic AP machinery to solve the SG system to achieve the desired AP goals, hence   minimize ``intrusion'' to the legacy deterministic codes. To this aim, the notion of {\it stochastic asymptotic preserving (sAP)} was introduced in \cite{JXZ15}. A scheme is sAP if an SG method for the random kinetic equation becomes an SG approximation for the limiting macroscopic, random (hydrodynamic or diffusion) equation  as $\varepsilon \rightarrow 0$, with $K$, mesh size and time step all held fixed. Such schemes guarantee that  for $\varepsilon \rightarrow 0$, {\it all}
numerical parameters, including $K$, can be chosen only for accuracy requirement, but {\it independent} of $\varepsilon$.

We now use the  micro-macro decomposition:
\begin{equation}
\label{micro-macro}
   \hat f(t,x,v, z) = \hat{\rho}(t,x, z) + \eps \hat{g}(t,x,v, z),
\end{equation}
where $\hat \rho=[\hat f]$ and $[\hat{g}]=0$, in~(\ref{eq-gPC-1}) to get
\begin{subequations} \label{eq-mM-gPC}
\begin{align}
& \partial_t \hat \rho + \partial_x \langle v \hat g \rangle = -\Sigma^a \hat \rho + \hat S ,\label{eq-mMrho-gPC}  \\
& \partial_t \hat g + \frac{1}{\eps} (I-\langle . \rangle ) (v\partial_x \hat g) =
-\frac{1}{\eps^2} \Sigma \hat g - \Sigma^a \hat g
- \frac{1}{\eps^2}v \partial_x \hat \rho,\label{eq-mMg-gPC}
\end{align}
\end{subequations}
with initial data
$$
  \hat \rho(0,x,z) =  \hat \rho_0(x,z), \quad
  \hat g(0,x,v,z) =  \hat g_0(x,v,z)\,.
$$

As $\e \to 0$,  system~(\ref{eq-mM-gPC}) formally approaches  the  diffusion limit
\begin{equation}  \label{eq-diff-gPC}
\partial_t  \hat \rho =  \partial_{x} \left( \frac{1}{3} \Sigma^{-1} \partial_{x}\hat \rho \right) \,.
\end{equation}
This is the SG approximation to the random diffusion equation
(\ref{eq-diff1}).
 Thus the SG approximation is sAP in the sense of~\cite{JXZ15}.

%Homogenization QUESTION.
%Derive a close equation for $\bar \rho$ and prove the convergence.

The following result was proved in \cite{JLM}.

\begin{theorem}\label{uni_conv}
  If for some integer $m\ge 0$,
  \begin{equation}
    \|\sigma(z)\|_{H^k}\le C_\sigma, \quad \|D^k f_{0}\|_{\Gamma}\le C_0, \quad \|D^k(\partial_x f_0)\|_{\Gamma} \le C_x, \quad k=0,\dots,m,
  \end{equation}
  then for $t\le T$, the error of the sG method is
  \begin{equation}\label{err}
  \|f - f^K \|_{\Gamma}\le \frac{C(T)}{K^k},
  \end{equation}
  where $C(T)$ is a constant independent of $\eps$.
\end{theorem}

Theorem~\ref{uni_conv} gives a uniform in $\eps$ spectral convergence rate, thus one can choose $K$ independent of $\eps$, a very strong sAP property.
Such a result is also obtained with the anisotropic scattering case, for the
linear semiconductor Boltzmann equation \cite{JinLiuL17}.

\subsubsection{A full discretization}
\label{sec:full}

 Here, we adopt the micro-macro decomposition based fully discrete  scheme  for the SG system~(\ref{eq-mM-gPC}).

%The fully discrete scheme for system (\ref{eq-mM}) for $i = 0, 1, \cdots N-1$
 %is \cite{LM-08}:
%\begin{subequations} \label{eq-scheme}
%\begin{align}
%& \frac{\rhonpi-\rhoni}{\Dt} +
%\av{v\frac{\gnpi-\gnpim}{\Dx}} = -\sigma^a_i \rho_i^{n+1} + S_i,
% \label{eq-schrho}  \\
%& \frac{\gnpi-\gni}{\Dt} + \frac{1}{\eps\Dx} (I-\av{.})
%\left(v^+(\gni-\gnim)+v^-(\gnip-\gni)\right)
%\label{eq-schg}  \\
%& \qquad\qquad =
%- \frac{\sigma(z)}{\eps^2} \gnpi -
%\frac{1}{\eps^2}v\frac{\rhonip-\rhoni}{\Dx},
%\nonumber
%\end{align}
%\end{subequations}
%where $v^\pm=\frac{v\pm|v|}{2}$.

%For simplicity we use periodic BCs: $ \rho^n_N=\rho^n_0$, $g^{n}_{-1/2}=g^{n}_%{N-1/2}$.  ({\bf {is this necessary?}}, yes)

%For each $z$, as $\eps \to 0$, this scheme formally has the following random diffusion limit
%\cite{LM-08}:
%\begin{equation}  \label{eq-sdiff}
%\frac{\rhonpi-\rhoni}{\Dt} +
%\kappa(z) \, \frac{\rhonip-2\rhoni+\rhonim}{\Dx^2}= -\sigma^a_i \rho_i + S_i.
%\end{equation}
%with $\kappa(z) =\frac{1}{3} \sigma(z)^{-1}$.
Corresponding to \eqref{LM-1} and \eqref{LM-2}, one has 
\begin{subequations} \label{eq-scheme}
\begin{align}
& \frac{\hat{\rho}^{n+1}_i-\hat{\rho}^{n}_i}{\Delta t} +
\langle {v\frac{\hat{g}^{n+1}_{i+1/2}-\hat{g}^{n+1}_{i-1/2}}{\Delta x}} \rangle  = 0 ,
 \label{eq-schrho}  \\
& \frac{\hat{g}^{n+1}_{i+1/2}-\hat{g}^n_{i+1/2}}{\Delta t} + \frac{1}{\eps\Delta x} (I-\langle \Pi \rangle)
\left(v^+(\hat{g}^n_{i+1/2}-\hat{g}^n_{i-1/2})+v^-(\hat{g}^n_{i+3/2}-\hat{g}^n_{i+1/2})\right)
\label{eq-schg}  \\
& \qquad\qquad =
- \frac{1}{\eps^2}\Sigma_i \hat{g}^{n+1}_{i+1/2}   -
\frac{1}{\eps^2}v\frac{\hat{\rho}^n_{i+1}-\hat{\rho}^n_{i}}{\Delta x}. \nonumber
\end{align}
\end{subequations}
Its formal  limit, when $\eps \rightarrow 0$, is given by
\begin{equation}  \label{eq-sdiff}
\frac{\hat{\rho}^{n+1}_{i}-\hat{\rho}^n_{i}}{\Dt} -
 \frac{1}{3}  \Sigma^{-1}\, \frac{\hat{\rho}^n_{i+1}-2\hat{\rho}^n_{i}+\hat{\rho}^n_{i-1}}{\Delta x^2}= 0 .
\end{equation}
 This is the fully discrete sG scheme for~(\ref{eq-diff-gPC}). Thus the fully discrete  scheme is  sAP.

One important property for an AP scheme is to have a stability  condition
independent of
$\varepsilon$, so one can take $\Delta t \gg O(\varepsilon)$.
The next theorem from \cite{JLM} confirms this.

\begin{theorem} \label{theo:stab} 
If $\Delta t$ satisfies the following CFL condition
\begin{equation}  \label{eq-CFL}
\Delta t\leq  \frac{\sigma_{\mathrm{min}}}{3}(\Delta x)^2 + \frac{2 \eps}{3}\Delta x ,
\end{equation}
then the solution obtained  by scheme~(\ref{eq-scheme}) satisfies the energy estimate
\begin{equation*}
%  \rnorm{\hat{\rho}^n_i}^2+ \eps^2\gnorm{\hat{g}^n_{i+1/2}}^2 \leq \rnorm{\rho^0}^2+ \eps^2\gnorm{g^0}^2
 \sum_{i=0}^{N-1} \left(  \left( \hat{\rho}^{n+}\right)^2 + \frac{\eps^2}2 \int_{-1}^{1} \left(\hat{g}^n_{i+1/2}\right)^2\, dv \right)
\leq
\sum_{i=0}^{N-1}
\left(  \left(\hat \rho^0_i\right)^2 + \frac{\eps^2}2 \int_{-1}^{1}\left(\hat g^0_{i+\half}\right)^2 dv \right)
\end{equation*}
for every $n$, and hence the scheme~(\ref{eq-scheme}) is stable.
\end{theorem}

Since the right hand side of~(\ref{eq-CFL}) has a lower bound,
which is essentially  a stability condition of the
discrete diffusion equation~(\ref{eq-sdiff})), when $\eps \to
0$, the scheme is asymptotically
stable and  $\Delta t$ remains finite even if $\eps \to 0$.

Next we consider a numerical example from \cite{JLM}.  Consider a random coefficient with one dimensional random parameter:
$$
\sigma(z) = 2 + z, \quad z \mbox{ is uniformly distributed in } (-1, 1).
$$
The limiting random diffusion equation  is
\begin{equation}\label{eq-limiting} 
   \partial_t\rho = \frac1{3\sigma(z)} \partial_{xx}\rho\,,
\end{equation}
with initial condition and boundary conditions:
$$
  \rho(t,0,z) =1, \quad \rho(t,1,z) =0, \quad  \rho(0,x,z) =0.
$$
The analytical solution for~(\ref{eq-limiting}) with the given initial and boundary conditions  is
\begin{equation}\label{limiting solution}
   \rho(t,x,z) = 1 - \mbox{ erf } \left( \frac{x}{\sqrt{\dfrac{4}{3\sigma(z)} t}} \right).
\end{equation}
When $\eps$ is small,  this can be used as the reference solution.  For large $\eps$ or in the case one can not get an analytic solution, we will use the collocation method (see~\cite{GWZ-Acta}) with the same time and spatial discretization to the micro-macro system~(\ref{eq-scheme}) as a comparison in the following examples. In addition, the standard 30-points Gauss-Legendre quadrature set is used for the velocity space to compute $\rho$.

To examine the accuracy,  two error norms are used: the differences in the mean solutions and in the corresponding standard deviation, with $\ell^2$ norm in $x$:
\begin{equation*}
  e_{mean}(t) = \big\|\mathbb{E}[u^h]-\mathbb{E}[u]\big\|_{\ell^2}\,,
\end{equation*}
\begin{equation*}
  e_{std}(t) = \big\|\sigma[u^h]-\sigma[u]\big\|_{\ell^2}\,,
\end{equation*}
where $u^h,u$ are the numerical solutions and the reference solutions, respectively.
\begin{figure}[htbp]  
  \centering
  \includegraphics[width=0.8\textwidth]{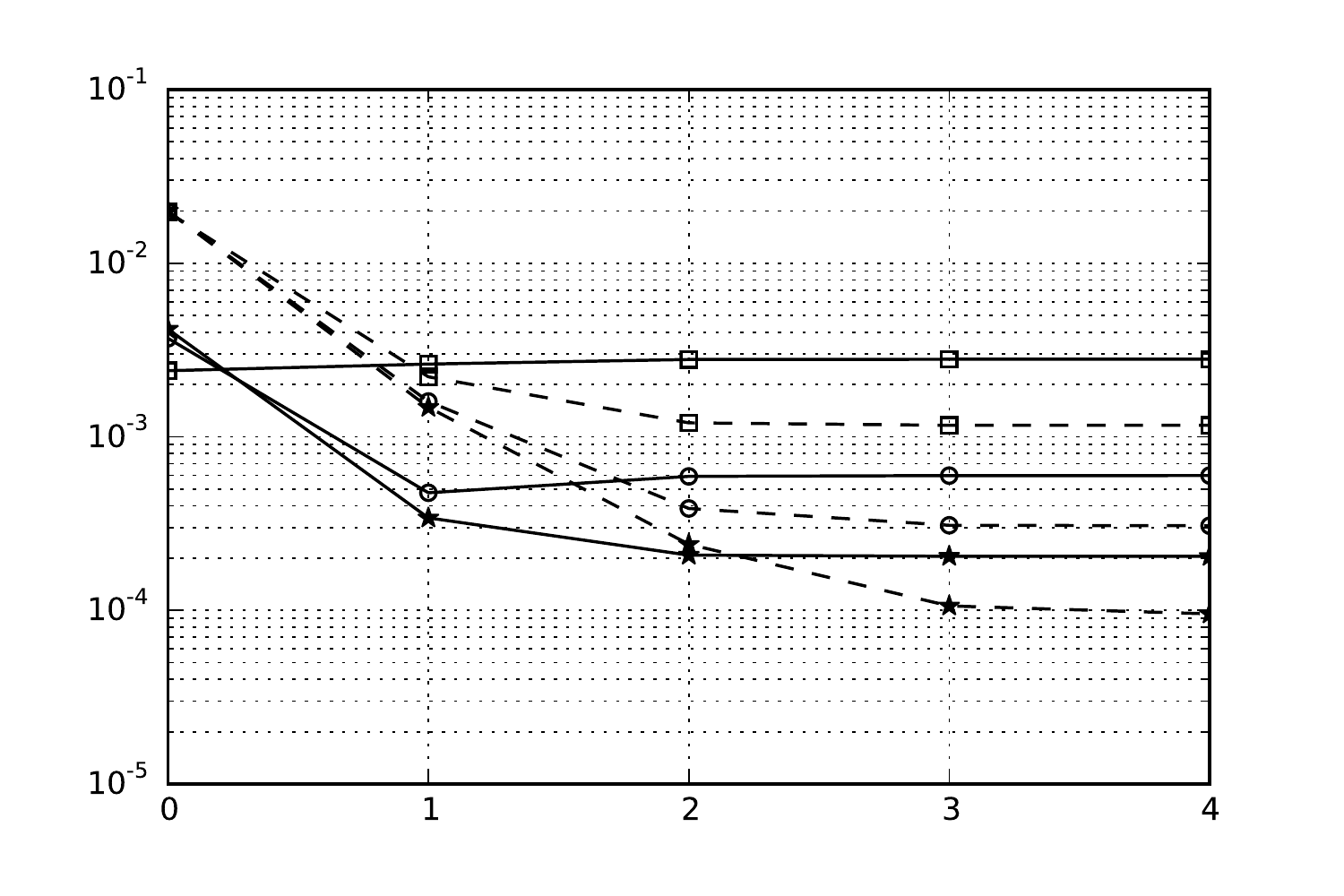}
  \caption{Example 1. Errors of the mean (solid line) and standard deviation (dash line) of $\rho$ with respect to the polynomial chaos order $K$ at $\e = 10^{-8}$: $\Delta x = 0.04$ (squares), $\Delta x = 0.02$ (circles), $\Delta x = 0.01$ (stars).}
  \label{Fig:1}
\end{figure}

In Figure~\ref{Fig:1},  the errors in mean and standard deviation of the SG  solutions at $t = 0.01$ with different $K$ are plotted.  Three sets of results are included: solutions with $\Delta x = 0.04$ (squares), $\Delta x = 0.02$ (circles), $\Delta x = 0.01$ (stars), with $\Delta t = 0.0002/3$ always used. One observes that the errors become smaller with finer mesh, and  the solutions decay rapidly in $K$ and then saturate where spatial discretization error dominates.

\begin{figure}[htbp]  
  \centering
  \includegraphics[width=0.8\textwidth]{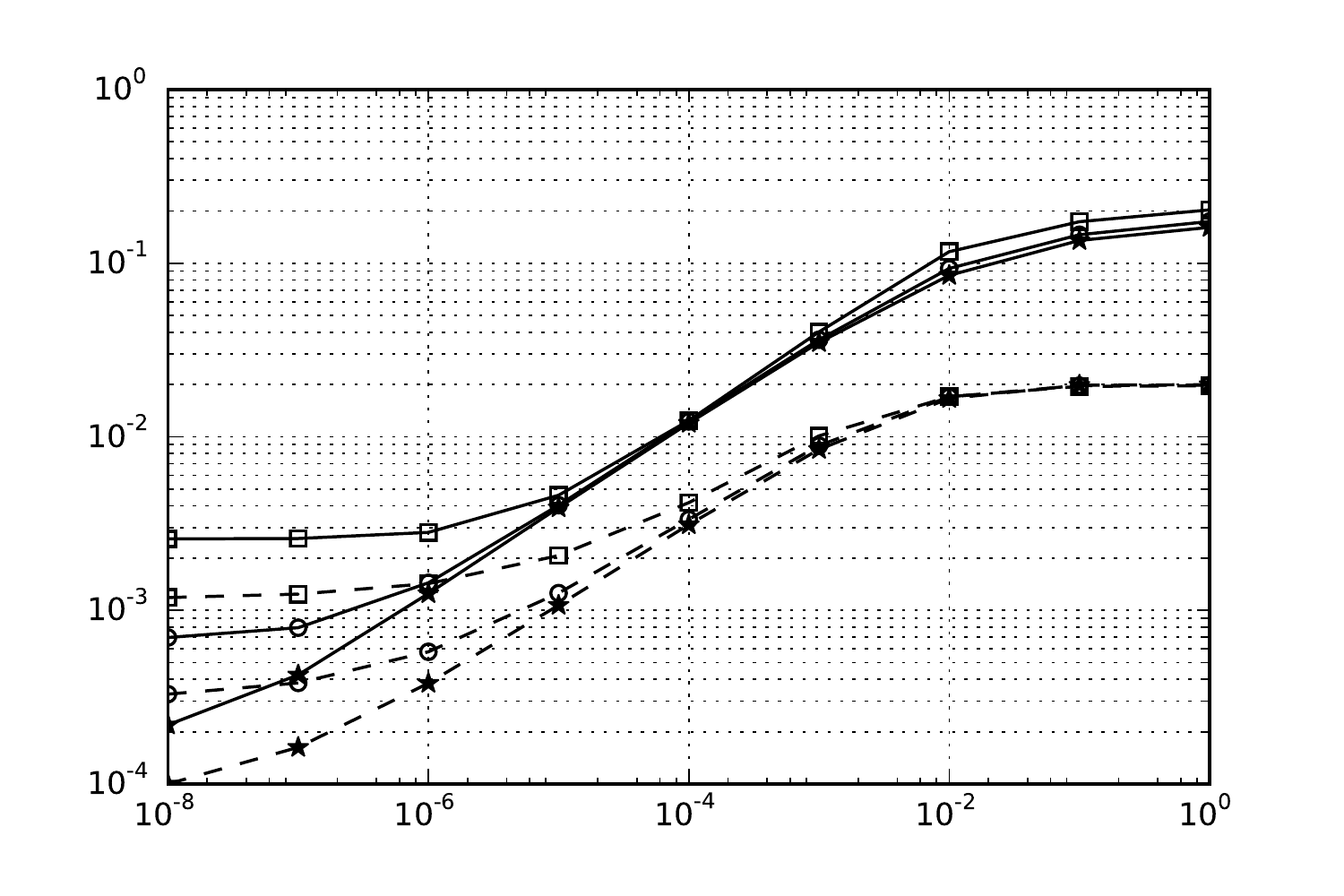}
  \caption{Example 1. Differences in the mean (solid line) and standard deviation (dash line) of $\rho$ with respect to $\eps^2$, between the limiting analytical solution~(\ref{limiting solution}) and the SG solution with $K=4$, $\Delta x = 0.04$ (squares), $\Delta x = 0.02$ (circles) and $\Delta x = 0.01$ (stars).}
  \label{Fig:3}
\end{figure}

In Figure~\ref{Fig:3}, we examine the difference between the solution $t = 0.01$ obtained by SG  with $K=4$, $\Delta x = 0.01$, $\Delta t = \Delta x^2/12$ and the limiting analytical solution~(\ref{limiting solution}). One can  observe the differences become smaller as $\eps$ is smaller in a quadratic fashion, before the numerical errors become dominant. Therefore the method works for all range of $\e$.

A discontinuous Galerkin method based sAP scheme for the same problem
 was developed in
\cite{ChenLiuMu}, where uniform stability and  rigorous sAP property were
also proven.

\subsection{Stochastic Galerkin  methods for general nonlinear  kinetic equations with uncertainties}\label{sec:3}

Consider a general nonlinear kinetic equation with multi-dimensional uncertainties:

\begin{align} \label{gPCsys-0}
\left\{
\begin{array}{l}
\e^\alpha \displaystyle \partial_t f+v\cdot \nabla_{x} f-\nabla_{\bx}\phi\cdot \nabla_{v}f=\frac{1}{\e}\mathcal{Q}(f),\quad t> 0, \ x\in  \Omega,\  v\in \mathbb{R}^d,  \ z\in \mathbb{R}^n,\\[6pt]
\displaystyle  f(0,x,v)=f^0(x,v).
\end{array}\right.
\end{align}
Here $\alpha=0$ or $1$ corresponds to the Euler (acoustic) or incompressible Navier-Stokes scaling respectively \cite{BGL91}. 

We again use the generalized polynomial chaos approximation
\begin{align} \label{gPCexp}
&f(t,x,v,z)\approx \sum_{|k|=0}^{K}f_{k}(t,x,v)\Phi_{k}(z):=f^K (t,x,v,z),
\end{align}
where $k=(k_1,\dots,k_n)$ is a multi-index with $|k|=k_1+\dots+k_n$. $\{\Phi_{k}(z)\}$ are orthornomal polynomials from $\mathbb{P}_K^n$, the set of all $n$-variate polynomials of degree up to $M$ and satisfy
\begin{equation*}
\langle \Phi_{k}, \Phi_{j} \rangle_\omega= \int_{I_{z}}\Phi_{k}(z)\Phi_{j}(z)\omega(z)\,\rd{z}=\delta_{kj}, \quad 0\leq|k|,|j|\leq K.
\end{equation*}
Here  $\delta_{kj}$ is the Kronecker delta function. 

Now inserting (\ref{gPCexp}) into \eqref{gPCsys-0}. Upon a standard Galerkin projection, one obtains, for each $0\leq k\leq M$,
\begin{align} \label{gPCsys}
\left\{
\begin{array}{l}
\e^\alpha \displaystyle \partial_t f_{k}+v\cdot \nabla_{x} f_{k}-\sum_{|j|=0}^K\nabla_{x}\phi_{kj}\cdot \nabla_{v}f_{j}=\frac{1}{\e}\mathcal{Q}_{k}(f^K ), \\[6pt]
  \displaystyle  f_{k}(0,x,v)=f^0_{k}(x,v), 
\end{array}\right.
\end{align}
with
\begin{align}
&\mathcal{Q}_{k}(f^K ):=\int_{I_{z}}\mathcal{Q}(f^K )(t,x,v,z)\Phi_{k}(z)\omega(z)\,dz, \nonumber\\ &\phi_{kj}:=\int_{I_{z}}\phi(t,x,z)\Phi_{k}(z)\Phi_{j}(z)\omega(z)\,dz, \nonumber \\
&f_{k}^0:=\int_{I_{z}}f^0(x,v,z)\Phi_{k}(z)\omega(z)\,dz.
   \nonumber
\end{align}
 We also assume that the potential $\phi(t,x,z)$ is given a priori for simplicity (the case that it is coupled to a Poisson equation can be treated similarly \cite{ZhuJin17}).

%------------------------------------------------------------------
\subsubsection{Hypocoercivity estimate of the SG system}
\label{gPC1}

The hypocoercivity theory \cite{villani-Hypo} can be used to
study the properties of the SG methods.  For general linear transport with uncertainty see
 \cite{LiWang-hypo}. For nonlinear problems one needs to consider the perturbative form \cite{ZhuJin18}\cite{LiuJin18}
\begin{equation}
  \label{ans}
f_{k}=\mathcal M +\varepsilon\, M h_{k}\,,
\end{equation}
  where $h_{k}$ is the coefficient of the following gPC expansion
  \begin{equation}
  \label{h_k}
    h(t,x,v,z)\approx \sum_{|k|=0}^{M}h_{k}(t,x,v)\Phi_{k}(z):=h^K (t,x,v,z)\,.
\end{equation}   
 Inserting ansatz (\ref{ans}) and (\ref{h_k}) into (\ref{gPCsys})  and conducting a standard Galerkin projection, one obtains the SG
system for $h_{k}$ (consider the case of $\phi=0$) \cite{HuJin-UQ}:
\begin{align}
\label{h_gPC1}
\left\{
\begin{array}{l}
\displaystyle \partial_t h_{k} + \frac{1}{\varepsilon^\alpha}v\cdot\nabla_x h_{k} =\frac{1}{\varepsilon^{1+\alpha}}\mathcal L_{k}(h^K) +
\frac{1}{\varepsilon^\alpha}\mathcal F_k(h^K, h^K),
\\[2pt]
\displaystyle  h_{k}(0,x,v)=h_{k}^{0}(x,v), \qquad x\in\Omega\subset\mathbb T^d,  \, v\in \mathbb R^d,
\end{array}\right.
\end{align}
for each $1\leq |k|\leq K$, with initial data given by
$$h_{k}^{0} := \int_{I_z}\, h^{0}(x,v,z)\psi_{k}(z)\, \pi(z)dz. $$
For the Boltzmann equation, the collision parts are given by
\begin{eqnarray*}
&& \mathcal L_{k}(h^K)=
 \mathcal L^{+}_{k}(h^K)=\sum_{|i|=1}^{K}\,\int_{\mathbb R^d\times\mathbb S^{d-1}}\,\widetilde S_{ki}\,
                          \phi(|v-v_{\ast}|)\, (h_{i}(v^{\prime}) M(v_{\ast}^{\prime}) \\
                          && \qquad\qquad\qquad\qquad\qquad + h_{i}(v_{\ast}^{\prime})M(v^{\prime}))\, M(v_{\ast})\, dv_{\ast}d\sigma \\
  && \quad
- M(v)\, \sum_{|i|=1}^{K}\,\int_{\mathbb R^d\times\mathbb S^{d-1}}\, \widetilde S_{ki}\,
 \phi(|v-v_{\ast}|)\, h_{i}(v_{\ast}) M(v_{\ast})\, dv_{\ast}d\sigma
     - \sum_{|i|=1}^K \nu_{ki} h_{i} \,,
  \\
&& \mathcal F_{k}(h^K, h^K) (t,x,v)=\sum_{|i|, |j|=1}^{K}\,\int_{\mathbb R^d\times\mathbb S^{d-1}}\, S_{kij}\,
\phi(|v-v_{\ast}|)\, M(v_{\ast})\\
&& \qquad\qquad\qquad\qquad\qquad \cdot (h_{i}(v^{\prime})h_{j}(v_{\ast}^{\prime}) - h_{i}(v)h_{j}(v_{\ast}))\, dv_{\ast}d\sigma,
\end{eqnarray*}
with
\begin{align*}
&\displaystyle\widetilde S_{ki}:=\int_{I_z}\, b(\cos\theta,z)\, \psi_{k}(z)\psi_{i}(z)\, \pi(z)dz, \\
&\displaystyle\widetilde \nu_{ki}:=\int_{\mathbb R^d\times\mathbb S^{d-1}}\, \widetilde S_{ki}\, \phi(|v-v_{\ast}|)\, \mathcal M(v_{\ast})\, dv_{\ast}d\sigma,  \\[2pt]
&\displaystyle \qquad\text{and}\qquad\quad S_{kij}:=\int_{I_z}\,  b(\cos\theta,z)\, \psi_{k}(z)\psi_{i}(z)\psi_{j}(z)\, \pi(z)dz.
\end{align*}

For technical reasons, we assume $z\in I_z$ is one dimensional and $I_z$ has finite support $|z|\leq C_z$ (which is the case, for example, for the uniform and Beta distribution). Define
\[
\|h\|^2_{H^s_{x,v}}=\sum_{|j|+|l|\le s} \| \partial_l^j h\|^2_{L^2_{x,v}}\,, \quad
\|h\|^2_{H_z^s}=\int_{I_z} \| h \|^2_{H^s_{x,v}} \pi(z) dz\,.
\]
In \cite{LiuJin18}\cite{DausJL} the following results are given, under some suitable assumption on $b$:

\begin{theorem}
\label{thm3}
Assume the collision kernel $B$ is given by
\[
B(|v-v_*|), \cos\theta, z)=\phi(|v-v_*|)b(\cos\theta, z), \quad 
\phi(\xi)=C_\phi \xi^\gamma \quad {\text {with}} \gamma\in [0,1], C_\phi>0
\]
\[
\forall \eta\in [-1,1]. \qquad  |\partial_\eta b(\eta, z)|\le C_b, \quad |\partial_z^k b(\eta, z)| \le C_b^*, \quad \forall \, 0\le k\le r,
\]
where $b$ is linear in $z$, given in the form
\begin{equation}\label{b1} b(\cos\theta, z)= b_0(\cos\theta)+  b_1(\cos\theta)z\,.  \end{equation}
 Assume some upper and positive lower boundedness on $b$ and its derivatives.
In addition, assume \cite{JinShu17}
\begin{equation}\label{basis}
||\psi_k||_{L^{\infty}} \leq C k^p, \qquad \forall\, k,  \end{equation}
with a parameter $p>0$. Let $q>p+2$, define the energy $E^{K}$ by
\begin{equation}\label{E^K}
  E^{K}(t) = E_{s,q}^{K}(t) = \sum_{k=1}^{K}\, ||k^q h_k||_{H_{x,v}^s}^2,
\end{equation}
with the initial data satisfying $E^{K}(0) \leq \eta$.
Then for all $s\geq s_0$, $0\leq \varepsilon_d\leq 1$, such that
for $0\leq \varepsilon\leq \varepsilon_d$, if $h^K$ is an SG solution  (\ref{h_gPC1}) in
$H_{x,v}^s$, then: \\
(i) Under the incompressible Navier-Stokes scaling ($\alpha=1$),
$$E^{K}(t)\leq \eta \, e^{-\tau t}\,. $$
(ii) Under the acoustic scaling ($\alpha=0$),
$$E^{K}(t) \leq \eta \, e^{-\varepsilon \tau t}\,, $$
where $\eta$, $\tau$ are all positive constants that only depend on $s$ and $q$, independent of $K$ and $z$.
\end{theorem}
%----------------------------------------------------------------------------------------

From here, one also concludes that, $\displaystyle||h^K||_{H_{x,v}^s L_z^{\infty}}$  decays exponentially in time, with the same rate as $E^K(t)$, namely
\begin{equation}\label{h_K2}||h^K||_{H^s_{x,v}L_z^{\infty}} \leq \eta\,  e^{-\tau t}\end{equation}
in the incompressible Navier-Stokes scaling, and
$$ ||h^K||_{H^s_{x,v}L_z^{\infty}} \leq \eta\,  e^{-\varepsilon\, \tau t}$$ in the acoustic scaling.

\cite{LiuJin18} also gives the following error estimates on the SG method.

\begin{theorem}
\label{thm4}
Suppose the assumptions on the collision kernel and basis functions in Theorem \ref{thm3} are satisfied, and the initial data are the same as those in
Theorem \ref{thm3}, then \\
(i)\, Under the incompressible Navier-Stokes scaling,
\begin{equation}\label{Thm4_1} ||h-h^K||_{H^s_z} \leq C_{e}\,  \frac{e^{-\lambda t}}{K^r}, \end{equation}
(ii)\, Under the acoustic scaling, \begin{equation}\label{Thm4_2}||h-h^K||_{H^s_z} \leq C_{e}\,  \frac{e^{-\varepsilon\lambda t}}{K^r}\,, \end{equation}
with the constants $C_{e}, \,\lambda>0$ independent of $K$ and $\varepsilon$.
\end{theorem}

The above results show that the regularity of the SG solutions is the same as the initial data. Furthermore,  the numerical fluctuation
$h^K$ converges to $h$ with spectral accuracy, and the numerical error
will  decay exponentially in time in the random space.

%%%%%%%%%%%%%%%%%%%%%%%%%%%%%%%%%%%%%%%%%%%%%%%%%%%%%%%%%%%%%%%%%%%%%%%%%%%%%%%%%%%%%%%
\subsection{Asymptotic-preserving neural network approximation}

%\section{Asymptotic-Preserving Neural Networks}

%For the classical Asymptotic-preserving schemes, we focus on the %properties of numerical discretizations as the scaling parameter %changes.
%AP schemes mean to achieve the goal of perserving the numerical %discretizations of the limiting equations as well as uniform %stability.
%Though the design of AP schemes requires understandings of the %specific equation itself and is often case by case, there are some %general ideas that one can follow to construct the schemes %especially kinetic equations.

%However, when we design an APNN method, things change dramatically.
%Using DNNs to solve PDEs follows totally different strategies.
%In this paper, we use the PINNs framework to present the idea of %building an APNN method, however, it is not restricted to only %PINNs and can be generalized to DeepRitz, DeepONet, etc.

Kinetic equations have curse-of-dimensionality since it solves equations in the phase space. While this survey mainly concentrates on dealing with multiscale issues, it will be interesting to also dealt with the issue of high dimensionality together with multiple scales.
To this aim the deep neural network (DNN) offers a possible direction, since there have been examples in which DNNs offer
some advantages for high dimensional PDEs  \cite{E2018}\cite{raissi2019physics}\cite{lu2021learning}\cite{li2020fourier}.

Unlike classical numerical schemes,  a neural network uses a non-polynomial approximation  to approximate the training data through an optimization of an  empirical loss/risk. For multiscale
kinetic equations it is essential to construct a neural network that is AP \cite{Yang-APDNN} (referred to APNNs).
%\begin{figure}[ht]
%    \centering
%    \includegraphics[width=0.45\textwidth]{Figure/apnns.png}
%    \caption{Illustration of APNNs. $\mathcal{F^{\varepsilon}}$ is the %microscopic equation that depends on the small scale parameter $\varepsilon$ %and $\mathcal{F}^{0}$ is its macroscopic limit as $\varepsilon \to 0$, which %is independent of $\varepsilon$. The latent solution of %$\mathcal{F^{\varepsilon}}$ is approximated by neural networks with its %measure denoted by $\mathcal{R}(\mathcal{F^{\varepsilon}})$. The asymptotic %limit of $\mathcal{R}(\mathcal{F^{\varepsilon}})$ as $\varepsilon \to 0$, if %exists, is denoted by $\mathcal{R}(\mathcal{F}^{0})$. If %$loss(\mathcal{F}^{0})$ is a good measure of $\mathcal{F}^{0}$, then it is %called asymptotic-preserving (AP).}
%    \label{fig:apnns}
%\end{figure}
%By comparing this with the traditional AP diagram \cite{GJL, %jin2010asymptotic}, we can think the loss/risk $\mathcal{R}$ as the %discretization or schemes used in the traditional numerical schemes %and think $\mathcal{R}\to 0$ line, that is, the training loss tends %to $0$, as the numerical consistency.

%\subsection{Solving the linear transport equation by APNNs}

We first introduce  conventional notations for deep neural networks (DNNs). An $L$-layer feed forward neural network is defined recursively as,
\begin{equation}
    \begin{aligned}
        f_{\theta}^{[0]}(x) & = x,                                                                              \\
        f_{\theta}^{[l]}(x) & = \sigma \circ (W^{[l-1]} f_{\theta}^{[l-1]}(x) + b^{[l-1]}), \, 1 \le l \le L-1, \\
        f_{\theta}(x)       & = f_{\theta}^{[L]}(x) = W^{[L-1]} f_{\theta}^{[L-1]}(x) + b^{[L-1]},
    \end{aligned}
\end{equation}
where $W^{[l]} \in  \mathbb{R}^{m_{l+1}\times m_l}, b^{l}\in  \mathbb{R}^{m_{l+1}}, m_0 = d_{in} = d$ is the input dimension, $m_{L} = d_0$ is the output dimension, $\sigma$ is a scalar function and $"\circ"$ means entry-wise operation. We denote the set of parameters by $\theta$.  The layers are denoted by a list, i.e., $[m_0, \cdots, m_L]$.

Consider the linear transport equation 
 with initial and boundary conditions over a bounded domain $ \mathcal{T} \times \mathcal{D} \times \Omega$:
\begin{equation}\label{eq:linear-trasport-eqn}
    \left\{
    \begin{array}{ll}
        \varepsilon \partial_t f + v \cdot \nabla_x f = \frac{1}{\varepsilon} \mathcal{L} f, & (t, x, v) \in \mathcal{T} \times \mathcal{D} \times \Omega,\,          \\
        \mathcal{B}f = F_{\text{B}},                                                                              & (t, x, {v}) \in \mathcal{T} \times \partial \mathcal{D} \times \Omega,\, \\
        \mathcal{I}f = f_{0},                                                                                     & (t, {x}, {v}) \in \{t = 0\} \times \mathcal{D} \times \Omega,
    \end{array}
    \right.
\end{equation}
where $F_{\text{B}}, f_{0}$ are given functions; $\partial \mathcal{D}$ is the boundary of $\mathcal{D}$, and $\mathcal{B}, \mathcal{I}$ are initial and boundary operators, respectively.
$\mathcal{L}=\sigma (\rho -f)$.

\subsubsection{The failure of PINNs  to resolve small scales}
 PINNs is a standard neural network to solve PDEs. There  the density function $f(t, x, v)$ is approximated by a neural network
\begin{equation}
    \text{NN}_{\theta}(t, x, v) \approx f(t, x, v).
\end{equation}
The inputs of DNN are $(t, {x}, {v})$, i.e., $m_0 = 3, 5$ for $1$-d and $2$-d respectively. The output is a scalar which represents the value of $f$ at $(t, x, v)$. Since $f$ is always non-negative, we put an exponential function at the last output layer of the DNN:
\begin{equation}
    f^{\text{NN}}_{\theta}(t, x, v) := \exp \left(-\tilde{f}^{\text{NN}}_{\theta}(t, x, v)\right) \approx f(t, x, v)
\end{equation}
to represent the numerical solution of $f$. Then the least square of the residual of the original transport equation~\eqref{eq:linear-trasport-eqn} is used as the target loss function, together with boundary and initial conditions as penalty terms, 
\begin{equation}\label{eq:loss-pinn}
    \begin{aligned}
        \mathcal{R}^{\varepsilon}_{\text{PINN}} = & \frac{1}{|\mathcal{T} \times \mathcal{D} \times \Omega|} \int_{\mathcal{T}} \int_{\mathcal{D}} \int_\Omega \left| \varepsilon^2 \partial_t f^{\text{NN}}_{\theta} + \varepsilon {v} \cdot \nabla_x f^{\text{NN}}_{\theta} - \mathcal{L} f^{\text{NN}}_{\theta}  \right|^2 d{v} d{{x}} \,d{t} \\
                                                  & +  \frac{\lambda_1}{|\mathcal{T} \times \partial \mathcal{D} \times \Omega|}  \int_{\mathcal{T}} \int_{\partial \mathcal{D}} \int_\Omega |\mathcal{B}f^{\text{NN}}_{\theta} - F_{\text{B}}|^2 dv\,dx\,dt                                                                                                 \\
                                                  & +  \frac{\lambda_2}{|\mathcal{D} \times \Omega|} \int_{\mathcal{D}} \int_\Omega |\mathcal{I}f^{\text{NN}}_{\theta} - f_{0}|^2 dv\, dx,
    \end{aligned}
\end{equation}
where $\lambda_1$ and $\lambda_2$ are the penalty weights  to be tuned. Then a standard stochastic gradient method (SGD) or Adam optimizer is used to find the global minimum of this loss.

Now let us check whether this PINN method is AP.  One only needs to focus on the first term of~\eqref{eq:loss-pinn}
\begin{equation}
    \mathcal{R}^{\varepsilon}_{\text{PINN, residual}} := \frac{1}{|\mathcal{T} \times \mathcal{D} \times \Omega|} \int_{\mathcal{T}} \int_{\mathcal{D}} \int_\Omega \left| \varepsilon^2 \partial_t f^{\text{NN}}_{\theta} + \varepsilon {v} \cdot \nabla_x f^{\text{NN}}_{\theta} - \mathcal{L} f^{\text{NN}}_{\theta}  \right|^2 dv\, dx\, dt.
\end{equation}
Taking $\varepsilon \to 0$, formally this will lead to
\begin{equation}
    \mathcal{R}_{\text{PINN, residual}} := \frac{1}{|\mathcal{T} \times \mathcal{D} \times \Omega|} \int_{\mathcal{T}} \int_{\mathcal{D}} \int_\Omega \left| - \mathcal{L} f^{\text{NN}}_{\theta} \right|^2 d{{v}}\, dx\,dt,
\end{equation}
which can be viewed as the PINN loss of the equilibrium equation
\begin{equation}
    \mathcal{L} f = 0.
\end{equation}
This shows that when $\varepsilon$ is very small, to the leading order we are solving equation $\mathcal{L} f = 0$ which gives
$f=\rho$. This does not give the desired   diffusion equation~\eqref{televen}. This explains why PINN will fail when $\varepsilon$ is small.

 The APNN, introduced in \cite{JinMaWu},  puts the micro-macro system~\eqref{mM-1} or \eqref{mM-2} into the loss, instead of the original equation~\eqref{eq:linear-trasport-eqn}.

First the DNN needs to parametrize two functions $\rho(x,v)$ and $g(t, x, v)$. So here two networks are used. First
\begin{equation}
    \rho^{\text{NN}}_{\theta}(t, x) := \exp \left( -\tilde{\rho}^{\text{NN}}_{\theta}(t, x)\right) \approx \rho(t, x).
\end{equation}
Notice here $\rho$ is non-negative. Second,
\begin{equation}\label{g-net}
    g^{\text{NN}}_{\theta}(t, x, v) := \tilde{g}^{\text{NN}}_{\theta}(t, x, v) - \langle \tilde{g}^{\text{NN}}_{\theta} \rangle (t, x) \approx g(t, x, v).
\end{equation}
Here $\tilde{\rho}$ and $\tilde{g}$ are both fully-connected neural networks. Notice that by choosing $g^{\text{NN}}_{\theta}(t, x, v)$ as in ~\eqref{g-net} it will automatically satisfy the constraint
\begin{equation}\label{g-cons}
    \langle g \rangle = 0\,,
\end{equation}
because
\begin{equation}
    \langle g^{\text{NN}}_{\theta} \rangle = \langle \tilde{g}^{\text{NN}}_{\theta} \rangle - \langle \tilde{g}^{\text{NN}}_{\theta} \rangle = 0, \quad \forall \; t, x.
\end{equation}

Then the  APNN loss is defined as
\begin{equation}\label{eq:loss-ap}
    \begin{aligned}
        \mathcal{R}^{\varepsilon}_{\text{APNN}} = & \frac{1}{|\mathcal{T} \times \mathcal{D}|} \int_{\mathcal{T}} \int_{\mathcal{D}} | \partial_t \rho^{\text{NN}}_{\theta} + \nabla_x \cdot \left \langle   {v} g^{\text{NN}}_{\theta} \right \rangle  |^2 dx\,dt                                     \\
                                                  & + \frac{1}{|\mathcal{T} \times \mathcal{D} \times \Omega|} \int_{\mathcal{T}} \int_{\mathcal{D}} \int_\Omega | \varepsilon^2 \partial_t g^{\text{NN}}_{\theta}  + \varepsilon (I - \Pi)({v} \cdot \nabla_x g^{\text{NN}}_{\theta})                                      \\
                                                  & \quad  + {v} \cdot  \nabla_{{x}} \rho^{\text{NN}}_{\theta} - {\mathcal L} g^{\text{NN}}_{\theta}|^2 dv\,dx\,dt                                                                                                      \\
                                                  & +  \frac{\lambda_1}{\mathcal{T} \times\partial \mathcal{D} \times \Omega|}  \int_{\mathcal{T}} \int_{\partial \mathcal{D}} \int_\Omega |\mathcal{B}(\rho^{\text{NN}}_{\theta} + \varepsilon g^{\text{NN}}_{\theta}) - F_{\text{B}}|^2 dv\,dx\,dt \\
                                                  & +  \frac{\lambda_2}{|\mathcal{D} \times \Omega|} \int_{\mathcal{D}} \int_\Omega |\mathcal{I}(\rho^{\text{NN}}_{\theta} + \varepsilon g^{\text{NN}}_{\theta}) - f_{0}|^2 dv\,dx.
    \end{aligned}
\end{equation}

\begin{figure}[ht]
    \centering
    \includegraphics[width=0.8\textwidth]{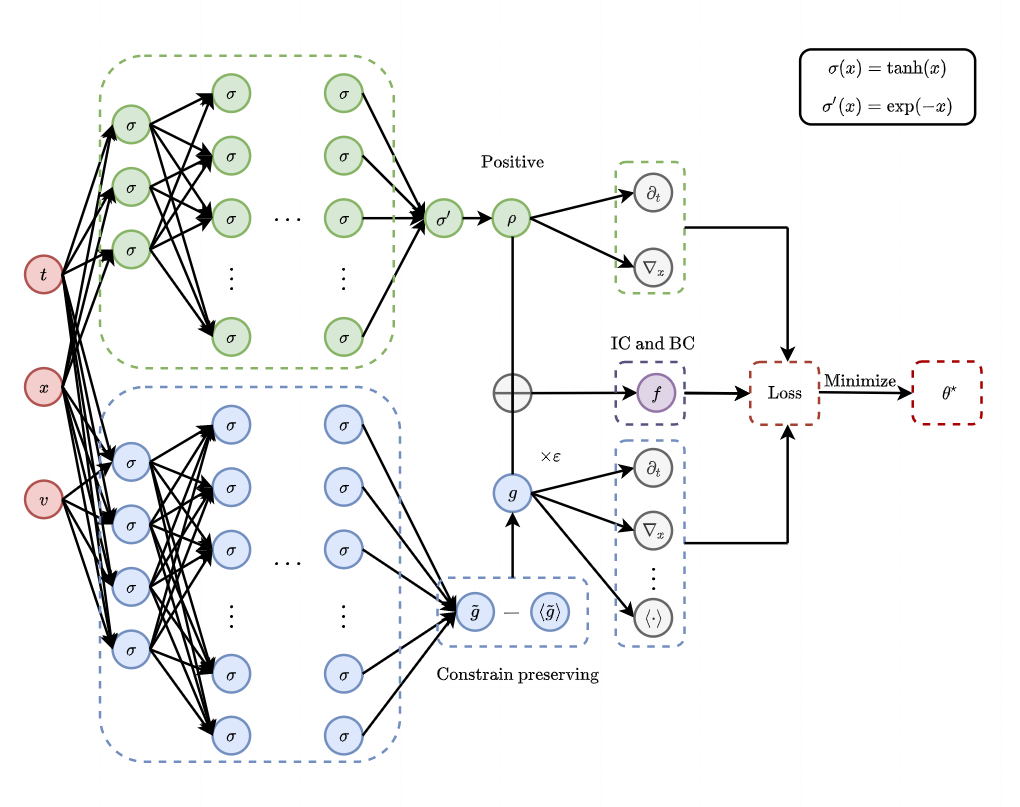}
    \caption{Schematic of APNNs for solving the linear transport equation with initial and boundary conditons.}
    \label{fig:APNNs}
\end{figure}

A schematic plot of the method is given in Figure~\ref{fig:APNNs}.

To show formally the AP property of  this loss, one only needs to focus on the first two terms of~\eqref{eq:loss-ap}
\begin{equation}
    \begin{aligned}
        \mathcal{R}^{\varepsilon}_{\text{APNN, residual}} = & \frac{1}{|\mathcal{T} \times \mathcal{D}|} \int_{\mathcal{T}} \int_{\mathcal{D}} | \partial_t \rho^{\text{NN}}_{\theta} + \nabla_x \cdot \left \langle   {v} g^{\text{NN}}_{\theta} \right \rangle - Q |^2 dx\,dt                                                               \\
                                                            & + \frac{1}{|\mathcal{T} \times \mathcal{D} \times \Omega|} \int_{\mathcal{T}} \int_{\mathcal{D}} \int_\Omega \Big | \varepsilon^2 \partial_t g^{\text{NN}}_{\theta} + \varepsilon (I - \Pi)({v} \cdot \nabla_x g^{\text{NN}}_{\theta})  \\
                                                            &  \qquad + {v} \cdot  \nabla_{{x}} \rho^{\text{NN}}_{\theta}- {\mathcal L} g^{\text{NN}}_{\theta}-(I -\Pi) \varepsilon Q \Big|^2 dv\,dx\,dt\,.
    \end{aligned}
\end{equation}
Taking $\varepsilon \to 0$, formally this will lead to
\begin{equation}
    \begin{aligned}
        \mathcal{R}_{\text{APNN, residual}} = & \frac{1}{|\mathcal{T} \times \mathcal{D}|} \int_{\mathcal{T}} \int_{\mathcal{D}} | \partial_t \rho^{\text{NN}}_{\theta} + \nabla_x \cdot \left \langle   {v} g^{\text{NN}}_{\theta} \right \rangle - Q |^2 dx\,dt                           \\
                                              & + \frac{1}{|\mathcal{T} \times \mathcal{D} \times \Omega|} \int_{\mathcal{T}} \int_{\mathcal{D}} \int_\Omega \Big | {v} \cdot  \nabla_{{x}} \rho^{\text{NN}}_{\theta} - {\mathcal L} g^{\text{NN}}_{\theta} \Big|^2 dv\,dx\,dt,
    \end{aligned}
\end{equation}
which is the least square loss of equations
\begin{equation}
    \left\{
    \begin{aligned}
         & \partial_t \rho + \nabla_{{x}} \cdot \left \langle  {{v}}g \right \rangle = Q, \\
         & {v} \cdot  \nabla_{{x}} \rho = \mathcal{L} g .
    \end{aligned}
    \right.
\end{equation}
The second equation above yields $g=\mathcal{L}^{-1}({v} \cdot  \nabla_{{x}} \rho)$, which, when plugging into the first equation and integrating over $v$, gives  the diffusion equation~\eqref{televen}.
Hence this proposed method is an APNN method.

\begin{example}\label{ex-DNN}

We present a numerical example from \cite{JinMaWu}. Let $\sigma=1$. Consider 
the initial data as follows
\begin{equation}
    f_0(x, v) = \frac{\rho(x)}{\sqrt{2\pi}}e^{-\frac{v^2}{2}},
\end{equation}
where
\begin{equation}
    \rho(x) = 1 + \cos (4 \pi x),
\end{equation}
and the isotropic in-flow boundary conditions:
\begin{equation}
     F_L(v) = 1, \quad F_R(-v) = 0, \quad {\text {for}} \, v>0.
\end{equation}
The results are shown in Figure \ref{fig:ex5}. Clearly PINN fails to conserve the mass, and  for small $\varepsilon$  does not give accurate results, while APNN gives quite accurate results even when $\e$ is very small.  

\begin{figure}[ht]
    \centering
    \includegraphics[width=0.8\textwidth]{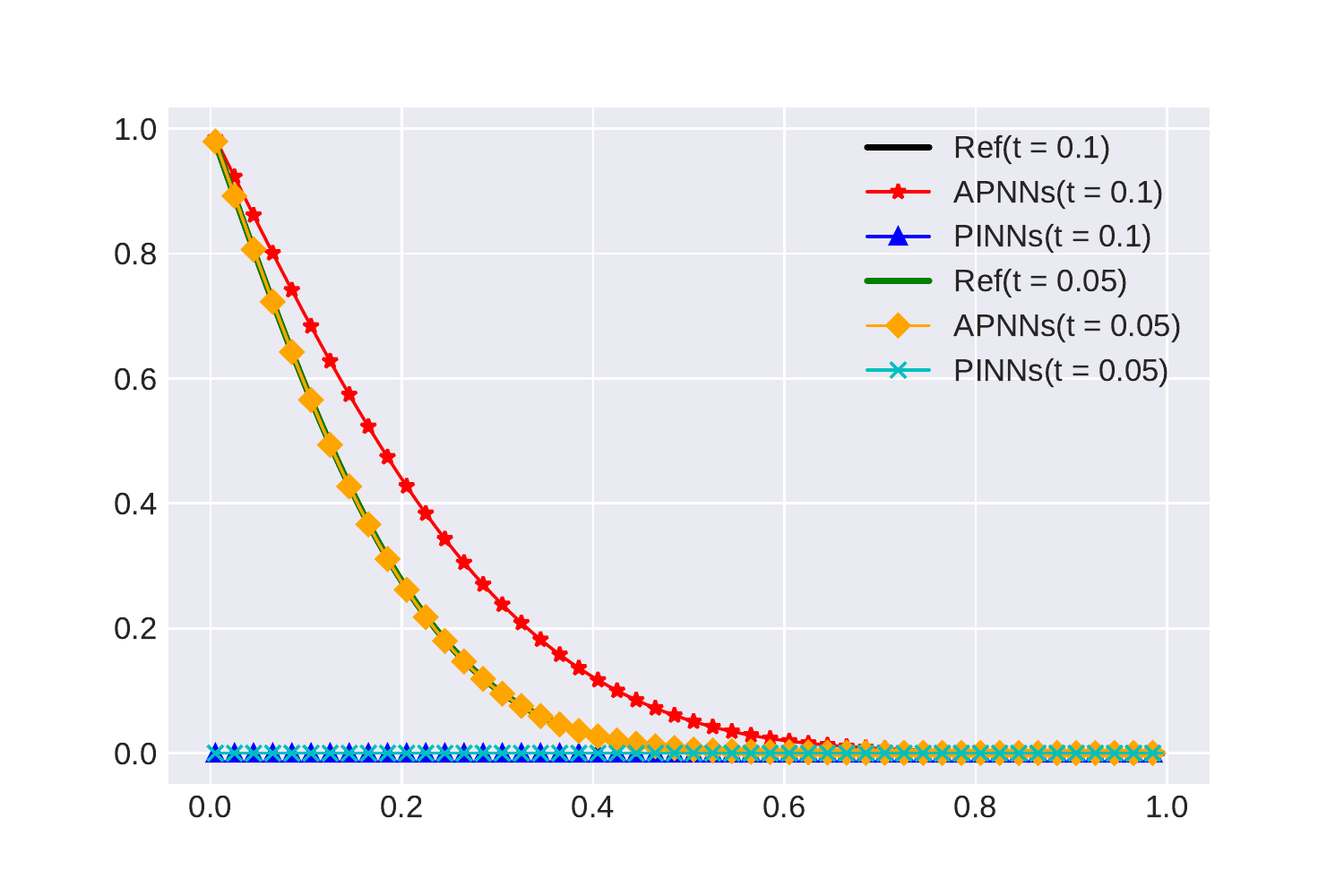}
    \caption{Example \ref{ex-DNN}. $\varepsilon = 10^{-8}$. Plots of density $\rho$ at $t = 0.05, 0.1$ by PINN, APNN and  Reference solutions. The neural networks are $[2, 128, 128, 128, 128, 1]$ for $\rho$ and $[3, 256, 256, 256, 256, 1]$ for $g, f$. Batch size is $1000$ in domain, $400 \times 2$ with penalty $\lambda_1 = 10$ for boundary condition and $1000$ with penalty $\lambda_2 = 10$ for initial condition, the number of quadrature points is $30$. Relative $\ell^2$ errors of PINNs and APNNs are $9.40 \times 10^{-1}, 2.76  \times 10^{-3}$ respectively.}
    \label{fig:ex5}
\end{figure}
\end{example}

\begin{remark}\label{rmk2} Not all AP schemes yield an APNN method when put into the loss. For example it was shown in \cite{JinMaWu} the  parity formulation (\ref{4.7})-(\ref{thsix}) does not give an APNN network. 

A similar loss function also based on micro-macro decomposition, but with a constraint (\ref{g-cons})
as a penalty, was proposed in \cite{lu2021solving} for stationary problem.
\end{remark}

\section{Other related multiscale problems}\label{sec5}

\subsection{Nonlinear hyperbolic systems with stiff source terms}

Numerical methods for nonlinear hyperbolic systems with stiff relaxation 
terms were among the earliest AP schemes for time-dependent problems.
A prototype equation is given by:
\begin{align} 
\label{Hyp-1}
\left\{
\begin{array}{l}
\displaystyle \partial_tu + \partial_x g(u,v)=0, \\[6pt]
\displaystyle \partial_tv+\partial_x h(u,v) =\frac{1}{\varepsilon}R(u,v),
\end{array}\right.
\end{align}
where the term $R$ is dissipative: $\partial_vR\leq 0$ and possesses a unique local equilibrium: $R(u,v)=0$, which implies $v=f(u)$. Then when $\varepsilon\rightarrow 0$, one has the macroscopic limit
\begin{equation*}
\label{relax-limit}
\partial_tu+\partial_xg(u,f(u))=0.
\end{equation*}
This is an analogy of the
Euler limit of the Boltzmann equation.  If one uses the Chapman-Enskog
expansion to $O(\varepsilon)$ term, then (\ref{Hyp-1}) can be approximated by
a parabolic equation
\begin{equation}
\label{relax-NS}
\partial_t u + \partial_x g(u, f(u)) = \varepsilon \partial_x[h(u)\partial_x
u]\,.
\end{equation}
Here one needs $h(u)>0$, which requires the characteristic speeds of the
original system (\ref{Hyp-1}) interwine with that of the limiting equation
(\ref{relax-limit}). This condition  is called the {\it subcharacteristic condition}
\cite{CLL}.  (\ref{relax-NS}) is an analogy of the Navier-Stokes approximation
to the Boltzmann equation.

Numerical study of system of the type (\ref{Hyp-1}) began
in the works \cite{jin1996numerical} \cite{jin1995runge}\cite{caflisch1997uniformly}, where the AP principle was applied
to design numerical schemes to handle the stiff relaxation term. High order IMEX type schemes were
developed in \cite{PR-IMEX}\cite{Toro2008}. Combining AP and positivity preserving property was done in
\cite{HuShu-APPP}. The relation between AP and well-balanced scheme is revealed in \cite{GosseTos}.
For AP schemes for gas dynamics with external force and frictions see
\cite{Bouchut-Perthame}\cite{Chalons-10}. A rigorous uniform accuracy proof of AP schemes for linear problems  was recently made in \cite{hu2021uniform}.

\subsection{Quasi-neutral limit in plasma}

In many plasma applications, one can disregard charge separations, and then the quasi-neutral approximation can be used. 
However,  near the plasma boundary,  electrostatic sheathes may appear, then one needs to
consider more complex models \cite{degond-JCPReview}.

Consider the one-species recaled  Euler-Poisson (EP) equations for charged 
particles:
\begin{eqnarray}
&&\partial_t n + \nabla \cdot q=0\,,
\label{eq:mass}
\\
&& \partial_t q + \nabla \cdot \left( \frac{q \otimes q}{n}\right)+\nabla p(n)
=n 
\nabla \phi\,,
\label{eq:momentum}
\\
&& \varepsilon^2 \Delta \phi=n-1\,,
\label{Poisson}
\end{eqnarray}
 where  $n=n(x,t)$ is the particle number density, 
 $q=q(x,t)=nu$ is the momentum ($u$ is the average velocity),
$p(n)=n^\gamma$ is the pressure law with $\gamma \ge 1$, and $\phi=\phi(x,t)$ 
is the electric potential.
Here the negatively charged electrons with scaled charge equal to
$-1$ is considered, with a
uniform ion background density equal to 1.
The dimensionless parameter $\varepsilon=\lambda_D/L$ is the scaled Debye length,
i.e., the ratio of the actual Debye length $\lambda_D$ to the macroscopic length
scale $L$. 
In the quasineutral regime,  $\varepsilon\ll 1$. 
 Letting $\varepsilon \to 0$ in (\ref{Poisson}), one has $n=1$, 
and the following
limiting equations arise \cite{Brenier-QN}:
\begin{eqnarray}
& \nabla \cdot q=0\,,
\label{div-free}
\\
\label{div-free-1}
& \partial_t q + \nabla \cdot ( q \otimes q)= \nabla \phi\,.
\label{momentum-incomp}
\end{eqnarray}
This is the incompressible Euler equations.

A typical 
stable time discretization of the Euler-Poisson system requires
$\Delta t \le \varepsilon$, which is quite restrictive. 
%\label{s-cond}

The main difficulty here is that when $\e \to 0$,  the Poisson equation (\ref{Poisson}) becomes degenerate, hence a naive  discretization would lead to poor numerical performance for small $\e$.
 A key idea introduced by 
 Degond etc. in \cite{CDV1} is to reformulate the system to a new one that
remains {\it uniformly elliptic}.  
 Taking $\partial_t$ on (\ref{eq:mass}), $\nabla \cdot$ on
(\ref{eq:momentum}) and $\partial_{tt}$ on (\ref{Poisson}) give
\begin{eqnarray}
&&
\partial_{tt}n + \nabla \cdot \partial_t q=0\,,
\label{n-eq}
\\
&& \nabla \cdot \partial_t q + \nabla^2: \left(\frac{q\otimes q}{n} + p(n) {\rm I} \right)
=\nabla \cdot (n\nabla \phi)\,,
\label{m-eq}
\\
&& 
\varepsilon^2 \Delta \partial_{tt} \phi=\partial_{tt} n\,.
\label{phi-eq}
\end{eqnarray}
Eliminating $\nabla \cdot \partial_t q$ by combining (\ref{n-eq}) and (\ref{m-eq}) and using
(\ref{phi-eq}), one gets
\begin{equation}
 - \nabla \cdot [(n+\varepsilon^2 \partial_{tt}) \nabla \phi] + \nabla^2 : 
\left[\frac{q\otimes q}{n} + p(n) {\rm I}\right] 
=0\,.
\label{2.7}
\end{equation}
Although this system is equivalent to the original Euler-Poisson system,  Equation (\ref{2.7}) now is  {\it uniformly elliptic} in $\e$, discretizing it suitably in time   will guarantee
 the asymptotic stability with respect to $\varepsilon$.
 
This framework is quite general, and has  been generalized to two-fluid model \cite{CDV2}, 
Particle-in-Cell method for Vlasov-Poisson system \cite{DDNSV}, Euler-Maxwell system \cite{Degond-EM}, and Vlasov-Maxwell system \cite{Degond-VM}, among
other plasma models. See a recent comprehensive review \cite{Degond-AP}.

\subsection{High-field limits}

In kinetic equations, often there is a strong external field, such as the
electric or magnetic field, that  balances the collision
term, leading to the so-called high field limit \cite{CGL97}.  

\subsubsection{High electric field}

Consider for
example  the interaction between the electrons 
and a surrounding bath through Coulomb force in electrostatic plasma, where the
electron density distribution $f(t,x,v)$ is governed by the Vlasov-Poisson-Fokker-Planck system:
\begin{eqnarray}
&&\partial_t f + v\cdot \nabla_x f - \frac{1}{\varepsilon}\nabla_x  \phi \cdot \nabla_v f =\frac{1}{\varepsilon} \nabla_v \cdot (vf + \nabla_v f) , \label{VPFPNondim1}
\\&& -\Delta \phi = \rho - h, \label{VPFPNondim2}
 \label{VPFPNondim}
\end{eqnarray}
where $\varepsilon = \left(\frac{l_e}{\Lambda}\right)^2$ is the ratio
between the mean free path and the Debye length.
Let $\varepsilon \rightarrow 0$ in (\ref{VPFPNondim1}), 
 one obtains the high-field limit \cite{nieto2001high}\cite{GNPS}:
\begin{eqnarray}
&&\partial_t\rho - \nabla_x \cdot(\rho\nabla_x \phi) = 0,
\\&& -\Delta \phi = \rho - h(x).
\label{HFL}
\end{eqnarray}

One can combine the force term with the
Fokker-Planck term as
\begin{equation}
\partial_t f + v\cdot \nabla_x  f = \frac{1}{\varepsilon}\nabla_v \cdot
\left[M \nabla_v (Mf)\right]
\label{VPFPNondim12}
\end{equation}
where $M=e^{-|v+\nabla_x \phi|^2/2}$. This form is convenient for designing AP schemes \cite{JinWang11}\cite{crouseilles2011asymptotic}, based on which one can easily use  other
well-developed AP frameworks.  For more general collision 
operator, for example the semiconductor Boltzmann collision operator, this 
trick does not apply and one needs other ideas, for example the BGK penalization \cite{JinWang13}. 
A variational approach was recently proposed in \cite{CWXY}, using 
the Wasserstein gradient structure, to get  positivity and AP easily. 

\subsubsection{High magnetic field}

Magnetized plasmas are encountered in a wide variety of astrophysical situations, but also in magnetic fusion devices such as Tokamaks, where a large external magnetic field needs to be applied in order to keep the particles on the desired tracks. The Vlasov equation for such problems takes the following form:
\begin{equation*}
\displaystyle{\varepsilon \partial_t f+ v \cdot\nabla_x f \,+\,\left(
  {E}(t,x) \,+\, \frac{1}{\eps}v \wedge B_{\rm ext}(t,x) \right)\cdot\nabla_v f
\,=\, 0.}
\end{equation*}
Here, for simplicity we set all physical constants to one and consider that $\varepsilon>0$ is a small parameter related to the ratio between the reciprocal Larmor frequency and the advection time scale. We consider a constant magnetic field acting in the vertical $z$-direction, hence it yields the two-dimensional Vlasov-Poisson system with an external strong force:
\begin{equation*}
%\label{eq:vlasov2d}
\left\{
\begin{array}{l} 
\displaystyle{\varepsilon \partial_t f+v\cdot\nabla_x f \,+\,\left(
  {E}(t,x) \,-\, \frac{b\,v^\perp}{\eps} \right)\cdot\nabla_v f
\,=\, 0,}
\\
\,
\\
\displaystyle{E = - \nabla_x \phi, \quad -\Delta_{x} \phi=
  \rho^\eps,\quad \rho=\int_{\RR^2} f dv,}
\end{array}\right.
\end{equation*}
where we use notation $v^\perp=(-v_y,v_x)$.

Most of the numerical schemes for the Vlasov-Poisson system are based on particle methods, which consist of approximating the distribution function by a finite number of macro-particles. The trajectories of these particles are determined from the characteristic curves corresponding to the Vlasov equation
\begin{equation}
\label{traj:00}
\left\{
\begin{array}{l}
\ds{\varepsilon\frac{\dd X^ \eps}{\dd t} \,=\, V^\eps,} 
\\
\,
\\
\ds{\varepsilon\frac{\dd V^\eps}{\dd t} \,=\, -\frac{b\,V^{\eps\perp}}{\varepsilon}  \,+\, E^\eps(t,X^\eps), }
\end{array}\right.
\end{equation}
where we use the conservation of $f$ along the characteristic curves  
$$
f(t, X^\eps(t),V^\eps(t)) = f(t^0, X^0, V^0).
$$
In the limit $\eps\to0$ one expects oscillations occurring on typical time scales $O(1/\eps^2)$ to coexist with a slow dynamics evolving on a time scale $O(1)$. We sketch now how to identify a closed system describing in the leading order the slow evolution.  To begin with, note that from the second line of system~\eqref{traj:00} it does follow that $V^\eps$ oscillates at order $1/\eps^2$ thus remains bounded and converges weakly\footnote{Though we do not want to be too precise here, let us mention that in the present discussion \emph{weakly} and \emph{strongly} refer to the weak-* and strong topologies of $L^\infty$ and that the weak convergences that we encounter actually correspond to strong convergence in $W^{-1,\infty}$.} to zero. As we detail below, one may also combine both lines of the system to obtain
$$
\frac{\dd}{\dd t}\left( X^\eps - \eps
  \,\frac{V^{\eps\perp}}{b}\right)  \,=\, -\frac{1}{b}\, E^{\perp}(t,X^\eps).
  $$
  This shows that $X^\eps$ evolves slowly but, as such, does not provide a closed asymptotic evolution in the limit $\eps\to 0$ and the corresponding asymptotic model is an equation for the density $\rho$
  $$
  \left\{
\begin{array}{l} 
\displaystyle {\partial_t \rho -\, b\,E^\perp\cdot\nabla_x \rho\,=\, 0,}
\\
\,
\\
\displaystyle{E = - \nabla_x \phi, \quad -\Delta_{x} \phi =
  \rho.}
\end{array}\right.
  $$

When $b$ is not constant, one also needs to know what happens to expressions that are quadratic in $V^\eps$ and this does not follow readily from the weak convergence of $V^\eps$  \cite{FR17}.
  
One of the oldest of these strategies \cite{FSS:09} is directly inspired by theoretical results on two-scale convergence and relies on the fact that at the limit $\eps\to0$ the $\tau$-dependence ($\tau=t/\e$) may be explicitly filtered out. Its main drawback is probably that it computes only the leading order term in the limit $\eps\to 0$. In particular it is only available when $\eps$ is very small. 

This may be fixed by keeping besides the stiff term to which a two-scale treatment is applied a non-stiff part that is smaller in the limit $\eps\to0$ but becomes important when $\eps$ is not small. Such a decomposition may be obtained by using a micro-macro
approach as in \cite{CFHM:13} and some references therein. This does allow to switch from one regime to another without any treatment of the transition between those but results in relatively heavy schemes,. 

Another approach with similar advantages, developed in \cite{CLM:13} and \cite{FHLS}, consists in explicitly doubling time variables and seeking higher-dimensional partial differential equations and boundary conditions in variables $(t,\tau,x,v)$ that contains the original system at the $\eps$-diagonal $(t,\tau)=(t,t/\eps)$. While the corresponding methods are extremely good at capturing oscillations their design requires a deep \emph{a priori} understanding of the detailed structure of oscillations. Also a class of semi-implicit schemes  have been proposed \cite{FR16,FR17}\cite{HLB}\cite{RC} to capture accurately the non stiff part of the evolution while allowing for coarse discretization parameters. It allows to capture the asymptotic limit of the two dimensional Vlasov-Poisson system with a \emph{uniform} magnetic field \cite{FR16}\cite{FRZ21}. In many respects those schemes are remarkably natural and simple and can be adapted to totroidal geometry as in \cite{FR20}.

\subsection{Highly anisotropic diffusion}

In magnetized plasma simulations, magnetic field confines the particles around the field lines, which leads to highly anisotropic problems. The model problem writes
\begin{equation}\label{eq:anisotropic}
\begin{aligned}
&-\nabla\cdot (A\nabla u)=f,\qquad&\mbox{on $\Omega$},\\
& \mathbf{n}\cdot A\nabla u=0,\qquad &\mbox{on $\Gamma_N$,}\\
& u=g,\qquad &\mbox{on $\Gamma_D$,}
\end{aligned}
\end{equation}
where $\Omega\subset\mathbb{R}^2$
 or $\Omega\subset\mathbb{R}^3$
 is a bounded domain with boundary $\partial \Omega = \Gamma_D \cup \Gamma_N$ and outward normal $\mathbf{n}$. The direction of the anisotropy is given by a unit vector field $\mathbf{b}$ and the anisotropic diffusion matrix is then given by
\begin{equation}\label{eq:anisotropicA}
A=\frac{1}{\e}A_{||}\mathbf{b}\otimes \mathbf{b}+(Id-\mathbf{b}\otimes\mathbf{b})A_{\perp}(Id-\mathbf{b}\otimes\mathbf{b})
\end{equation}
$A_{||}>0$ is a scalar and $A_\perp$ is a symmetric positive definite matrix field, both of them are of order one. The problem becomes highly anisotropic when $\epsilon\ll 1$. Let $(\xi, \eta)$ be the aligned coordinate system, the formal limit of $\epsilon\to 0$ leads to
\[\begin{aligned}
&-\partial_\xi (A_{||}\partial_\xi u)=f,\qquad&\mbox{on $\Omega$},\\
& \partial_{\xi} u=0,\qquad &\mbox{on $\Gamma_N$,}\\
& u=g,\qquad &\mbox{on $\Gamma_D$.}
\end{aligned}
\]
Any function that remains constant along the $\mathbf{b}$ field solves the above limit equation, which indicates that the limit model is not well-posed.
Due to the existence of infinitely many solutions when $\epsilon\to 0$, standard numerical discretizations suffer from large condition numbers when $\epsilon$ is small and usually lose convergence when $\e \ll h$ ($h$ is the mesh size).

To avoid the aforementioned problem, the common approach is to use magnetic field aligned coordinates, which may run into problems when
there are magnetic re-connections or highly fluctuating field directions. The other approach is to design methods whose condition numbers do not scale with the anisotropy
strength and the convergence orders are uniform with respect to $\e$. The main idea is to construct new systems that keep well-posed when $\e\to 0$. 
In a series of papers by Degond, Narski, Negulescu, {\it et. al.} \cite{degond2010asymptotic} \cite{degond2012asymptotic}\cite{degond2010duality}\cite{narski2014asymptotic}, various reformulations and discretization strategies based on macro-micro decomposition are proposed. Another idea is based on field line integration.  By substituting \eqref{eq:anisotropicA} into \eqref{eq:anisotropic}, one gets
\begin{equation}\label{eq:anisotropicxi}
\begin{aligned}
&-\partial_{\xi}\big(\frac{1}{\epsilon}A_{||}\partial_{\xi} u\big)-\nabla\cdot\Big((Id-\mathbf{b}\otimes\mathbf{b})A_{\perp}(Id-\mathbf{b}\otimes\mathbf{b})\nabla u\Big)=f,\qquad&\mbox{on $\Omega$},\\
& \frac{1}{\epsilon}A_{||}\partial_{\xi} u+
\mathbf{n}\cdot(Id-\mathbf{b}\otimes\mathbf{b})A_{\perp}(Id-\mathbf{b}\otimes\mathbf{b})\nabla u=0,\qquad &\mbox{on $\Gamma_N$,}\\
& u=g,\qquad &\mbox{on $\Gamma_D$,}
\end{aligned}
\end{equation}
Taking the integration of the first equation in Eq. \eqref{eq:anisotropicxi} along a field line and using the second equation of the Neumann boundary condition, one can get 
\begin{equation}\label{eq:anisotropiclimit}
\begin{aligned}
&\int_0^{L}\Big(\nabla\cdot\Big((Id-\mathbf{b}\otimes\mathbf{b})A_{\perp}(Id-\mathbf{b}\otimes\mathbf{b})\nabla u\Big)+f\Big)\mathrm{d}\xi\\
=&
\Big(\mathbf{n}\cdot(Id-\mathbf{b}\otimes\mathbf{b})A_{\perp}(Id-\mathbf{b}\otimes\mathbf{b})\nabla u\Big)\Big|_{0}^{L},
\end{aligned}\end{equation}
where $0$ and $L$ are the two end points of a field line. \eqref{eq:anisotropiclimit} is an equation independent of $\epsilon$ and provides the information of how to determine the limit solution. Similar idea can be extended to more complex models like the closed field line \cite{wang2018uniformly}\cite{narski2014asymptotic} and high order differential operators arising in plasma physics \cite{yang2019numerical}.

\subsection{Low Mach number limit of compressible flows}

Recently there has been increasing research activities in developing Mach 
number uniform fluid solvers.  Consider the case of isentropic Navier-Stokes
equations:
\begin{eqnarray}\label{low-Mach}
&&\partial_t \rho + \nabla_x\cdot(\rho u) = 0\,,\\
&& \partial_t(\rho u) + \nabla_x\cdot(\rho u\otimes u) +\frac{1}{\varepsilon^2}\nabla_xp = \frac{1}{\text{Re}}\Delta u\,.
\end{eqnarray}
Here $\rho$ is density, $u$ velocity, $p = \rho^\gamma$ pressue,
 $\varepsilon$ being the Mach number and ${\text {Re}}$ the Reynolds number.
When $\varepsilon\ll1$, one seeks the asymptotic expansions: $\rho = \rho^{(0)}+\varepsilon^2 \rho^{(2)}+\cdots$ and $p = p^{(0)}+\varepsilon^2 p^{(2)}+\cdots$
which then yieds \cite{MK81}:
\begin{eqnarray}
&&\nabla\cdot u^{(0)}=0\,,\\
&&\partial_t u^{(0)} + \left(u^{(0)}\cdot\nabla\right)u^{(0)} + \frac{1}{\rho^{(0)}}\nabla p^{(2)} = \frac{1}{\text{Re}}\Delta u^{(0)}\,.
\end{eqnarray}

The characteristic speeds of system (\ref{low-Mach}) are of $O(1/\e)$, corresponding to fast acoustic waves.
One would think the low Mach number problem is mainly a numerically stiff problem hence
 a small time step of $O(\e)$ is needed  if an explicit method is used. In fact,  the constraints are more severe.  
For  shock capturing methods, the numerical viscosity, a necessary ingredient to suppress artificial oscillations, is inversely proportional to the speed of sound hence $\Delta x = \mathcal{O}(\varepsilon)$  is needed to reduce numerical dissipation
\cite{GV99}\cite{Dell}. One then needs $\Delta t = \mathcal{O}(\varepsilon\Delta x)$ for numerical stability in an explicit scheme.

In developing a numerical scheme that is efficient for all Mach numbers,
ideally one hopes to  use mesh size and
time step {\it independent} of $\varepsilon$, namely the scheme is AP. This is usually achieved by
splitting the flux into the fast moving (corresponding to the acoustic waves) one and  the slowly moving one. 
An earlier attempt in this direction is a splitting by Klein \cite{Klein1995}, which was further improved in  \cite{noelle2014weakly}. Here we mention an 
approach introduced in \cite{HJL12} (see a related
approach in \cite {DegTang}\cite{CDK12}) which takes the following splitting:
\begin{equation}
\begin{cases}
\partial_t\rho +\alpha\nabla\cdot (\rho u) + (1-\alpha)\nabla\cdot(\rho u)=0\\
\partial_t(\rho u) + \nabla\cdot(\rho u\otimes u)+\nabla\left(\frac{p(\rho) - a(t)\rho}{\varepsilon^2}\right) +\frac{a(t)}{\varepsilon^2}\nabla\rho = \frac{1}{\text {Re}}\Delta u
\end{cases}\,,
\end{equation}
where $\alpha$ and $a(t)$ are artificial parameters.
By choosing $a(t)$ well approximating $p'(\rho)$, the third term in the second equation is non-stiff, thus will be treated explicitly. An implicit treatment on term $\nabla\rho$ is necessary, but thanks to its linearity this can  be done easily,
since only Poisson solvers are needed. The scheme is shock-capturing in
the high Mach number regime, and reduces to a projection method--a popular method for incompressible Navier-Stokes equations 
\cite{Temam1969} \cite{chorin1968} when $\varepsilon
\to 0$, hence the AP property for $\e \to 0$ is justified.

This  direction  is  still rapidly evolving. One can find other
techniques such as a  Langrange-Projection scheme \cite{chalons2016}\cite{Zak17}, a modification of the Roe solver \cite{Roepko}\cite{Klingen17} with applications to
astrophysics problems, careful choice of numerical viscosity \cite{DLV2017}, and error estimates on AP schemes for low-Mach number flows
\cite{feireisl2018}.

\bigskip
\centerline{\bf Acknowledgements}
\medskip

Many of the presented works were partially supported by NSF Divison of Mathematical Sciences Focused Research Group: "Collaborative
Research on Kinetic Description of Multiscale Phenomena: Modeling, Theory
and Computation" (2008-2011),  and NSF Division of Mathematical Sciences Research Network "KI-Net: Kinetic description of emerging
challenges in multiscale problems of natural sciences" (2012-2019).  The author thanks Eitan Tadmor for his tireless effort in making these projects successful,  and for his long-time support. 

The author also thanks Francis Filbet for contributing to subsection 6.3, and Min Tang for contributing subsection 6.4.

%\section{Conclusions}

%%%%%%%%%%%%%%%%%%%%%%%%%%%%%%%%%%%%%%%%%%%%%%%%%%%%%%%%%%%%%%%%%%%%%%%%%

%%%%%%%%%%%%%%%%%%%%%%%%%%%%%%%%%%%%%%%%%%%%%%%%%%%%%%%%%%%%%%%%%%%%%%%%

\bibliography{biblio}

%%%%%%%%%%%%%%%%%%%%%%%%%%%%%%%%%%%%%%%%%%%%%%%%%%%%%%%%%%%%%%%%%%%%%%%%

\label{lastpage}

\end{document}